# UNIVERSITY OF CRETE

## DEPARTMENT OF APPLIED MATHEMATICS

Doctoral Dissertation:

**HARDY INEQUALITIES IN GENERAL DOMAINS**

by Georgios Psaradakis

Supervisor: Professor Stathis Filippas

Heraklion
13 December 2011

*to Rena*

*Kyriakos*

*Maria*

*Froso*

*Lefteris*

# Committees

## Thesis committee

- Stathis Filippas (supervisor), Professor (Dept. of Appl. Math., Univ. of Crete)
- Georgia Karali, Assistant Professor (Dept. of Appl. Math., Univ. of Crete)
- Achilles Tertikas, Professor (Dept. of Math., Univ. of Crete)

## Thesis defence committee

- Gerassimos Barbatis, Assistant Professor (Dept. of Math., National and Kapodistrian Univ. of Athens)
- Stathis Filippas, Professor (Dept. of Appl. Math., Univ. of Crete)
- Spyridon Kamvissis, Professor (Dept. of Appl. Math., Univ. of Crete)
- Georgia Karali, Assistant Professor (Dept. of Appl. Math., Univ. of Crete)
- Alkis Tersenov, Professor (Dept. of Math., Univ. of Crete)
- Achilles Tertikas, Professor (Dept. of Math., Univ. of Crete)
- Athanasios Tzavaras, Professor (Dept. of Appl. Math., Univ. of Crete)

# Preface

This work is about improvements to multidimensional Hardy inequalities. We focus in two main directions:

(i) to obtain sharp homogeneous remainders to $L^1$ weighted Hardy inequalities, and

(ii) to obtain optimal Sobolev-type remainder terms to $L^p$ Hardy inequalities for $p > n$.

More precisely, regarding (i) we obtain homogeneous remainder terms to $L^1$ weighted Hardy inequalities (abbreviated HI) involving distance to the boundary of subsets of $\mathbb{R}^n$ with *finite inner radius*. The improvements obtained are in some sense optimal.

For a domain $\Omega$ of class $\mathcal{C}^2$ for which $\inf_{y \in \partial\Omega} \mathcal{H}(y) > 0$, where $\mathcal{H}(y)$ is the mean curvature of the boundary at $y \in \partial\Omega$, we obtain an optimal homogeneous improvement to the $L^1$-weighted HI with a lower estimate for the best constant of the remainder term. To this end we prove that the distributional Laplacian of the distance function to the boundary of a $\mathcal{C}^2$-smooth set is a signed Radon measure with nonnegative singular part, and absolutely continuous part satisfying $(-\Delta d)_{ac} \geq (n-1) \inf_{y \in \partial\Omega} \mathcal{H}(y)$ a.e. in $\Omega$. This leads to the fact that a $\mathcal{C}^2$-smooth domain is *mean convex* if and only if

$$-\Delta d \geq 0, \quad \text{in the sense of distributions in } \Omega. \tag{$\mathfrak{C}$}$$

Note that condition ($\mathfrak{C}$) has been shown in [BFT1] to be enough for the $L^p$ Hardy inequality to hold (without any smoothness condition on the boundary).

An upper bound involving the *total mean curvature* is obtained also for the best constant of the remainder term in the $L^1$-weighted HI mentioned above. These upper and lower bounds coincide when the domain is a ball, and in this particular case we achieve a finite series of optimal remainder terms with best possible constants. This is in contrast to the results for the $L^p$ case with $p > 1$, where an infinite series involving optimal logarithmic terms can be added.

Regarding (ii), we improve the sharp HI involving distance to the origin in case $p > n$, by adding an optimal weighted Hölder seminorm. To achieve this we first obtain a local improvement. Partial results are obtained when the distance is taken from the boundary. We also obtain a refinement of both the Sobolev inequality for $p > n$ and the HI involving either distance to the origin, or distance to the boundary, where the HI appears with best constant.

# Acknowledgements

This thesis is the result of four years of work during which I have been encouraged, supported and guided by my teacher Stathis Filippas. I would like to thank him once more, and also to express my sincere gratitude to him for providing knowledge and disposal of his time trying to make me take part in mathematical research and to feel the joy of a mathematical proof.

I would like to thank Achilles Tertikas for his encouragement and advise, and for his lectures in the graduate class "Partial Differential Equations - Weak Solutions" which introduced me to Sobolev inequalities.

Special thanks goes also to Ioannis Athanasopoulos for his seminar lectures on issues of Partial Differential Equations, to Themis Mitsis for five semesters of collaboration in my teaching assistance duty and to my friends K. Gkikas and C. Saroglou for our many discussions. I wish to thank my teachers and other members of the Departments of Mathematics and Applied Mathematics of the University of Crete for the nice atmosphere.

I was funded for a year by the Maria M. Manassaki bequest scholarship and for three years by the State Scholarship Foundation of Greece (I.K.Y.).

# Notation

Throughout this thesis:
- The letter $n$ will always denote a positive integer.
- $\mathbb{R}^n$ stands for the $n$-dimensional Euclidean space.
- $|x|$ is the Euclidean length of $x \in \mathbb{R}^n$.
- $\overline{\Omega}$ is the closure of the set $\Omega \subset \mathbb{R}^n$.
- If $\Omega \subsetneq \mathbb{R}^n$ is open, then by $\partial\Omega$ we denote the boundary of $\Omega$, i.e. $\partial\Omega = \overline{\Omega} \setminus \Omega$.
- The characteristic function of a set $\Omega \subset \mathbb{R}^n$ is denoted by $\chi_\Omega$.
- $\nabla u$ is the Gradient and $\Delta u$ is the Laplacian of a real valued smooth function $u$.
- $\operatorname{div} \vec{T}$ is the Divergence of a smooth vector field $\vec{T}$.
- $\operatorname{sprt}\{u\}$ is the support of a continuous function $u$.
- A *domain* is an open and connected subset of $\mathbb{R}^n$.
- The Lebesgue measure on $\mathbb{R}^n$ is denoted by $\mathcal{L}^n$.

- By $B_r(x)$ (resp. $\overline{B}_r(x)$) we denote an open (resp. closed) ball ("disk" if $n = 2$) having radius $r > 0$ and center at $x \in \mathbb{R}^n$. When the center is of no importance we simply write $B_r$.

- The volume of $B_1$ is denoted by $\omega_n$, i.e. $\omega_n := \mathcal{L}^n(B_1) = \pi^{n/2}/\Gamma(1 + n/2)$, where $\Gamma$ is the Gamma Function defined by $\Gamma(x) := \int_0^\infty t^{x-1} e^{-x} dt; x \in \mathbb{R}$. Thus the surface area of $\partial B_1$ is $n\omega_n$.

- If $\Omega \subseteq \mathbb{R}^n$ is open then by $C^\infty(\Omega)$ we denote the set of functions defined in $\Omega$ which are infinite times differentiable with respect to any variable.

- By $C_c^\infty(\Omega)$ we denote the subset of $C^\infty(\Omega)$ comprised of functions having compact support in $\Omega$. We will always extend every function in $C_c^\infty(\Omega)$ to be zero outside it's support.

- If $\Omega \subseteq \mathbb{R}^n$ is open and $p \in [1, \infty]$ then $L^p(\Omega)$ is the set of all measurable functions $u : \Omega \to \mathbb{R}$, such that

$$\|u\|_{L^p(\Omega)} := \begin{cases} \left( \int_\Omega |u|^p dx \right)^{1/p}, & \text{if } p \geq 1 \\ \operatorname{ess\,sup}_\Omega |u|, & \text{if } p = \infty. \end{cases}$$

- $L^p_{loc}(\Omega)$ is the set of all measurable functions $u : \Omega \to \mathbb{R}$, such that $u \in L^p(V)$ for each compact set $V \subset \Omega$.

- If $\Omega \subseteq \mathbb{R}^n$ is open and $p \in [1, \infty)$, then we say that a function $u \in L^1_{loc}(\Omega)$ belongs to the *Sobolev space* $W^{1,p}(\Omega)$ if $u \in L^p(\Omega)$ and the weak first partial derivatives with respect to any variable exist and belong to $L^p(\Omega)$ as well. It is a Banach space under the norm

$$\|u\|_{W^{1,p}(\Omega)} := \left( \|u\|^p_{L^p(\Omega)} + \|| \nabla u |\|^p_{L^p(\Omega)} \right)^{1/p}.$$

- If $\Omega \subsetneq \mathbb{R}^n$ is open and $p \in [1, \infty)$, then we denote by $W^{1,p}_0(\Omega)$ the closure of $C^\infty_c(\Omega)$ under the norm of $W^{1,p}(\Omega)$.

- Let $\Omega$ be an open subset of $\mathbb{R}^n$ and $u \in L^1(\Omega)$. We say $u$ is of *bounded variation* in $\Omega$ if it has finite *total variation* in $\Omega$, where

$$\text{total variation of } u \text{ in } \Omega := \sup \left\{ \int_\Omega u \operatorname{div} \phi \, \mathrm{d}x; \ \phi \in C^\infty_c(\Omega; \mathbb{R}^n) \text{ with } \|| \phi |\|_{L^\infty(\Omega)} \leq 1 \right\}.$$

We denote by $BV(\Omega)$ the set of all functions having bounded variation in $\Omega$. It is a Banach space under the norm

$$\|u\|_{BV(\Omega)} := \|u\|_{L^1(\Omega)} + \text{total variation of } u \text{ in } \Omega.$$

- The symbol ■ means *end of proof* while the symbol □ denotes *end of example*, or *end of remark*.

# Contents

# Chapter 1

# Introduction

Scaling invariant inequalities involving integrals (possibly weighted) of functions and their derivatives in various powers, constitute a basic tool in analysis, in the theory of partial differential equations and in the calculus of variations, see [Br], [Gq], [HLP], [LbL], [LSU], [Mz2], [Mrr], [N] and [St]. In addition they find various applications in several branches of geometry and physics; see for instance [Lb], [S-C], [Zh]. Some essays on this kind of inequalities and the spaces that they indicate are [AdF], [Ln], [Mz1], [OpK].

In the three last decades a lot of attention has been given to improved versions of some of the aforementioned type inequalities. By *improved* it is usually meant that one considers the initial inequality and adds a positive term in the least hand side. Of course this is not always possible. Nevertheless, when possible, resulted inequalities play an essential role in the theory of partial differential equations and nonlinear analysis. They are used for instance in the study of the stability of solutions of elliptic and parabolic equations (see for example [BrV], [FT], [DN], [Gk1], [PV]), in the study of existence and asymptotic behavior of solutions of heat equations with singular potentials; see for instance [BrN], [CM], [DD], [GP], [VZ], as well as in the study of the stability of eigenvalues in elliptic problems (see [D3], [FlHTh]).

## 1.1 Sobolev inequalities with remainder terms

Two important inequalities of the kind described above were proved by S. L. Sobolev:

(I) Let $1 < p < n;\ n \geq 2$. There exists a constant $S(n,p) > 0$ such that for all $u \in C_c^\infty(\mathbb{R}^n)$

$$\left(\int_{\mathbb{R}^n} |\nabla u|^p \mathrm{d}x\right)^{1/p} \geq S(n,p) \left(\int_{\mathbb{R}^n} |u|^{np/(n-p)} \mathrm{d}x\right)^{(n-p)/np}. \tag{1.1}$$

This inequality is optimal in the sense that it fails if $np/(n-p)$ is replaced by any number $q > np/(n-p)$. The best constant

$$S(n,p) = n^{1/p}\sqrt{\pi}\left[\frac{\Gamma(n/p)\Gamma(1+n-n/p)}{\Gamma(1+n/2)\Gamma(n)}\right]^{1/n}\left(\frac{n-p}{p-1}\right)^{1-1/p},$$

and the family of extremal functions

$$U_{\alpha,\beta}(x) = \alpha[1 + \beta|x - x_0|^{1-1/p}]^{1-n/p}; \quad \alpha \neq 0, \beta > 0, x_0 \in \mathbb{R}^n, \tag{1.2}$$

for $x \in \mathbb{R}^n$, have been found simultaneously in [A] and [Tl1].

(II) If $p > n \geq 1$ and $\Omega$ is an open set in $\mathbb{R}^n$ with finite volume $\mathcal{L}^n(\Omega)$, then there exists a constant $s(n,p) > 0$ such that for all $u \in C_c^\infty(\Omega)$

$$\sup_{x \in \Omega}|u(x)| \leq s(n,p)[\mathcal{L}^n(\Omega)]^{1/n - 1/p}\left(\int_\Omega |\nabla u|^p \mathrm{d}x\right)^{1/p}. \tag{1.3}$$

This inequality is also optimal in the sense that it fails if $1/n - 1/p$ is replaced by any exponent $q < 1/n - 1/p$ on $\mathcal{L}^n(\Omega)$. The best constant

$$s(n,p) = n^{-1/p}\omega_n^{-1/n}\left(\frac{p-1}{p-n}\right)^{1-1/p},$$

and the family of extremal functions

$$V_{\alpha,\beta}(x) = \alpha[\beta^{(p-n)/(p-1)} - |x - x_0|^{(p-n)/(p-1)}]; \quad \alpha \neq 0, \beta > 0, x_0 \in \mathbb{R}^n, \tag{1.4}$$

for $x \in B_\beta(x_0)$ and $V_{\alpha,\beta}(x) = 0$ otherwise, have been found in [Tl2].

We can easily verify that the family of functions $U_{\alpha,\beta}$ defined in (1.2) belong to the space $W^{1,p}(\mathbb{R}^n)$ and thus the constant $S(n,p)$ in (1.1) is achieved. If we replace $\mathbb{R}^n$ by an open set $\Omega \subsetneq \mathbb{R}^n$ then (1.1) still holds for any $u \in C_c^\infty(\Omega)$. In this case it is well known that: *the constant remains optimal but it is no more attained by some function in the corresponding Sobolev space $W_0^{1,p}(\Omega)$*. This fact was interpreted in a quantitative form by H. Brezis and L. Nirenberg in their highly cited paper [BrN], where they estimate from below the difference of the two sides of (1.1) in the case $p = 2$. Their result can be stated as follows:

BREZIS, NIRENBERG (1983) - *Assume $\Omega$ is a domain in $\mathbb{R}^n$; $n \geq 3$, with $\mathcal{L}^n(\Omega) < \infty$ and let $1 < q < n/(n-2)$. There exists a positive constant $C(n,q)$ such that for all $u \in C_c^\infty(\Omega)$*

$$\left[\int_\Omega |\nabla u|^2 \mathrm{d}x - S^2(n,2)\left(\int_\Omega |u|^{2n/(n-2)}\mathrm{d}x\right)^{(n-2)/n}\right]^{1/2}$$
$$\geq \frac{C(n,q)}{[\mathcal{L}^n(\Omega)]^{1/q - (n-2)/2n}}\left(\int_\Omega |u|^q \mathrm{d}x\right)^{1/q}. \tag{1.5}$$

*Moreover,* (i) *the above estimate fails for the critical case* $q = n/(n-2)$, (ii) *the sharp constant* $C(3,2)$ *in case* $\Omega$ *is any ball in* $\mathbb{R}^3$ *and* $q = 2$, *is* $\sqrt[3]{\pi^4/6}$ *and is not achieved in* $W_0^{1,2}(\Omega)$.

This surprising result shows the existence of *remainder terms* in Sobolev's inequality (1.1). The direct generalization of (1.5) for $p \in (1, n)$ was given in [EPTr]. In [BrL] the authors substitute the $L^q$ norm on right hand side of (1.5) by the weak-$L^{n/(n-2)}$ norm, and also (even more strongly) by the weak-$L^{n/(n-1)}$ norm of the length of the gradient (see also [Alv] for some weak-type remainder terms). In addition, a direction initiated in [BrL] for $p = 2$ is to consider functions not necessarily vanishing on $\partial\Omega$, so that a trace remainder term appears. The full picture for other values of $p \in (1, n)$, was completed in [MV].

Another open question in [BrL] (for $p = 2$) asked if one can bound in some natural way the difference of the two sides of (1.1) in terms of the *distance from the set of the extremal functions* $U_{\alpha,\beta}$. In [BEgn], a quantity that measures the distance between $\nabla u$ and $\nabla U_{\alpha,\beta}$ was added on the right hand side of (1.1) for $p = 2$. Recently a result for all values of $p \in (1, n)$ was obtained in [CFMPr] where the *functional asymmetry* of $u$, a quantity that measures the distance between $u$ and the family $U_{\alpha,\beta}$, was added on the right hand side of (1.1). The analogous result for (1.3) had been given earlier in [C].

Finally, a different direction in strengthening Sobolev's inequality without remainder terms being involved was given in [LYZh]. In that work the $L^p$-norm of the length of the gradient of $u$ in (1.1) is replaced by a smaller quantity called *the affine energy* of $u$, which is additionally invariant under affine transformations of $\mathbb{R}^n$. The analogous result for (1.3) is in [CLYZh]. For a simplified approach we refer to [AlBB].

## 1.2 Hardy inequalities with remainder terms

Another well known family of scale invariant inequalities consists of Hardy inequalities. These involve the distance function usually taken from a point or from the boundary of a set; see for instance [HLP], [CKN], [D], [D1], [D2]), [Mz1], [N] and [OpK].

**Distance to the boundary.** We denote a generic point in $\mathbb{R}^n$ by $x = (x', x_n)$, where $x' = (x_1, ..., x_{n-1})$. Hardy's inequality in the half space $\mathbb{R}_+^n := \mathbb{R}^n \cap \{x_n > 0\}$; $n \geq 2$, asserts that if $p > 1$ then for all $u \in C_c^\infty(\mathbb{R}_+^n)$

$$\int_{\mathbb{R}_+^n} |\nabla u|^p \mathrm{d}x \geq \left(\frac{p-1}{p}\right)^p \int_{\mathbb{R}_+^n} \frac{|u|^p}{x_n^p} \mathrm{d}x. \tag{1.6}$$

It is well known that the constant appearing on the right hand-side is sharp and is not attained in $W_0^{1,2}(\mathbb{R}_+^n)$. V. G. Maz'ya in [Mz1]-§2.1.6 replaced the $L^2$-norm of the length of the gradient in (1.1) for $p = 2$ by the *sharp Hardy difference* on $\mathbb{R}_+^n$ (implied by (1.6)) for $p = 2$. More precisely

MAZ'YA (1985) - *There exists a constant $M(n) > 0$ depending only on $n$, such that for any $u \in C_c^\infty(\mathbb{R}_+^n)$*

$$\left( \int_{\mathbb{R}_+^n} |\nabla u|^2 \mathrm{d}x - \frac{1}{4} \int_{\mathbb{R}_+^n} \frac{|u|^2}{x_n^2} \mathrm{d}x \right)^{1/2} \geq M(n) \left( \int_{\mathbb{R}_+^n} |u|^{2n/(n-2)} \mathrm{d}x \right)^{(n-2)/2n}. \tag{1.7}$$

The existence of extremals for $n \geq 4$ has been established in [TT] where it is proved that the best constant $M(n)$ in (1.7) satisfies $M(n) < S(n,2)$ for all $n \geq 4$, where $S(n,2)$ is the best constant in Sobolev's inequality (1.1) for $p = 2$. The best constant for $n = 3$ is found independently in [BFrL] and [MncS], and surprisingly one has $M(3) = S(3,2)$.

A version of (1.7) in more general domains was given in [FMzT3]. Let us recall first that in [BFT1] the authors proved the following Hardy inequality involving distance to the boundary that generalizes (1.6),

$$\int_\Omega |\nabla u|^p \mathrm{d}x \geq \left( \frac{p-1}{p} \right)^p \int_\Omega \frac{|u|^p}{d^p} \mathrm{d}x, \tag{1.8}$$

for all $u \in C_c^\infty(\Omega)$. Here $d(x) := \mathrm{dist}(x, \mathbb{R}^n \setminus \Omega)$ and $\Omega \subsetneq \mathbb{R}^n$ is any domain satisfying

$$-\Delta d \geq 0 \quad \textit{in the sense of distributions in } \Omega. \tag{$\mathfrak{C}$}$$

They also proved that: *under condition ($\mathfrak{C}$), the constant on the right hand-side of (1.8) is optimal.* Note here that no smoothness assumption on the boundary is imposed. In [FMzT1] and [FMzT3] it is proved that if $\Omega$ is smooth enough and has finite inner radius, then one can replace the $L^p$-norm of the length of the gradient in (1.1) for $2 \leq p < n$ by the sharp Hardy difference on $\Omega$. Their precise statement reads as follows: *Suppose condition ($\mathfrak{C}$) is valid for a domain $\Omega$ having finite inner radius and boundary of class $\mathcal{C}^2$. Then letting $2 \leq p < n$, there exists a positive constant $C(n, p, \Omega)$ depending on $n, p$ and $\Omega$, such that for all $u \in C_c^\infty(\Omega)$*

$$\left( \int_\Omega |\nabla u|^p \mathrm{d}x - \left( \frac{p-1}{p} \right)^p \int_\Omega \frac{|u|^p}{d^p} \mathrm{d}x \right)^{1/p} \geq C(n,p,\Omega) \left( \int_\Omega |u|^{np/(n-p)} \mathrm{d}x \right)^{(n-p)/np}. \tag{1.9}$$

For an analogous condition and results for domains having infinite inner radius see [Gk2]. Recently, in [FrL] the dependence of the constant $C(n,p,\Omega)$ on $\Omega$ was dropped when the domain is convex. Let us note here that condition ($\mathfrak{C}$) is more general than convexity for $n \geq 3$, and equivalent to convexity for $n = 2$.

Instead of having a Sobolev norm as a remainder term, one may look for different kind of improvements. In this direction H. Brezis and M. Marcus have considered *homogeneous* remainder terms:

BREZIS, MARCUS (1997) - *Assume $\Omega$ is a bounded convex domain in $\mathbb{R}^n$; $n \geq 1$, and set $D := \operatorname{diam}(\Omega)$. Then for all $u \in C_c^\infty(\Omega)$*

$$\int_\Omega |\nabla u|^2 \mathrm{d}x - \frac{1}{4} \int_\Omega \frac{|u|^2}{d^2} \mathrm{d}x \geq \frac{1}{4} \int_\Omega \frac{|u|^2}{d^2} X^2(d/D) \mathrm{d}x, \tag{1.10}$$

*where $X(t) = (1-\log t)^{-1}$, $t \in (0,1]$. Moreover, the weight function $X^2$ is optimal in the sense that the power $2$ cannot be decreased, and the coefficient $1/4$ on the right hand side is sharp.*

This result has been generalized in many aspects in [BFT1], [BFT2]: all values of $p \in (1,\infty)$ were considered, domains merely satisfying condition ($\mathfrak{C}$) and having finite inner radius were allowed, and an infinite series involving iterated logarithmic potentials (in some sense optimal) was added.

**Distance to a point.** Another well-known version of Hardy's inequality asserts that if $n,p \geq 1$ with $p \neq n$, then for all $u \in C_c^\infty(\mathbb{R}^n \setminus \{0\})$

$$\int_{\mathbb{R}^n} |\nabla u|^p \mathrm{d}x \geq \left|\frac{n-p}{p}\right|^p \int_{\mathbb{R}^n} \frac{|u|^p}{|x|^p} \mathrm{d}x, \tag{1.11}$$

where the constant $|(n-p)/p|^p$ is sharp. In contrast to (1.1), the sharp constant is not attained by some function in the corresponding Sobolev space (which is $W^{1,p}(\mathbb{R}^n)$ when $1 \leq p < n$ and $W^{1,p}(\mathbb{R}^n \setminus \{0\})$ when $p > n$). Motivated by [BrN] and [BrL], it is natural to ask if one can have remainder terms on the right hand side of (1.11). More precisely, one might ask wether an inequality of the form

$$\int_{\mathbb{R}^n} |\nabla u|^p \mathrm{d}x \geq \left|\frac{n-p}{p}\right|^p \int_{\mathbb{R}^n} \frac{|u|^p}{|x|^p} \mathrm{d}x + C \left(\int_{\mathbb{R}^n} |u|^q V(|x|) \mathrm{d}x\right)^{p/q},$$

holds for some nontrivial potential function $V \geq 0$, some $q > 0$ and some positive constant $C = C(n,p,q)$. The authors in [dPDFT] explain that this is not true by testing the above inequality with the function

$$U_\varepsilon(x) = \begin{cases} |x|^{(p-n)/p+\varepsilon}, & |x| \leq 1 \\ |x|^{(p-n)/p-\varepsilon}, & |x| > 1, \end{cases}$$

and taking the limit $\varepsilon \downarrow 0$.

If we replace $\mathbb{R}^n$ by a bounded open set $\Omega \subset \mathbb{R}^n$ containing $0$, then (1.11) continues to hold for any $u \in C_c^\infty(\Omega \setminus \{0\})$ and the constant $|(n-p)/p|^p$ remains optimal. It turns out that in this case one can have remainder terms.

In their pioneering work [BrV], H. Brezis and J. L. Vazquez improved Hardy's inequality (1.11) in the case $p = 2$ as follows

BREZIS, VAZQUEZ (1997) - *Let $\Omega$ be a domain in $\mathbb{R}^n$; $n \geq 3$, with $\mathcal{L}^n(\Omega) < \infty$ and $0 \in \Omega$. If $1 < q < 2n/(n-2)$, then there exists a constant $C(n, q) > 0$ such that for all $u \in C_c^\infty(\Omega)$*

$$\left(\int_\Omega |\nabla u|^2 \mathrm{d}x - \left(\frac{n-2}{2}\right)^2 \int_\Omega \frac{|u|^2}{|x|^2} \mathrm{d}x\right)^{1/2} \geq \frac{C(n,q)}{[\mathcal{L}^n(\Omega)]^{1/q - (n-2)/2n}} \left(\int_\Omega |u|^q \mathrm{d}x\right)^{1/q}. \quad (1.12)$$

*Moreover, the constant $C(n, q)$ equals the first eigenvalue of the Dirichlet Laplacian for the unit disk in $\mathbb{R}^2$ and it is optimal when $\Omega$ is a ball in $\mathbb{R}^n$ centered at the origin and $q = 2$, independently of the dimension $n \geq 2$.*

This result motivated many mathematicians to look for remainder terms to Hardy's inequality (1.11); [AbdCP], [AdChR], [AlvVV] [BFT1], [BFT2], [CF], [CP], [GGrM], [RSW], [VZ] and [WW]. In Problem 2 of [BrV], the question of whether there is a further improvement in the direction of the inequality (1.12) is posed. An optimal answer was given in [FT], where it was shown that the critical exponent $q = 2n/(n-2)$ is possible after considering a logarithmic correction weight for which the sharp exponent was given. More precisely it is proved that: *if $\Omega$ is a bounded domain in $\mathbb{R}^n$; $n \geq 3$, containing the origin, then there exists a positive constant $C(n)$ such that for all $u \in C_c^\infty(\Omega)$*

$$\left(\int_\Omega |\nabla u|^2 \mathrm{d}x - \left(\frac{n-2}{2}\right)^2 \int_\Omega \frac{|u|^2}{|x|^2} \mathrm{d}x\right)^{1/2}$$
$$\geq C(n) \left(\int_\Omega |u|^{2n/(n-2)} X^{1+n/(n-2)}(|x|/D) \mathrm{d}x\right)^{(n-2)/2n}, \quad (1.13)$$

*where $D = \sup_{x \in \Omega} |x|$ and $X(t) = (1 - \log t)^{-1}$, $t \in (0, 1]$. Moreover, the weight function $X^{1+n/(n-2)}$ is optimal in the sense that the power $1 + n/(n-2)$ cannot be decreased*; see [AdFT] for a second proof where in addition the best constant $C(n)$ is obtained.

## 1.3 Main results

The contribution of this thesis to the study of remainder terms in Hardy inequalities splits in two parts.

### 1.3.1 Part I: $L^1$ Hardy inequalities with weights

Recall that Hardy's inequality involving distance from the boundary of a convex set $\Omega \subsetneq \mathbb{R}^n$; $n \geq 1$, asserts that

$$\int_\Omega |\nabla u|^p \mathrm{d}x \geq \left(\frac{p-1}{p}\right)^p \int_\Omega \frac{|u|^p}{d^p} \mathrm{d}x, \quad p > 1, \quad (1.14)$$

for all $u \in C_c^\infty(\Omega)$, where $d \equiv d(x) := \mathrm{dist}(x, \mathbb{R}^n \setminus \Omega)$. Due to [HLP], [MtskS] and [MMP] the constant appearing in (1.14) is optimal. After the pioneering results in [Mz1] (see inequality (1.7)) and [BrM] (see inequality (1.10)), a sequence of papers have improved (1.14) by adding extra terms on its right hand side, see for instance [H-OL], [Tdb1], [Tdb2], [EH-S], [BFT2], [BFT3], [FMzT3], [FTT], [FrL], [BFrL] and primarily [BFT1] and [FMzT1], [FMzT2] where it was also noted that (1.14) remains valid with the sharp constant in more general sets than convex ones, and in particular in sets that satisfy $-\Delta d \geq 0$ in the distributional sense (condition ($\mathfrak{C}$)).

In the case $p = 1$, (1.14) reduces to a trivial inequality, at least for sets having non positive distributional Laplacian of the distance function. However, in the one dimensional case, the following $L^1$ weighted Hardy inequality is well known:

$$\int_0^\infty \frac{|u'(x)|}{x^{s-1}} \mathrm{d}x \geq (s-1) \int_0^\infty \frac{|u(x)|}{x^s} \mathrm{d}x; \quad s > 1, \tag{1.15}$$

for all absolutely continuous functions $u : [0, \infty) \to \mathbb{R}$, such that $u(0) = 0$. This is the special case $p = 1$ of Theorem 330 in the classical treatise of G. Hardy, J. E. Littlewood and G. Pólya, [HLP]. Inequality (1.15) becomes equality for $u$ increasing and thus the constant on the right hand side is sharp.

Assume now that $\Omega$ is a domain satisfying condition ($\mathfrak{C}$) :

$$-\Delta d \geq 0 \quad \textit{in the sense of distributions in } \Omega. \tag{$\mathfrak{C}$}$$

Without much effort (see Remark 2.7 in the present thesis) one can see that both (1.14) and (1.15) can be generalized as follows

$$\int_\Omega \frac{|\nabla u|^p}{d^{s-p}} \mathrm{d}x \geq \left(\frac{s-1}{p}\right)^p \int_\Omega \frac{|v|^p}{d^s} \mathrm{d}x, \tag{1.16}$$

valid for any $s > 1, p \geq 1$ and any $u \in C_c^\infty(\Omega)$. Further, we may follow [BFT1] to prove that if $\Omega$ has in addition finite inner radius, then for all $s > 1, p > 1$ and any $u \in C_c^\infty(\Omega)$ we have

$$\int_\Omega \frac{|\nabla u|^p}{d^{s-p}} \mathrm{d}x - \left(\frac{s-1}{p}\right)^p \int_\Omega \frac{|u|^p}{d^s} \mathrm{d}x \geq \frac{1}{2} \frac{p-1}{p} \left(\frac{s-1}{p}\right)^{p-2} \int_\Omega \frac{|u|^p}{d^s} X^2(d/D) \mathrm{d}x, \tag{1.17}$$

where $D = B(n, p, s) \sup_{x \in \Omega} d(x)$ and $X(t) = (1 - \log t)^{-1}; t \in (0, 1]$. The weight function $X^2$ is optimal in the sense that the power 2 cannot be decreased, and the constant on the right-hand side is the best possible (see §2.2 of the present thesis). We point out that as $p \downarrow 1$ the right hand side vanishes. This motivated us to search for remainder terms in the limit case $p = 1$, that is to search for an inequality of the type

$$\int_\Omega \frac{|\nabla u|}{d^{s-1}} \mathrm{d}x \geq (s-1) \int_\Omega \frac{|u|}{d^s} \mathrm{d}x + \mathcal{B}_1 \int_\Omega V(d)|u| \mathrm{d}x; \quad s \geq 1, \tag{1.18}$$

valid for all $u \in C_c^\infty(\Omega)$. Here $\mathcal{B}_1 \in \mathbb{R}$ and $V$ is a potential function, i.e. nonnegative and $V \in L^1_{loc}(\mathbb{R}^+)$. Questions concerning optimal inverse distance power potentials, sharp constants for the remainder term and possible further improvements will be studied. In this direction our first result states

**Theorem I** (i) *Let $\Omega$ be a domain in $\mathbb{R}^n$ with boundary of class $\mathcal{C}^2$ satisfying a uniform interior sphere condition and we denote by $\underline{\mathcal{H}} := \inf_{y \in \partial\Omega} \mathcal{H}(y)$ the infimum of the mean curvature of the boundary. Then there exists $\mathcal{B}_1 \geq (n-1)\underline{\mathcal{H}}$ such that for all $u \in C_c^\infty(\Omega)$ and all $s \geq 1$*

$$\int_\Omega \frac{|\nabla u|}{d^{s-1}} \mathrm{d}x \geq (s-1) \int_\Omega \frac{|u|}{d^s} \mathrm{d}x + \mathcal{B}_1 \int_\Omega \frac{|u|}{d^{s-1}} \mathrm{d}x. \qquad (1.19)$$

(ii) *Let $s \geq 2$. If $\Omega$ is a bounded domain in $\mathbb{R}^n$ with boundary of class $\mathcal{C}^2$ having strictly positive mean curvature, then the constant $s-1$ in the first term as well as the exponent $s-1$ on the distance function on the remainder term in (1.19), are optimal. In addition, we have the following estimates*

$$(n-1)\underline{\mathcal{H}} \leq \mathcal{B}_1 \leq \frac{n-1}{|\partial\Omega|} \int_{\partial\Omega} \mathcal{H}(y) \mathrm{d}S_y, \qquad (1.20)$$

*where $\mathcal{H}(y)$ is the mean curvature of the boundary at $y \in \partial\Omega$, and $\underline{\mathcal{H}} := \min_{y \in \partial\Omega} \mathcal{H}(y)$ is its minimum value.*

We stress that in part (i) no other condition than smoothness of the boundary is imposed on $\Omega$. Thus the second term on the right hand side of (1.19) is a remainder term in case $\underline{\mathcal{H}} > 0$ only.

The following result, which is of independent interest, played a key role in establishing Theorem I and is proved in §3 of this thesis:

**Theorem II** *Let $\Omega \subset \mathbb{R}^n$ be a domain with boundary of class $\mathcal{C}^2$ satisfying a uniform interior sphere condition. Then $\mu := (-\Delta d)\mathrm{d}x$ is a signed Radon measure on $\Omega$. Let $\mu = \mu_{ac} + \mu_s$ be the Lebesgue decomposition of $\mu$ with respect to $\mathcal{L}^n$, i.e. $\mu_{ac} \ll \mathcal{L}^n$ and $\mu_s \perp \mathcal{L}^n$. Then $\mu_s \geq 0$ in $\Omega$, and $\mu_{ac} \geq (n-1)\underline{\mathcal{H}}\mathrm{d}x$ a.e. in $\Omega$, where $\underline{\mathcal{H}} := \inf_{y \in \partial\Omega} \mathcal{H}(y)$.*

An easy consequence of Theorem II is

**Corollary** *Let $\Omega$ be a domain with boundary of class $\mathcal{C}^2$ satisfying a uniform interior sphere condition. Then $\Omega$ is mean convex, i.e. $\mathcal{H}(y) \geq 0$ for all $y \in \partial\Omega$, if and only if $-\Delta d \geq 0$ holds in $\Omega$, in the sense of distributions.*

We emphasize that a set $\Omega \subsetneq \mathbb{R}^n$ with distance function having non positive distributional Laplacian, is shown in [BFT1], [BFT2], [BFT3] and [FMzT1], [FMzT2], [FMzT3] to be the natural assumption for the validity of various Hardy inequalities.

In special geometries we are able to compute the best constant $\mathcal{B}_1$ in (1.19):

In case $\Omega$ is a ball of radius $R$ then the upper and lower estimates (1.20) coincide, yielding $\mathcal{B}_1 = (n-1)/R$, i.e.

$$\int_{B_R} \frac{|\nabla u|}{d^{s-1}} \mathrm{d}x \geq (s-1) \int_{B_R} \frac{|u|}{d^s} \mathrm{d}x + \frac{n-1}{R} \int_{B_R} \frac{|u|}{d^{s-1}} \mathrm{d}x. \qquad (1.21)$$

One then may ask if (1.21) can be further improved. We provide a full answer to this question by showing that for $s \geq 2$ one can add a finite series of $[s] - 1$ terms on the right hand side before adding an optimal logarithmic correction. More precisely we prove the following

**Theorem III** *Let $B_R$ be a ball of radius $R$, then: $(i)$ For all $u \in C_c^\infty(B_R)$, all $s \geq 2, \gamma > 1$, there holds*

$$\int_{B_R} \frac{|\nabla u|}{d^{s-1}} \mathrm{d}x \geq (s-1) \int_{B_R} \frac{|u|}{d^s} \mathrm{d}x + \sum_{k=1}^{[s]-1} \frac{n-1}{R^k} \int_{B_R} \frac{|u|}{d^{s-k}} \mathrm{d}x + \frac{C}{R^{s-1}} \int_{B_R} \frac{|u|}{d} X^\gamma \left(\frac{d}{R}\right) \mathrm{d}x, \quad (1.22)$$

*where $X(t) := (1 - \log t)^{-1}$, $t \in (0, 1]$ and $C \geq \gamma - 1$. The exponents $s$ and $s - k$; $k = 1, 2, ..., [s] - 1$, on the distance function, as well as the constants $s - 1$, $(n - 1)/R^k$; $k = 1, 2, ..., [s] - 1$, in the first and the summation terms respectively, are optimal. The last term in (1.22) is optimal in the sense that if $\gamma = 1$, there is not positive constant $C$ such that (1.22) holds.*

*$(ii)$ For all $u \in C_c^\infty(B_R)$, all $1 \leq s < 2, \gamma > 1$, there holds*

$$\int_{B_R} \frac{|\nabla u|}{d^{s-1}} \mathrm{d}x \geq (s-1) \int_{B_R} \frac{|u|}{d^s} \mathrm{d}x + \frac{C}{R^{s-1}} \int_{B_R} \frac{|u|}{d} X^\gamma \left(\frac{d}{R}\right) \mathrm{d}x, \quad (1.23)$$

*where $X(t) := (1 - \log t)^{-1}$, $t \in (0, 1]$ and $C \geq \gamma - 1$. The last term in (1.23) is optimal in the sense that if $\gamma = 1$, there is not positive constant $C$ such that (1.23) holds.*

**Remark.** Let us make clear in what sense the added terms on the right hand side of (1.22) are optimal. For any $s \geq 1$ set

$$I_0[u] := \int_{B_R} \frac{|\nabla u|}{d^{s-1}} \mathrm{d}x - (s-1) \int_{B_R} \frac{|u|}{d^s} \mathrm{d}x.$$

Then for any $s \geq 2$ we prove that

$$\inf_{u \in C_c^\infty(B_R) \setminus \{0\}} \frac{I_0[u]}{\int_{B_R} \frac{|u|}{d^\beta} \mathrm{d}x} = \begin{cases} (n-1)/R, & \text{if } \beta = s - 1 \\ 0, & \text{if } \beta > s - 1. \end{cases}$$

Thus,

$$I_0[u] \geq \frac{n-1}{R} \int_{B_R} \frac{|u|}{d^{s-1}} \mathrm{d}x, \quad (1.24)$$

and this is the optimal inverse distance power remainder term with best constant, to the $L^1$ weighted Hardy inequality $I_0[u] \geq 0$. Next, for any $s \geq 3$ we prove that

$$\inf_{u \in C_c^\infty(B_R) \setminus \{0\}} \frac{I_0[u] - \frac{n-1}{R} \int_{B_R} \frac{|u|}{d^{s-1}} \mathrm{d}x}{\int_{B_R} \frac{|u|}{d^\beta} \mathrm{d}x} = \begin{cases} (n-1)/R^2, & \text{if } \beta = s - 2 \\ 0, & \text{if } \beta > s - 2. \end{cases}$$

Thus,

$$I_0[u] - \frac{n-1}{R}\int_{B_R}\frac{|u|}{d^{s-1}}\mathrm{d}x \geq \frac{n-1}{R^2}\int_{B_R}\frac{|u|}{d^{s-2}}\mathrm{d}x,$$

and this is the optimal inverse distance power remainder term with best constant, to the $L^1$ weighted improved Hardy inequality (1.24). We proceed in the same fashion for $s \geq 4$, and for precisely $[s] - 1$ steps.

Note that this is in contrast with the results in case $p > 1$, where an infinite series involving optimal logarithmic terms can be added (see [BFT2]) and ([BFT3]).

In case $\Omega$ is an infinite strip, using a more general upper bound on $\mathcal{B}_1$ (see Theorem 4.22), we prove that $\mathcal{B}_1 = 0$. As a matter of fact the finite series structure of (1.22) disappears and only the final logarithmic correction term survives. More precisely

**Theorem IV** *Let $S_R$ be an infinite strip of inner radius $R$. For all $u \in C_c^\infty(S_R)$, all $s \geq 1$, $\gamma > 1$, there holds*

$$\int_{S_R}\frac{|\nabla u|}{d^{s-1}}\mathrm{d}x \geq (s-1)\int_{S_R}\frac{|u|}{d^s}\mathrm{d}x + \frac{C}{R^{s-1}}\int_{S_R}\frac{|u|}{d}X^\gamma\Big(\frac{d}{R}\Big)\mathrm{d}x, \tag{1.25}$$

*where $C \geq \gamma - 1$. The last term in (1.25) is optimal in the sense that if $\gamma = 1$, there is not positive constant $C$ such that (1.25) holds.*

### 1.3.2 Part II: Hardy-Sobolev type inequalities for $p > n$

Motivated by the aforementioned results of [BrV], [FT], [FMzT3] and [FrL] and the fact that other type of improvements in Sobolev's inequality (like functional asymmetry and affine energy) have been extended to all values of $p > 1$, it is natural to consider the cases where $1 < p < n$ for (1.13) and the counterpart to (1.13) and (1.9) inequalities for $p > n$.

Some results in the direction of extending (1.12) and (1.13) in the range $1 < p < n$, were obtained in [AdChR]-Theorem 1.1, [BFT1]-Theorems B, C and 6.4 and also in [AbdCP]-Theorem 1.1. We focus in the case where $p > n$. Our aim is to provide optimal improvements to the sharp Hardy inequalities (1.11) and (1.8) for $p > n$.

**Distance to a point.** First we improve both (1.11) and (1.3) by replacing the $L^p$-norm of the length of the gradient in (1.3) with the sharp $L^p$ Hardy difference involving distance to the origin:

**Theorem V** *Suppose $\Omega$ is a domain in $\mathbb{R}^n$; $n \geq 1$, containing the origin and having finite volume $\mathcal{L}^n(\Omega)$. Letting $p > n$, there exists a constant $C(n,p) > 0$ such that for all $u \in C_c^\infty(\Omega \setminus \{0\})$*

$$\sup_{x\in\Omega}|u(x)| \leq C(n,p)[\mathcal{L}^n(\Omega)]^{1/n-1/p}\left(\int_\Omega |\nabla u|^p \mathrm{d}x - \Big(\frac{p-n}{p}\Big)^p \int_\Omega \frac{|u|^p}{|x|^p}\mathrm{d}x\right)^{1/p}. \tag{1.26}$$

We go even further. Let us first recall that C. B. Morrey sharpened Sobolev's inequality for $p > n \geq 1$ by replacing the supremum norm in (1.3) with an optimal Hölder semi-norm. What he showed is that there exists a positive constant $C(n,p)$ depending only on $n, p$, such that

$$\sup_{\substack{x,y \in \mathbb{R}^n \\ x \neq y}} \left\{ \frac{|u(x) - u(y)|}{|x-y|^{1-n/p}} \right\} \leq C(n,p) \left( \int_{\mathbb{R}^n} |\nabla u|^p \mathrm{d}x \right)^{1/p}, \qquad (1.27)$$

for all $u \in C_c^\infty(\mathbb{R}^n)$, and the modulus of continuity $1 - n/p$ is optimal. It is possible to replace the $L^p$-norm of the length of the gradient in (1.27) with the sharp $L^p$ Hardy difference implied by (1.11), only after considering a logarithmic correction weight for which we obtain the sharp exponent. The central result of [Ps2] is the following optimal Hardy-Morrey inequality

**Theorem VI** *Suppose $\Omega$ is a bounded domain in $\mathbb{R}^n$; $n \geq 1$, containing the origin and let $p > n$. There exist constants $B = B(n,p) \geq 1$ and $C = C(n,p) > 0$ such that for all $u \in C_c^\infty(\Omega \setminus \{0\})$*

$$\sup_{\substack{x,y \in \Omega \\ x \neq y}} \left\{ \frac{|u(x) - u(y)|}{|x-y|^{1-n/p}} X^{1/p}\left(\frac{|x-y|}{D}\right) \right\} \leq C \left( \int_\Omega |\nabla u|^p \mathrm{d}x - \left(\frac{p-n}{p}\right)^p \int_\Omega \frac{|u|^p}{|x|^p} \mathrm{d}x \right)^{1/p}, \quad (1.28)$$

*where $D = B \operatorname{diam}(\Omega)$ and $X(t) = (1 - \log t)^{-1}$; $t \in (0,1]$. Moreover, the weight function $X^{1/p}$ is optimal in the sense that the power $1/p$ cannot be decreased.*

Note that since $p > n$ one is forced to consider functions in $C_c^\infty(\mathbb{R}^n \setminus \{0\})$, i.e. supported away from the origin. This excludes symmetrization techniques as a method of proof. Thus we turn to multidimensional arguments and in particular in Sobolev's integral representation formula. The first step is to show that (1.28) is equivalent to it's counterpart inequality with one point in the Hölder semi-norm taken to be the origin; see Proposition 5.4. In establishing Proposition 5.4 a crucial step is obtaining estimates on balls $B_r$ intersecting $\Omega$ and with arbitrarily small radius. To this end the following local improvement of the sharp Hardy inequality which is of independent interest is proved:

**Theorem VII** *Suppose $\Omega$ is a bounded domain in $\mathbb{R}^n$; $n \geq 2$, containing the origin and let $p > n$ and $1 \leq q < p$. There exist constants $\Theta = \Theta(n, p, q) \geq 0$ and $C = C(n, p, q) > 0$ such that for all $u \in C_c^\infty(\Omega \setminus \{0\})$, any open ball $B_r$ with $r \in (0, \operatorname{diam}(\Omega))$, and any $D \geq e^\Theta \operatorname{diam}(\Omega)$*

$$r^{n/p} X^{1/p}(r/D) \left( \frac{1}{|B_r|} \int_{B_r} \frac{|u|^q}{|x|^q} \mathrm{d}x \right)^{1/q} \leq C \left( \int_\Omega |\nabla u|^p \mathrm{d}x - \left(\frac{p-n}{p}\right)^p \int_\Omega \frac{|u|^p}{|x|^p} \mathrm{d}x \right)^{1/p}, (1.29)$$

*where $X(t) = (1 - \log t)^{-1}$; $t \in (0,1]$.*

The exponent $1/p$ on the logarithmic factor $X^{1/p}$ is translated as the optimal exponent in (1.28). To obtain this exponent in (1.29) we carefully estimate a trace term on the boundary of $B_r$. Let us note here that if one restricts to the family of open balls $B_r$ containing the origin, then Theorem VII remains valid for $p < n$ (with the factor $((n-p)/p)^p$ instead of $((p-n)/p)^p$ in the Hardy difference).

We finally note that the second important ingredient in the proof of Theorem VI is to show the optimality of the exponent $1/p$. This is done by finding a suitable family of functions that plays the role of a minimizing sequence for inequality (1.28).

**Distance to the boundary.** Next we improve both (1.11) and (1.3) by replacing the $L^p$-norm of the length of the gradient in (1.3) with the sharp $L^p$ Hardy difference involving distance to the boundary:

**Theorem VIII** *Suppose $\Omega$ is a domain in $\mathbb{R}^n$; $n \geq 1$ having finite volume $\mathcal{L}^n(\Omega)$ and such that condition $(\mathfrak{C})$ is satisfied. Letting $p > n$, there exists a constant $C(n,p) > 0$ such that for all $u \in C_c^\infty(\Omega)$*

$$\sup_{x \in \Omega} |u(x)| \leq C(n,p)[\mathcal{L}^n(\Omega)]^{1/n - 1/p} \left( \int_\Omega |\nabla u|^p \mathrm{d}x - \left(\frac{p-1}{p}\right)^p \int_\Omega \frac{|u|^p}{d^p} \mathrm{d}x \right)^{1/p}. \quad (1.30)$$

The analogous to Theorem VI in the case where $\Omega$ is the ball $B_R$ reads as follows

**Theorem IX** *Let $p > n \geq 2$. There exist constants $b = b(n,p) \geq 1$ and $c = c(n,p) > 0$ such that for all $u \in C_c^\infty(B_R)$*

$$\sup_{\substack{x,y \in B_R \\ x \neq y}} \left\{ \frac{|u(x) - u(y)|}{|x-y|^{1-n/p}} X^{1/p}\left(\frac{|x-y|}{D}\right) \right\} \leq$$

$$c \left( \int_{B_R} |\nabla u|^p \mathrm{d}x - \left(\frac{p-1}{p}\right)^p \int_{B_R} \frac{|u|^p}{d^p} \mathrm{d}x \right)^{1/p}, \quad (1.31)$$

*where $D = 2bR$ and $X(t) = (1 - \log t)^{-1}$; $t \in (0,1]$.*

We do not know if the exponent $1/p$ cannot be decreased. In fact this is the case in the one dimensional case. More precisely we have the the following result

**Theorem X** *Let $p > 1$. There exist constants $b = b(p) \geq 1$ and $c = c(p) > 0$ such that for all $u \in C_c^\infty((0,R))$*

$$\sup_{\substack{x,y \in (0,R) \\ x \neq y}} \left\{ \frac{|u(x) - u(y)|}{|x-y|^{1-1/p}} X^{1/p}\left(\frac{|x-y|}{D}\right) \right\} \leq$$

$$c \left( \int_0^R |u'|^p \mathrm{d}x - \left(\frac{p-1}{p}\right)^p \int_0^R \frac{|u|^p}{d^p} \mathrm{d}x \right)^{1/p}, \quad (1.32)$$

*where $D = bR$ and $X(t) = (1 - \log t)^{-1}$; $t \in (0,1]$. The weight function $X^{1/p}$ is optimal in the sense that the power $1/p$ cannot be decreased.*

# Chapter 2

# Preliminary results

In this chapter we gather various inequalities which will be used in this thesis. In particular in §2.1 we will prove various multidimensional weighted Hardy inequalities, most of which are weighted versions of the results in [BFT1]. The results will be used intensively throughout §5 and §6. In §2.2 we will prove the well-known Sobolev inequalities. For $1 \leq p < n$ we give a proof based on the one dimensional weighted Hardy inequality. For $p > n$ we give the classic proof, which is the one we adopt later in §5 and §6.

## 2.1 Hardy inequalities with weights

Our aim in this section is to obtain higher dimensional versions of the following theorems

**Theorem 2.1.** *Let $s \neq 1, q \geq 1$ and $R \in (0, \infty]$. For all $v \in W^{1,q}((0,R))$ such that $v(R) = 0$ if $s < 1$, or $v(0) = 0$ if $s > 1$, there holds*

$$\int_0^R \frac{|v'|^q}{t^{s-q}} dt \geq \left|\frac{s-1}{q}\right|^q \int_0^R \frac{|v|^q}{t^s} dt. \tag{2.1}$$

*The constant is the best possible.*

This is Theorem 330 in the classical book [HLP]. The non-weighted inequality, i.e. for $s = q > 1$, was established earlier by Hardy and Landau (see [Hrd] and [Lnd]). When $R < \infty$ and $s > 1$, we can easily obtain the following equivalent statement

**Theorem 2.2.** *Let $s > 1, q \geq 1$ and $R \in (0, \infty)$. For all $v \in W_0^{1,q}((0,R))$ there holds*

$$\int_0^R \frac{|v'|^q}{d^{s-q}} dt \geq \left(\frac{s-1}{q}\right)^q \int_0^R \frac{|v|^q}{d^s} dt, \tag{2.2}$$

*where $d \equiv d(t) = \min\{t, R-t\}$. The constant is the best possible.*

### 2.1.1 Trace inequalities

We start with the multidimensional counterpart of Theorem 2.1. It is a consequence of the following Lemma (see the remark that follows it) which states the trace Hardy inequality with weights involving distance to the origin. It's proof consists merely of an integration by parts and Young's inequality.

**Lemma 2.3.** *Let $V$ be a domain in $\mathbb{R}^n$; $n \geq 2$, having locally Lipschitz boundary. Denote by $\vec{\nu}(x)$ the exterior unit normal vector defined at almost every $x \in \partial V$. For all $q \geq 1$, all $s \neq n$ and any $v \in C_c^\infty(\mathbb{R}^n \setminus \{0\})$, there holds*

$$\int_V \frac{|\nabla v|^q}{|x|^{s-q}} \mathrm{d}x - \frac{s-n}{q} \left|\frac{s-n}{q}\right|^{q-2} \int_{\partial V} \frac{|v|^q}{|x|^s} x \cdot \vec{\nu}(x) \mathrm{d}S_x \geq \left|\frac{s-n}{q}\right|^q \int_V \frac{|v|^q}{|x|^s} \mathrm{d}x. \qquad (2.3)$$

**Proof.** Integration by parts gives

$$\int_V \nabla|v| \cdot \frac{x}{|x|^s} \mathrm{d}x = -\int_V |v| \operatorname{div}\left(\frac{x}{|x|^s}\right) \mathrm{d}x + \int_{\partial V} |v| \frac{x}{|x|^s} \cdot \vec{\nu} \mathrm{d}S_x,$$

and since $\operatorname{div}(x|x|^{-s}) = -(s-n)|x|^{-s}$ we get

$$\int_V \frac{|\nabla v|}{|x|^{s-1}} \mathrm{d}x \geq (s-n) \int_V \frac{|v|}{|x|^s} \mathrm{d}x + \int_{\partial V} \frac{|v|}{|x|^s} x \cdot \vec{\nu} \mathrm{d}S_x, \quad \text{if } s > n,$$

$$\int_V \frac{|\nabla v|}{|x|^{s-1}} \mathrm{d}x \geq -(s-n) \int_V \frac{|v|}{|x|^s} \mathrm{d}x - \int_{\partial V} \frac{|v|}{|x|^s} x \cdot \vec{\nu} \mathrm{d}S_x, \quad \text{if } s < n,$$

where we have also used the fact that $|\nabla|v(x)|| \leq |\nabla v(x)|$ for a.e. $x \in V$ (see [LbL]-Theorem 6.17). We may write both inequalities in one as follows

$$\int_V \frac{|\nabla v|}{|x|^{s-1}} \mathrm{d}x - \frac{s-n}{|s-n|} \int_{\partial V} \frac{|v|}{|x|^s} x \cdot \vec{\nu} \mathrm{d}S_x \geq |s-n| \int_V \frac{|v|}{|x|^s} \mathrm{d}x.$$

This is inequality (2.3) for $q = 1$. Substituting $v$ by $|v|^q$ with $q > 1$, we arrive at

$$\frac{q}{|s-n|} \int_V \frac{|\nabla v||v|^{q-1}}{|x|^{s-1}} \mathrm{d}x - \frac{s-n}{|s-n|^2} \int_{\partial V} \frac{|v|^q}{|x|^s} x \cdot \vec{\nu} \mathrm{d}S_x \geq \int_V \frac{|v|^q}{|x|^s} \mathrm{d}x. \qquad (2.4)$$

The first term on the left of (2.4) can be written as follows

$$\frac{q}{|s-n|} \int_V \frac{|\nabla v||v|^{q-1}}{|x|^{s-1}} \mathrm{d}x \;=\; \int_V \left\{\frac{q}{|s-n|} \frac{|\nabla v|}{|x|^{s/q-1}}\right\}\left\{\frac{|v|^{q-1}}{|x|^{s-s/q}}\right\} \mathrm{d}x$$
$$\leq \;\; \frac{1}{q}\left|\frac{q}{s-n}\right|^q \int_V \frac{|\nabla v|^q}{|x|^{s-q}} \mathrm{d}x + \frac{q-1}{q} \int_V \frac{|v|^q}{|x|^s} \mathrm{d}x,$$

by Young's inequality with conjugate exponents $q$ and $q/(q-1)$. Thus (2.4) becomes

$$\frac{1}{q}\left|\frac{q}{s-n}\right|^q \int_V \frac{|\nabla v|^q}{|x|^{s-q}} \mathrm{d}x - \frac{s-n}{|s-n|^2} \int_{\partial V} \frac{|v|^q}{|x|^s} x \cdot \vec{\nu} \mathrm{d}S_x \geq \frac{1}{q} \int_V \frac{|v|^q}{|x|^s} \mathrm{d}x.$$

Rearranging the constants we arrive at the inequality we sought for. ∎

**Remark 2.4.** Choosing $V = \mathbb{R}^n$ the trace term on the left-hand side vanishes and we get

$$\int_{\mathbb{R}^n} \frac{|\nabla v|^q}{|x|^{s-q}} \mathrm{d}x \geq \left|\frac{s-n}{q}\right|^q \int_{\mathbb{R}^n} \frac{|v|^q}{|x|^s} \mathrm{d}x. \tag{2.5}$$

The constant is the best possible (see for example [Mz2]-§1.3.1). □

The corresponding trace Hardy inequality with weights involving distance from the boundary reads as follows

**Lemma 2.5.** *Let $\Omega \subsetneq \mathbb{R}^n; n \geq 2$, and set $d \equiv d(x) := \operatorname{dist}(x, \mathbb{R}^n \setminus \Omega)$. Let also $V$ be a domain in $\mathbb{R}^n$ having locally Lipschitz boundary and such that $\Omega \cap V \neq \emptyset$. Denote by $\vec{\nu}(x)$ the exterior unit normal vector defined at almost every $x \in \partial V$. For all $q \geq 1$, all $s \neq 1$ and any $v \in C_c^\infty(\Omega)$, there holds*

$$\int_V \frac{|\nabla v|^q}{d^{s-q}} \mathrm{d}x - \frac{s-1}{q} \left|\frac{s-1}{q}\right|^{q-2} \left( \int_V \frac{|v|^q}{d^{s-1}} (-\Delta d) \mathrm{d}x + \int_{\partial V} \frac{|v|^q}{d^{s-1}} \nabla d \cdot \vec{\nu} \mathrm{d}S_x \right)$$
$$\geq \left|\frac{s-1}{q}\right|^q \int_V \frac{|v|^q}{d^s} \mathrm{d}x. \tag{2.6}$$

**Proof.** Integration by parts gives

$$\int_V \nabla |v| \cdot \frac{\nabla d}{d^{s-1}} \mathrm{d}x = -\int_V |v| \operatorname{div}\left(\frac{\nabla d}{d^{s-1}}\right) \mathrm{d}x + \int_{\partial V} |v| \frac{\nabla d}{d^{s-1}} \cdot \vec{\nu} \mathrm{d}S_x,$$

and since $\operatorname{div}(\nabla d/d^{s-1}) = (1-s)/d^s - (-\Delta d)/d^{s-1}$ in the sense of distributions in $\Omega$, we get

$$\int_V \frac{|\nabla v|}{d^{s-1}} \mathrm{d}x - \int_V \frac{|v|}{d^{s-1}}(-\Delta d) \mathrm{d}x - \int_{\partial V} \frac{|v|}{d^{s-1}} \nabla d \cdot \vec{\nu} \mathrm{d}S_x \geq (s-1) \int_V \frac{|v|}{d^s} \mathrm{d}x, \quad \text{if } s > 1,$$

$$\int_V \frac{|\nabla v|}{d^{s-1}} \mathrm{d}x + \int_V \frac{|v|}{d^{s-1}}(-\Delta d) \mathrm{d}x + \int_{\partial V} \frac{|v|}{d^{s-1}} \nabla d \cdot \vec{\nu} \mathrm{d}S_x \geq -(s-1) \int_V \frac{|v|}{d^s} \mathrm{d}x, \quad \text{if } s < 1,$$

where we have also used the fact that $|\nabla |v(x)|| \leq |\nabla v(x)|$ for a.e. $x \in V$ (see [LbL]-Theorem 6.17). We may write both inequalities in one as follows

$$\int_V \frac{|\nabla v|}{d^{s-1}} \mathrm{d}x - \frac{s-1}{|s-1|} \left( \int_V \frac{|v|}{d^{s-1}}(-\Delta d) \mathrm{d}x + \int_{\partial V} \frac{|v|}{d^{s-1}} \nabla d \cdot \vec{\nu} \mathrm{d}S_x \right) \geq |s-1| \int_V \frac{|v|}{d^s} \mathrm{d}x.$$

This is inequality (2.6) for $q = 1$. Substituting $v$ by $|v|^q$ with $q > 1$, we arrive at

$$\frac{q}{|s-1|} \int_V \frac{|\nabla v| |v|^{q-1}}{d^{s-1}} \mathrm{d}x - \frac{s-1}{|s-1|^2} \left( \int_V \frac{|v|^q}{d^{s-1}}(-\Delta d) \mathrm{d}x + \int_{\partial V} \frac{|v|^q}{d^{s-1}} \nabla d \cdot \vec{\nu} \mathrm{d}S_x \right)$$
$$\geq \int_V \frac{|v|^q}{d^s} \mathrm{d}x. \tag{2.7}$$

The first term on the left of (2.7) can be written as follows

$$\frac{q}{|s-1|}\int_V \frac{|\nabla v||v|^{q-1}}{d^{s-1}}\mathrm{d}x = \int_V\left\{\frac{q}{|s-1|}\frac{|\nabla v|}{d^{s/q-1}}\right\}\left\{\frac{|v|^{q-1}}{d^{s-s/q}}\right\}\mathrm{d}x$$
$$\leq \frac{1}{q}\left|\frac{q}{s-1}\right|^q\int_V\frac{|\nabla v|^q}{d^{s-q}}\mathrm{d}x + \frac{q-1}{q}\int_V\frac{|v|^q}{d^s}\mathrm{d}x,$$

by Young's inequality with conjugate exponents $q$ and $q/(q-1)$. Thus (2.7) becomes

$$\frac{1}{q}\left|\frac{q}{s-1}\right|^q\int_V\frac{|\nabla v|^q}{d^{s-q}}\mathrm{d}x - \frac{s-1}{|s-1|^2}\left(\int_V\frac{|v|^q}{d^{s-1}}(-\Delta d)\mathrm{d}x + \int_{\partial V}\frac{|v|^q}{d^{s-1}}\nabla d\cdot\vec{\nu}\mathrm{d}S_x\right)$$
$$\geq \frac{1}{q}\int_V\frac{|v|^q}{d^s}\mathrm{d}x.$$

Rearranging the constants we arrive at the inequality we sought for. ■

Choosing $V \supsetneq \Omega$, the trace term on the left-hand side vanishes and we get

$$\int_\Omega\frac{|\nabla v|^q}{d^{s-q}}\mathrm{d}x \geq \left|\frac{s-1}{q}\right|^q\int_\Omega\frac{|v|^q}{d^s}\mathrm{d}x + \frac{s-1}{q}\left|\frac{s-1}{q}\right|^{q-2}\int_\Omega\frac{|v|^q}{d^{s-1}}(-\Delta d)\mathrm{d}x. \tag{2.8}$$

A few remarks follow. We introduce first the geometric condition:

$$-\Delta d \geq 0, \quad \text{in the sense of distributions in } \Omega. \tag{$\mathfrak{C}$}$$

We will clear up this condition in the next chapter. At the moment note that convex sets satisfy condition ($\mathfrak{C}$).

**Remark 2.6.** If $s > 1$, $\Omega$ has finite measure and satisfies condition ($\mathfrak{C}$), the second constant appearing on the right hand side in (2.8) is optimal. To see this, we choose $v_\varepsilon(x) = (d(x))^{(s-1)/q+\varepsilon} \in W_0^{1,q}(\Omega; d^{-(s-q)}); \varepsilon > 0$, and after simple computations, involving an integration by parts in the denominator, we obtain

$$\frac{\int_\Omega\frac{|\nabla v_\varepsilon|^q}{d^{s-q}}\mathrm{d}x - (\frac{s-1}{q})^q\int_\Omega\frac{|v_\varepsilon|^q}{d^s}\mathrm{d}x}{\int_\Omega\frac{|v_\varepsilon|^q}{d^{s-1}}(-\Delta d)\mathrm{d}x} = \frac{(\frac{s-1}{q}+\varepsilon)^q - (\frac{s-1}{q})^q}{\varepsilon q}$$
$$\to \left(\frac{s-1}{q}\right)^{q-1},$$

as $\varepsilon \downarrow 0$. □

**Remark 2.7.** If $s > 1$ and condition ($\mathfrak{C}$) is satisfied we may cancel the last term on the right-hand side in (2.8), to deduce the multidimensional counterpart of Theorem 2.2

$$\int_\Omega\frac{|\nabla v|^q}{d^{s-q}}\mathrm{d}x \geq \left(\frac{s-1}{q}\right)^q\int_\Omega\frac{|v|^q}{d^s}\mathrm{d}x. \tag{2.9}$$

The constant is the best possible as we can check by testing the ratio left/right-hand side with $u_\varepsilon(x) = (d(x))^{(s-1)/q+\varepsilon}\phi(x) \in W_0^{1,q}(\Omega; d^{-(s-q)}); \varepsilon > 0$, and using the elementary inequality $|a+b|^q \leq |a|^q + c_q(|a|^{q-1}|b| + |b|^q); a, b \in \mathbb{R}^n$ and $q > 1$, in the numerator. Here, $\phi \in C_c^\infty(B_\delta(y))$, $0 \leq \phi \leq 1$ and $\phi \equiv 1$ in $B_{\delta/2}(y)$, for some small but fixed $\delta$. □

**Remark 2.8.** From the preceding remarks it turns out that if $s > 1$, $\Omega$ has finite measure and satisfies condition $(\mathfrak{C})$, then all constants appearing in (2.8) are optimal. □

**Remark 2.9.** Inequality (2.9) for $s = q = 2$ goes back at least to [D1] for convex domains, while for $s = q > 1$ it was established in [MtskS] and [MMP] for convex $\Omega$. The weighted case for convex $\Omega$ was given in [Avkh] using the techniques from [MMP]. The above generalization to domains that satisfy condition $(\mathfrak{C})$ was given in [BFT1] for $s = q$, where inequality (2.8) was also obtained in that case (see Lemma 3.3 of [BFT1] and also [FMzT1], [FMzT2]). □

### 2.1.2 A general lemma

Here we will generalize Lemma 3.3 of [BFT1] to include weights. The proof consists of a change of variables (see [BrM], [BrV], [FT], [GGrM] and [Mz1]) and also some pointwise inequalities.

**Lemma 2.10.** *Let* $s \in \mathbb{R}$, $p > 1$ *and set* $c_1 := (2^{p-1} - 1)^{-1}$ *and* $c_2 := 3p(p-1)/16$.
(i) *For any* $u \in C_c^\infty(\mathbb{R}^n); n \geq 1$, *the following inequalities hold with* $v := |x|^{(n-s)/p}u$

$$\int_{\mathbb{R}^n} \frac{|\nabla u|^p}{|x|^{s-p}}\mathrm{d}x - \left|\frac{s-n}{p}\right|^p \int_{\mathbb{R}^n} \frac{|u|^p}{|x|^s}\mathrm{d}x \geq c_1 \int_{\mathbb{R}^n} |x|^{p-n}|\nabla v|^p \mathrm{d}x, \tag{2.10}$$

$$\int_{\mathbb{R}^n} \frac{|\nabla u|^p}{|x|^{s-p}}\mathrm{d}x - \left|\frac{s-n}{p}\right|^p \int_{\mathbb{R}^n} \frac{|u|^p}{|x|^s}\mathrm{d}x \geq \frac{c_1}{2^{p-2}}\left|\frac{s-n}{p}\right|^{p-2} \int_{\mathbb{R}^n} |x|^{2-n}|v|^{p-2}|\nabla v|^2 \mathrm{d}x, \tag{2.11}$$

*both in case* $p \geq 2$. *If* $1 < p < 2$, *then*

$$\int_{\mathbb{R}^n} \frac{|\nabla u|^p}{|x|^{s-p}}\mathrm{d}x - \left|\frac{s-n}{p}\right|^p \int_{\mathbb{R}^n} \frac{|u|^p}{|x|^s}\mathrm{d}x \geq c_2 \int_{\mathbb{R}^n} \frac{|x|^{2-n}|\nabla v|^2}{(|x||\nabla v| + |\frac{s-n}{p}v|)^{2-p}}\mathrm{d}x. \tag{2.12}$$

(ii) *Let* $\Omega \subsetneq \mathbb{R}^n$; $n \geq 1$, *be open and set* $d \equiv d(x) := \mathrm{dist}(x, \mathbb{R}^n \setminus \Omega)$. *For any* $u \in C_c^\infty(\Omega)$, *the following inequalities hold with* $v := d^{(1-s)/p}u$

$$\begin{aligned}\int_\Omega \frac{|\nabla u|^p}{d^{s-p}}\mathrm{d}x - \left|\frac{s-1}{p}\right|^p \int_\Omega \frac{|u|^p}{d^s}\mathrm{d}x &\geq c_1 \int_\Omega d^{p-1}|\nabla v|^p \mathrm{d}x \\ &\quad + \frac{s-1}{p}\left|\frac{s-1}{p}\right|^{p-2}\int_\Omega |v|^p(-\Delta d)\mathrm{d}x,\end{aligned} \tag{2.13}$$

$$\begin{aligned}\int_\Omega \frac{|\nabla u|^p}{d^{s-p}}\mathrm{d}x - \left|\frac{s-1}{p}\right|^p \int_\Omega \frac{|u|^p}{d^s}\mathrm{d}x &\geq \frac{c_1}{2^{p-2}}\left|\frac{s-1}{p}\right|^{p-2}\int_\Omega d|v|^{p-2}|\nabla v|^2 \mathrm{d}x \\ &\quad + \frac{s-1}{p}\left|\frac{s-1}{p}\right|^{p-2}\int_\Omega |v|^p(-\Delta d)\mathrm{d}x,\end{aligned} \tag{2.14}$$

*both in case* $p \geq 2$. *If* $1 < p < 2$, *then*

$$\int_\Omega \frac{|\nabla u|^p}{d^{s-p}}\mathrm{d}x - \left|\frac{s-1}{p}\right|^p \int_\Omega \frac{|u|^p}{d^s}\mathrm{d}x \geq c_2 \int_\Omega \frac{d|\nabla v|^2}{(d|\nabla v| + |\frac{s-1}{p}v|)^{2-p}}\mathrm{d}x$$
$$+ \frac{s-1}{p}\left|\frac{s-1}{p}\right|^{p-2} \int_\Omega |v|^p(-\Delta d)\mathrm{d}x. \qquad (2.15)$$

**Proof.** Let $\Omega$ be a domain in $\mathbb{R}^n$; $n \geq 2$, and let either $k = 1$ and $\Omega \subsetneq \mathbb{R}^n$, or $k = n$ and $0 \in \Omega$. For any $x \in \Omega$ we define

$$\delta(x) := \begin{cases} d(x), & \text{if } k = 1 \text{ and } \Omega \subsetneq \mathbb{R}^n \\ |x|, & \text{if } k = n \text{ and } 0 \in \Omega. \end{cases}$$

We also introduce the notation $H := (s-k)/p$, where $p \geq 1$, $s \neq k$. By a straightforward calculation

$$\int_\Omega \left[\frac{|\nabla u|^p}{\delta^{s-p}} - |H|^p \frac{|u|^p}{\delta^s}\right]\mathrm{d}x = \int_\Omega \delta^{-k}\Big[(|Hv\nabla\delta - \delta\nabla v|^p - |Hv|^p)\Big]\mathrm{d}x$$
$$= \int_\Omega \delta^{-k}\Big[(|a-b|^p - |a|^p)\Big]\mathrm{d}x,$$

where $a := Hv\nabla\delta$ and $b := \delta\nabla v$. Note that since $|\nabla\delta| = 1$ a.e. in $\Omega$, we also have $|Hv| = |a|$ a.e. in $\Omega$. The first term on the right hand side of (2.10),(2.13) and of (2.11),(2.14) is evident from the first term on the right hand side in (ii)-(a) and (ii)-(b) of Lemma A.2 in the Appendix. The first term on the right hand side in (i) of Lemma A.2 is

$$\frac{3p(p-1)}{16}\frac{|b|^2}{(|a-b| + |a|)^{2-p}}.$$

Since $|a-b| + |a| \leq |b| + 2|a| \leq 2(|b| + |a|)$, we get

$$\frac{3p(p-1)}{16}\frac{|b|^2}{(|a-b| + |a|)^{2-p}} \geq c_2\frac{|b|^2}{(|b| + |a|)^{2-p}},$$

from which the first term on the right hand side in (2.12),(2.15) follows. The common term $-p|a|^{p-2}a \cdot b$ on the right hand side of (ii)-(a), (ii)-(b) and (i) in Lemma A.2, gives the term

$$pH|H|^{p-2} \int_\Omega \delta^{1-k}v|v|^{p-2}\nabla v \cdot \nabla\delta\mathrm{d}x = H|H|^{p-2} \int_\Omega \delta^{1-k}\nabla|v|^p \cdot \nabla\delta\mathrm{d}x$$
$$= -H|H|^{p-2} \int_\Omega |v|^p \mathrm{div}[\delta^{1-k}\nabla\delta]\mathrm{d}x$$
$$= H|H|^{p-2} \int_\Omega \frac{|v|^p}{\delta^k}(k - 1 - \delta\Delta\delta)\mathrm{d}x,$$

which equals $0$ if $k = n$, and equal to the second term on the right hand side of (2.13)-(2.15) if $k = 1$. ∎

Note that by canceling the first term on the right hand side of each one of (2.13), (2.14) and (2.15) and return to the original function in the second term we obtain (2.8).

### 2.1.3   Improved $L^p$; $p > 1$, **Hardy inequalities with weights**

In the work [BFT1] the authors obtained optimal homogeneous improvements with sharp constants for the sharp Hardy inequalities (2.5) for $s = q$ in bounded domains, and (2.9) for $s = q$ in domains with finite inner radius and satisfying condition $(\mathfrak{C})$. In this section we reproduce these improvements for other values of $s$.

**Theorem 2.11.** *Let $p > 1$ and suppose $\Omega$ be an open domain in $\mathbb{R}^n$; $n \geq 1$.*

(i) *Let $s \neq n$. If $\Omega$ is bounded and contains the origin, then there exists a constant $B = B(n,p,s) \geq 1$ such that for all $u \in C_c^\infty(\Omega \setminus \{0\})$*

$$\int_\Omega \frac{|\nabla u|^p}{|x|^{s-p}}\mathrm{d}x - \left|\frac{s-n}{p}\right|^p \int_\Omega \frac{|u|^p}{|x|^s}\mathrm{d}x \geq \frac{1}{2}\frac{p-1}{p}\left|\frac{s-n}{p}\right|^{p-2}\int_\Omega \frac{|u|^p}{|x|^s}X^2(|x|/D)\mathrm{d}x, \quad (2.16)$$

*where $D = B\sup_{x\in\Omega}|x|$ and $X(t) = (1-\log t)^{-1}$; $t \in (0,1]$. The weight function $X^2$ is optimal, in the sense that the power 2 cannot be decreased, and the constant on the right-hand side is the best possible.*

(ii) *Let $s > 1$. If $\Omega$ has finite inner radius and satisfies condition $(\mathfrak{C})$, then there exists a constant $B = B(n,p,s) \geq 1$ such that for all $u \in C_c^\infty(\Omega)$*

$$\int_\Omega \frac{|\nabla u|^p}{d^{s-p}}\mathrm{d}x - \left(\frac{s-1}{p}\right)^p \int_\Omega \frac{|u|^p}{d^s}\mathrm{d}x \geq \frac{1}{2}\frac{p-1}{p}\left(\frac{s-1}{p}\right)^{p-2}\int_\Omega \frac{|u|^p}{d^s}X^2(d/D)\mathrm{d}x, \quad (2.17)$$

*where $D = B\sup_{x\in\Omega}d(x)$ and $X(t) = (1-\log t)^{-1}$; $t \in (0,1]$. The weight function $X^2$ is optimal, in the sense that the power 2 cannot be decreased, and the constant on the right-hand side is the best possible.*

**Proof.** *Proof of* (i). We will give a more direct proof of (2.16) with a different constant using Lemma 2.10 of the previous subsection. To obtain the sharp constant one has to work as in the proof of (ii). Note however that for $p = 2$ the following proof yields also the best constant. Let $D = \sup_{x\in\Omega}|x|$ and setting $u(x) = |x|^{(s-n)/p}v(x)$ we have

$$\begin{aligned}
\int_\Omega \frac{|u|^p}{|x|^s}X^2(|x|/D)\mathrm{d}x &= \int_\Omega \frac{|v|^p}{|x|^n}X^2(|x|/D)\mathrm{d}x \\
&= \int_\Omega |v|^p \mathrm{div}\left\{\frac{X(|x|/D)}{|x|^n}x\right\}\mathrm{d}x \\
&= -p\int_\Omega \frac{|v|^{p-1}}{|x|^n}X(|x|/D)x\cdot\nabla|v|\mathrm{d}x \\
&\leq p\int_\Omega \frac{|v|^{p-1}}{|x|^{n-1}}X(|x|/D)|\nabla v|\mathrm{d}x \quad (2.18)\\
&= p\int_\Omega \left\{\frac{|v|^{p/2-1}}{|x|^{n/2-1}}|\nabla v|\mathrm{d}x\right\}\left\{\frac{|v|^{p/2}}{|x|^{n/2}}X(|x|/D)\right\}\mathrm{d}x \\
&\leq p\left(\int_\Omega \frac{|v|^{p-2}}{|x|^{n-2}}|\nabla v|^2 \mathrm{d}x\right)^{1/2}\left(\int_\Omega \frac{|u|^p}{|x|^s}X^2(|x|/D)\mathrm{d}x\right)^{1/2},
\end{aligned}$$

where we have integrated by parts, used Hölder's inequality and returned to the original function in the second integral. Rearranging and squaring we deduce

$$\int_\Omega \frac{|u|^p}{|x|^s} X^2(|x|/D) \mathrm{d}x \leq p^2 \int_\Omega |x|^{2-n} |v|^{p-2} |\nabla v|^2 \mathrm{d}x.$$

If $p \geq 2$, the result comes from (2.11). If $1 < p < 2$ we proceed from (2.18) as follows

$$\begin{aligned}
\int_\Omega \frac{|u|^p}{|x|^s} X^2(|x|/D) \mathrm{d}x &\leq p \int_\Omega \left\{ \frac{|\nabla v|}{|x|^{n/p-1}} X^{2/p-1}(|x|/D) \right\} \left\{ \frac{|v|^{p-1}}{|x|^{n-n/p}} X^{2-2/p}(|x|/D) \right\} \mathrm{d}x \\
&\leq p \left( \int_\Omega \frac{|\nabla v|^p}{|x|^{n-p}} X^{2-p}(|x|/D) \mathrm{d}x \right)^{1/p} \left( \int_\Omega \frac{|u|^p}{|x|^s} X^2(|x|/D) \mathrm{d}x \right)^{1-1/p},
\end{aligned}$$

where we have used Hölder's inequality and returned to the original function in the second integral. Rearranging and raising to $p$-th power we get

$$\int_\Omega \frac{|u|^p}{|x|^s} X^2(|x|/D) \mathrm{d}x \leq p^p \underbrace{\int_\Omega |x|^{p-n} |\nabla v|^p X^{2-p}(|x|/D) \mathrm{d}x}_{=:J[v]}. \tag{2.19}$$

To estimate $J[v]$ we argue as follows

$$\begin{aligned}
J[v] &= \int_\Omega \left\{ \frac{|x|^{(2-n)p/2} |\nabla v|^p}{(|x||\nabla v| + |\frac{s-n}{p} v|)^{(2-p)p/2}} \right\} \left\{ \frac{|x|^{-(2-p)n/2} X^{2-p}(|x|/D)}{(|x||\nabla v| + |\frac{s-n}{p} v|)^{(p-2)p/2}} \right\} \mathrm{d}x \tag{2.20} \\
&\leq \left( \int_\Omega \frac{|x|^{2-n} |\nabla v|^2}{(|x||\nabla v| + |\frac{s-n}{p} v|)^{2-p}} \mathrm{d}x \right)^{p/2} \left( \int_\Omega \frac{|x|^{-n} X^2(|x|/D)}{(|x||\nabla v| + |\frac{s-n}{p} v|)^{-p}} \mathrm{d}x \right)^{1-p/2} \\
&\leq c_2^{-p/2} (I[u])^{p/2} \left( \int_\Omega |x|^{-n} \left( |x||\nabla v| + |\frac{s-n}{p} v| \right)^p X^2(|x|/D) \mathrm{d}x \right)^{1-p/2},
\end{aligned}$$

where we have used Hölder's inequality and (2.12). Here, by $I[u]$ we denote the Hardy difference as appears in the left-hand side of (2.12). Minkowski's inequality and the fact that $X^2(t) \leq X^{2-p}(t)$ for all $t \in (0, 1]$, yields

$$\begin{aligned}
J[v] &\leq c_2^{-p/2} (I[u])^{p/2} \left\{ (J[v])^{1/p} + \left( \left| \frac{s-n}{p} \right|^p \int_\Omega \frac{|v|^p}{|x|^n} X^2(|x|/D) \mathrm{d}x \right)^{1/p} \right\}^{p(1-p/2)} \\
&\leq c_2^{-p/2} (I[u])^{p/2} \left\{ (J[v])^{1/p} + |s-n|^p (J[v])^{1/p} \right\}^{p(1-p/2)},
\end{aligned}$$

where in the last inequality we have used Lemma-3.2 of [BFT1] with $\alpha = 2$. It turns up that $J[v] \leq C(n, p, s) I[u]$, and the result is evident by (2.19).

*Proof of* (ii). Let $T$ be a vector field on $\Omega$. Integrating by parts and using elementary inequalities, we get

$$\int_\Omega \operatorname{div} T |u|^p \mathrm{d}x \leq p \int_\Omega |T||\nabla u||u|^{p-1} \mathrm{d}x$$
$$= p \int_\Omega \left\{\frac{|\nabla u|}{d^{s/p-1}}\right\} \left\{d^{s/p-1}|T||u|^{p-1}\right\} \mathrm{d}x,$$

for all $u \in C_c^\infty(\Omega)$. By Young's inequality: $ab \leq a^p/p + (p-1)b^{(p-1)/p}/p$, we arrive at

$$\int_\Omega \operatorname{div} T |u|^p \mathrm{d}x \leq \int_\Omega \frac{|\nabla u|^p}{d^{s-p}} \mathrm{d}x + (p-1) \int_\Omega d^{(s-p)/(p-1)}|T|^{\frac{p}{p-1}}|u|^p \mathrm{d}x.$$

Hence,

$$\int_\Omega \frac{|\nabla u|^p}{d^{s-p}} \mathrm{d}x \geq \int_\Omega \left[\operatorname{div} T - (p-1)d^{(s-p)/(p-1)}|T|^{\frac{p}{p-1}}\right]|u|^p \mathrm{d}x.$$

In view of this, inequality (2.17) will be proved if we could find a vector field $T$ such that the following inequality to hold (at least in the sense of distributions)

$$\operatorname{div} T - (p-1)d^{(s-p)/(p-1)}|T|^{\frac{p}{p-1}} \geq H^p \frac{1}{d^s}\left[1 + \frac{p-1}{2pH^2}X^2(d/D)\right],$$

where we have set $H := (s-1)/p$. To proceed we make a specific choice of $T$. We take

$$T = -H^{p-1}\left[1 - \frac{p-1}{s-1}X(d/D) + aX^2(d/D)\right]\frac{1}{d^{s-1}}\nabla d,$$

where $a$ is a free parameter to be chosen later. In any case $a$ will be such that the quantity inside brackets is nonnegative on $\Omega$. By a straightforward calculation we get

$$\operatorname{div} T = -H^{p-1}\left[-\frac{p-1}{s-1}X^2(d/D) + 2aX^3(d/D)\right]\frac{1}{d^s}$$
$$+ H^{p-1}\left[1 - \frac{p-1}{s-1}X(d/D) + aX^2(d/D)\right]\frac{s-1-d\Delta d}{d^s}$$
$$\geq -H^{p-1}\left[-\frac{p-1}{s-1}X^2(d/D) + 2aX^3(d/D)\right]\frac{1}{d^s}$$
$$+ H^{p-1}\left[1 - \frac{p-1}{s-1}X(d/D) + aX^2(d/D)\right]\frac{s-1}{d^s}.$$

where we have used condition $(\mathfrak{C})$. Thus, we have

$$\operatorname{div} T - (p-1)d^{(s-p)/(p-1)}|T|^{\frac{p}{p-1}} \geq pH^p\left[1 - \frac{p-1}{s-1}X(d/D) + aX^2(d/D)\right]\frac{1}{d^s}$$
$$+ H^{p-1}\left[\frac{p-1}{s-1}X^2(d/D) - 2aX^3(d/D)\right]\frac{1}{d^s}$$
$$- (p-1)H^p\left[1 - \frac{p-1}{s-1}X(d/D) + aX^2(d/D)\right]^{\frac{p}{p-1}}\frac{1}{d^s}.$$

It then follows that for (2.17) to hold, it is enough to establish the inequality

$$f(t) \geq 1 + \frac{p-1}{2pH^2}t^2, \qquad t \in [0, M], \tag{2.21}$$

where $M = M(D) := \sup_{x \in \Omega} X(d(x)/D) \leq 1$ and

$$f(t) := p\Big(1 - \frac{p-1}{pH}t + at^2\Big) + \frac{1}{H}\Big(\frac{p-1}{pH}t^2 - 2at^3\Big) - (p-1)\Big(1 - \frac{p-1}{pH}t + at^2\Big)^{\frac{p}{p-1}}.$$

From Taylor's formula we have that

$$f(t) = f(0) + f'(0)t + \frac{1}{2}f''(\xi_t)t^2; \quad 0 \leq \xi_t \leq t \leq M. \tag{2.22}$$

We have $f(0) = 1$. Moreover, after some simple calculations we find

$$f'(0) = 0, \quad f''(0) = \frac{p-1}{pH^2} \quad \text{and} \quad f'''(0) = -\frac{6a}{H} + \frac{(2-p)(p-1)}{p^2H^3}.$$

We choose $a$ so that $f'''(0) > 0$, that is, $a < (2-p)(p-1)/(6(s-1)^2)$. Hence $f''$ is an increasing function in some interval of the form $(0, M_0)$. Consequently, for $t \in (0, M_0)$

$$f''(\xi_t) \geq f''(0) = \frac{p-1}{pH^2}.$$

It then follows from (2.22)

$$f(t) \geq 1 + \frac{p-1}{2pH^2}t^2; \quad t \in [0, M_0].$$

It is clear that we can choose small enough $M_0$, depending only on $n, s, p$, such that $a < (2-p)(p-1)/(6(s-1)^2)$ and also $1 - \frac{p-1}{pH}t + at^2 \geq 0$ for all $0 < t < M_0$. Since $X(t) = (1 - \log(t))^{-1}$, the condition $X(d/D) \leq M_0$ is equivalent to $D \geq B\sup_{x \in \Omega} d(x)$, where $B = e^{1/M_0 - 1}$. ∎

The proofs for the optimality of the constant and the exponent on $X$ in the remainder terms of Theorem 2.11, are essentially the same to the proof given in [BFT1]-§5.

## 2.2   Sobolev inequalities

Here we state and prove the well-known Sobolev inequalities. For $1 \leq p < n$, although there are many other proofs, we chose an approach via Lorentz spaces based on Hardy's inequality (2.1). For $p > n$ we follow [GTr]-§7.7-7.8 with some minor modifications.

### 2.2.1 The case $1 \leq p < n$

We first recall the definition of Lorentz spaces and some of their basic properties.

**Definition 2.12.** Let $p, q > 0$. The *Lorentz space* $L^{p,q}$ is the collection of all measurable functions $f$ defined on $\mathbb{R}^n$, such that $[f]_{L^{p,q}} < \infty$, where

$$[f]_{L^{p,q}} := \left( \frac{q}{p} \int_0^\infty (f^*(t))^q t^{q/p-1} \mathrm{d}t \right)^{1/q}.$$

Here, $f^*$ denotes the *decreasing rearrangement* of $f$,

$$f^*(t) := \inf\{s \geq 0 \text{ s.t. } \mu_f(s) \leq t\},$$

where $\mu_f$ is the distribution function of $f$, i.e.

$$\mu_f(s) := \mathcal{L}^n(\{x \in \mathbb{R}^n \text{ s.t. } |f(x)| > s\}).$$

The first two properties that follow are elementary. The proof of the third can be found in [H] or [Tl2].

**Proposition 2.13.** (i) *We have* $[f]_{L^{p,p}} = [f]_{L^p}$.
 (ii) *If $\chi_\Omega$ is the characteristic function of a set of finite volume, then $[\chi_\Omega]_{L^{p,q}} = [\mathcal{L}^n(\Omega)]^{1/p}$ for all $p, q > 0$.*
 (iii) *Let $p > 0$ and $q_2 \geq q_1 > 0$. There holds $[f]_{L^{p,q_1}} \geq [f]_{L^{p,q_2}}$.*

We now prove the Sobolev inequality. Note that for $p = 1$ the constant we will obtain is $n\omega_n^{1/n}$, which is known to be the optimal constant in this case, see [FFl] and [Mz].

**Theorem 2.14.** *Let $1 \leq p < n$ and set $p^* := np/(n-p)$. For any $u \in C_c^\infty(\mathbb{R}^n)$ there holds*

$$\left( \int_{\mathbb{R}^n} |\nabla u(x)|^p \mathrm{d}x \right)^{1/p} \geq \left( \frac{n-p}{p} \right)^{1-1/p} \left( \frac{n}{p} \right)^{1/p} \omega_n^{1/n} \left( \int_{\mathbb{R}^n} |u(x)|^{p^*} \mathrm{d}x \right)^{1/p^*}.$$

**Proof.** Let $u \in C_c^\infty(\mathbb{R}^n)$. By the Pòlya-Szegö inequality (see [K] or [Fsc]-Theorem 3.1)

$$\int_{\mathbb{R}^n} |\nabla u(x)|^p \mathrm{d}x \geq \int_{\mathbb{R}^n} |\nabla u^\star(x)|^p \mathrm{d}x, \tag{2.23}$$

where $u^\star$ is the *symmetric rearrangement* of $u$, defined by

$$u^\star := \begin{cases} u^*(\omega_n |x|^n), & x \in B_R(0) \\ 0, & \text{otherwise,} \end{cases}$$

where $R := [\mathcal{L}^n(\text{sprt}\{u\})/\omega_n]^{1/n}$. We calculate

$$
\begin{aligned}
\int_{\mathbb{R}^n} |\nabla u^\star(x)|^p \mathrm{d}x &= \int_{\mathbb{R}^n} |u^{*\prime}(\omega_n|x|^n)\nabla(\omega_n|x|^n)|^p \mathrm{d}x \\
&= (n\omega_n)^p \int_{\mathbb{R}^n} |u^{*\prime}(\omega_n|x|^n)|^p |x|^{(n-1)p} \mathrm{d}x \\
&= (n\omega_n)^{p+1} \int_0^\infty |u^{*\prime}(\omega_n r^n)|^p r^{(n-1)(p+1)} \mathrm{d}r \\
&= (n\omega_n^{1/n})^p \int_0^\infty |u^{*\prime}(t)|^p t^{p-p/n} \mathrm{d}t.
\end{aligned}
\quad (2.24)
$$

Coupling (2.23), (2.24) and the one dimensional Hardy inequality (2.1) for $v = u^*$ and $q = p, s = p/n$, we obtain

$$
\begin{aligned}
\int_{\mathbb{R}^n} |\nabla u(x)|^p \mathrm{d}x &\geq \left(\frac{n-p}{p}\omega_n^{1/n}\right)^p \int_0^\infty |u^*(t)|^p t^{-p/n} \mathrm{d}t \\
&= \left(\frac{n-p}{p}\omega_n^{1/n}\right)^p \frac{n}{n-p} [u]^p_{L^{p^*,p}} \\
&\geq \left(\frac{n-p}{p}\omega_n^{1/n}\right)^p \frac{n}{n-p} [u]^p_{L^{p^*,p^*}} \\
&= \left(\frac{n-p}{p}\right)^{p-1} \frac{n}{p} \omega_n^{p/n} \left(\int_{\mathbb{R}^n} |u(x)|^{p^*} \mathrm{d}x\right)^{p/p^*},
\end{aligned}
\quad (2.25)
$$

where the last inequality follows from Proposition 2.13-(iii). ∎

**Remark 2.15.** The Hardy-Littlewood inequality (see [K] or [Ln]-Theorem 16.9) implies

$$
\begin{aligned}
\int_{\mathbb{R}^n} \frac{|u(x)|^p}{|x|^p} \mathrm{d}x &\leq \int_{\mathbb{R}^n} \frac{(u^\star(x))^p}{|x|^p} \mathrm{d}x \\
&= \int_{\mathbb{R}^n} \frac{(u^*(\omega_n|x|^n))^p}{|x|^p} \mathrm{d}x \\
&= n\omega_n \int_0^\infty (u^*(\omega_n r^n))^p r^{n-1-p} \mathrm{d}r \\
&= \omega_n^{p/n} \int_0^\infty (u^*(t))^p t^{-p/n} \mathrm{d}t.
\end{aligned}
$$

Hence, Hardy's inequality (2.5) for $s = q = p \in [1, n)$ with the sharp constant is obtained from (2.25). □

### 2.2.2 The case $p > n \geq 1$

We start with the case $n = 1$. Set $I = (\alpha, \beta)$; $\alpha < \beta$. If $u \in C_c^\infty(I)$, then for any $x \in I$ we have

$$|u(x)| \leq \left| \int_\alpha^x u'(t)\mathrm{d}t \right| \leq \int_\alpha^x |u'(t)|\mathrm{d}t,$$

and

$$|u(x)| \leq \left| -\int_x^\beta u'(t)\mathrm{d}t \right| \leq \int_x^\beta |u'(t)|\mathrm{d}t.$$

Adding we obtain

$$\begin{aligned} |u(x)| &\leq \frac{1}{2} \int_I |u'(t)|\mathrm{d}t \\ &\leq \frac{1}{2}[\mathcal{L}^1(I)]^{1-1/p} \left( \int_I |u'(t)|^p \mathrm{d}t \right)^{1/p}, \end{aligned} \qquad (2.26)$$

by Hölder's inequality. The constant $1/2$ is optimal and a standard scaling argument reveals that the exponent $1 - 1/p$ on $\mathcal{L}^1(I)$ cannot be decreased. Now we prove the direct multidimensional analog of this inequality. The steps are essentially the same. We need first an integral representation formula for the values of the function.

**Lemma 2.16.** *Let $\Omega$ be an open set in $\mathbb{R}^n$; $n \geq 2$. For any $u \in C_c^\infty(\Omega)$ and any $x \in \Omega$ there holds*

$$u(x) = \frac{1}{n\omega_n} \int_\Omega \frac{(x-z) \cdot \nabla u(z)}{|x-z|^n} \mathrm{d}z.$$

**Proof.** Using polar coordinates around $x$ and then changing variables by $z = x + \rho y$, we have

$$\begin{aligned} \int_\Omega \frac{(x-z) \cdot \nabla u(z)}{|x-z|^n} \mathrm{d}z &= \int_0^\infty \int_{\partial B_\rho(x)} \frac{(x-z) \cdot \nabla u(z)}{|x-z|^n} \mathrm{d}S_z \mathrm{d}\rho \\ &= -\int_0^\infty \int_{\partial B_1(0)} y \cdot \nabla u(x + \rho y) \mathrm{d}S_y \mathrm{d}\rho \\ &= -\int_{\partial B_1(0)} \int_0^\infty y \cdot \nabla u(x + \rho y) \mathrm{d}\rho \mathrm{d}S_y \\ &= -\int_{\partial B_1(0)} \int_0^\infty \frac{\mathrm{d}}{\mathrm{d}\rho}\Big[u(x + \rho y)\Big] \mathrm{d}\rho \mathrm{d}S_y \\ &= -\int_{\partial B_1(0)} \Big[0 - u(x)\Big] \mathrm{d}S_y \\ &= n\omega_n u(x), \end{aligned}$$

where we have also used Fubini's Theorem. ∎

**Theorem 2.17.** *Let $\Omega$ be an open set in $\mathbb{R}^n$; $n \geq 2$, having finite volume $\mathcal{L}^n(\Omega)$. If $p > n$, for any $u \in C_c^\infty(\Omega)$ there holds*

$$\sup_{x \in \Omega} |u(x)| \leq n^{-1/p} \omega_n^{-1/n} \left(\frac{p-1}{p-n}\right)^{1-1/p} [\mathcal{L}^n(\Omega)]^{1/n - 1/p} \left(\int_\Omega |\nabla u|^p \mathrm{d}x\right)^{1/p}.$$

*The constant appearing on the right hand side is optimal and the exponent $1/n - 1/p$ on $\mathcal{L}^n(\Omega)$ cannot be decreased.*

**Proof.** By Lemma 2.16 we get

$$|u(x)| \leq \frac{1}{n\omega_n} \int_\Omega \frac{|\nabla u(z)|}{|x-z|^{n-1}} \mathrm{d}z,$$

and applying Hölder's inequality

$$|u(x)| \leq \frac{1}{n\omega_n} \Bigg(\underbrace{\int_\Omega |x-z|^{-(n-1)p/(p-1)} \mathrm{d}z}_{=: M(x)}\Bigg)^{1-1/p} \left(\int_\Omega |\nabla u(z)|^p \mathrm{d}z\right)^{1/p}. \tag{2.27}$$

Note that $M(x)$ is finite since $(n-1)p/(p-1) < n$ if and only if $p > n$. To estimate $M(x)$ we set $R := [\mathcal{L}^n(\Omega)/\omega_n]^{1/n}$, so that the volume of a ball with radius $R$ to be equal to the volume of $\Omega$. Then $M(x)$ increases if we change the domain of integration from $\Omega$ to $B_R(x)$. Therefore

$$\begin{aligned} M(x) &\leq \int_{B_R(x)} |x-z|^{-(n-1)p/(p-1)} \mathrm{d}z \\ &= n\omega_n \frac{p-1}{p-n} [\mathcal{L}^n(\Omega)/\omega_n]^{(p-n)/n(p-1)}. \end{aligned}$$

Substituting this in (2.27) we are done. ∎

**Remark 2.18.** The best constant and the family of extremal functions (see (1.4)) have been found in [Tl2] using symmetrization techniques. However, the above method of proof has also revealed the best constant. Note also that a scaling argument shows that this inequality is also optimal in the sense that it fails if $1/n - 1/p$ is replaced by any exponent $q < 1/n - 1/p$ on $\mathcal{L}^n(\Omega)$. □

### 2.2.3 Morrey's inequality

Starting with the case $n = 1$, let $u \in C_c^\infty(\mathbb{R})$. Then for $y < x$ we have

$$|u(x) - u(y)| = \left|\int_y^x u'(t) \mathrm{d}t\right| \leq \int_y^x |u'(t)| \mathrm{d}t.$$

By Hölder's inequality we obtain

$$|u(x) - u(y)| = (x-y)^{1-1/p}\left(\int_y^x |u'(t)|^p dt\right)^{1/p}$$

$$\leq (x-y)^{1-1/p}\left(\int_\mathbb{R} |u'(t)|^p dt\right)^{1/p}$$

Rearranging we arrive at

$$\sup_{\substack{x,y\in\mathbb{R}\\ x\neq y}}\left\{\frac{|u(x)-u(y)|}{|x-y|^{1-1/p}}\right\} \leq \left(\int_\mathbb{R} |\nabla u|^p dx\right)^{1/p}.$$

As in the previous section we will prove the direct multidimensional analog of this inequality and the steps are essentially the same with the case $n = 1$. We first obtain a local substitute of the integral representation formula of Lemma 2.16.

**Lemma 2.19.** *For any $u \in C^\infty(B_r)$ and any $x \in B_r$ there holds*

$$|u(x) - u_{B_r}| \leq \frac{2^n}{n\omega_n}\int_{B_r} \frac{|\nabla u(z)|}{|x-z|^{n-1}} dz.$$

**Proof.** Letting $x, y \in B_r$ we have

$$u(x) - u(y) = -\int_0^{|y-x|} \frac{d}{d\rho}\left[u\left(x + \rho\frac{y-x}{|y-x|}\right)\right]d\rho$$

$$= -\int_0^{|y-x|} \nabla u\left(x + \rho\frac{y-x}{|y-x|}\right)\cdot\frac{y-x}{|y-x|}d\rho.$$

Integrating this with respect to $y$ in $B_r$ we get

$$u(x) - u_{B_r} = -\frac{1}{\omega_n r^n}\int_{B_r}\int_0^{|y-x|} \nabla u\left(x + \rho\frac{y-x}{|y-x|}\right)\cdot\frac{y-x}{|y-x|}d\rho dy.$$

Hence

$$|u(x) - u_{B_r}| \leq \frac{1}{\omega_n r^n}\int_{B_r}\int_0^{|y-x|} F\left(x + \rho\frac{y-x}{|y-x|}\right)d\rho dy,$$

where we have set

$$F(z) = \begin{cases} |\nabla u(z)|, & \text{if } z \in B_r \\ 0, & \text{if } z \notin B_r. \end{cases}$$

Since $x, y \in B_r$ we have $|y - x| \leq 2r$ and $B_r \subset B_{2r}(x)$. Thus

$$\begin{aligned}
|u(x) - u_{B_r}| &\leq \frac{1}{\omega_n r^n} \int_{B_{2r}(x)} \int_0^{2r} F\Big(x + \rho \frac{y-x}{|y-x|}\Big) d\rho dy \\
&= \frac{1}{\omega_n r^n} \int_0^{2r} \int_{\partial B_t(x)} \int_0^{2r} F\Big(x + \rho \frac{y-x}{t}\Big) d\rho dS_y dt \\
&= \frac{1}{\omega_n r^n} \int_0^{2r} \int_0^{2r} \int_{\partial B_t(x)} F\Big(x + \rho \frac{y-x}{t}\Big) dS_y dt d\rho \\
&= \frac{1}{\omega_n r^n} \int_0^{2r} \rho^{1-n} \int_0^{2r} t^{n-1} \int_{\partial B_\rho(x)} F(z) dS_z dt d\rho,
\end{aligned}$$

where we have used polar coordinates, Fubini's Theorem and the change of variables $z = x + \rho \frac{y-x}{t}$. Computing the $t$-integral and using polar coordinates once more

$$\begin{aligned}
|u(x) - u_{B_r}| &\leq \frac{2^n}{n\omega_n} \int_0^{2r} \rho^{1-n} \int_{\partial B_\rho(x)} F(z) dS_z d\rho \\
&= \frac{2^n}{n\omega_n} \int_0^{2r} \int_{\partial B_\rho(x)} \frac{F(z)}{|x-z|^{n-1}} dS_z d\rho \\
&= \frac{2^n}{n\omega_n} \int_{B_{2r}(x)} \frac{F(z)}{|x-z|^{n-1}} dz \\
&= \frac{2^n}{n\omega_n} \int_{B_r} \frac{|\nabla u(z)|}{|x-z|^{n-1}} dz,
\end{aligned}$$

by the definition of $F$. ∎

**Theorem 2.20.** *If $p > n \geq 2$, for any $u \in C_c^\infty(\mathbb{R}^n)$ there holds*

$$\sup_{\substack{x,y \in \mathbb{R}^n \\ x \neq y}} \left\{ \frac{|u(x) - u(y)|}{|x-y|^{1-n/p}} \right\} \leq 2^{n+1} (n\omega_n)^{-1/p} \Big(\frac{p-1}{p-n}\Big)^{1-1/p} \bigg( \int_{\mathbb{R}^n} |\nabla u|^p dx \bigg)^{1/p}.$$

**Proof.** Letting $x, y \in \mathbb{R}^n$ with $x \neq y$ we consider a ball $B_r$ of radius $r := |x - y|$ containing $x, y$. By Lemma 2.19 we get

$$\begin{aligned}
|u(x) - u(y)| &\leq |u(x) - u_{B_r}| + |u(y) - u_{B_r}| \\
&\leq \frac{2^n}{n\omega_n} \bigg( \underbrace{\int_{B_r} \frac{|\nabla u(z)|}{|x-z|^{n-1}} dz}_{=: J(x)} + \underbrace{\int_{B_r} \frac{|\nabla u(z)|}{|y-z|^{n-1}} dz}_{=: J(y)} \bigg). \quad (2.28)
\end{aligned}$$

We will bound $J(x)$ independently on $x$ so that the same estimate holds also for $J(y)$. Hölder's inequality gives

$$J(x) \leq \bigg( \int_{B_r} |x-z|^{-(n-1)p/(p-1)} dz \bigg)^{1-1/p} \bigg( \int_{B_r} |\nabla u|^p dz \bigg)^{1/p}.$$

Both integrals increase if we integrate over $B_r(x)$ and $\mathbb{R}^n$ respectively. Hence

$$\begin{aligned} J(x) &\leq \left(\int_{B_r(x)} |x-z|^{-(n-1)p/(p-1)} \mathrm{d}z\right)^{1-1/p} \left(\int_{\mathbb{R}^n} |\nabla u|^p \mathrm{d}z\right)^{1/p} \\ &= \left(n\omega_n \frac{p-1}{p-n}\right)^{1-1/p} r^{1-n/p} \left(\int_{\mathbb{R}^n} |\nabla u|^p \mathrm{d}z\right)^{1/p}. \end{aligned}$$

Substituting this into (2.28) twice we obtain

$$|u(x) - u(y)| \leq 2^{n+1}(n\omega_n)^{-1/p}\left(\frac{p-1}{p-n}\right)^{1-1/p}|x-y|^{1-n/p}\left(\int_{\mathbb{R}^n} |\nabla u|^p \mathrm{d}z\right)^{1/p}.$$

The proof follows by rearranging and taking the supremum over all $x, y \in \mathbb{R}^n$ with $x \neq y$. ∎

**Remark 2.21.** Morrey's inequality states a stronger result than Sobolev's inequality for $p > n$. Moreover it implies a more general version of it. We can use Morrey's inequality to obtain Sobolev's inequality for $p > n$ for a more general class of domains than ones having finite volume, and in particular for domains having finite inner radius. To see this assume that $\Omega$ has finite inner radius $R := \sup_{x \in \Omega} d(x) < \infty$. Let $y = \xi(x)$ in Morrey's inequality where $\xi(x) \in \mathbb{R}^n \setminus \Omega$ is such that $d(x) := \operatorname{dist}(x, \mathbb{R}^n \setminus \Omega) = |x - \xi(x)|$. Then $u(y) = 0$ and we obtain

$$\sup_{x \in \Omega} \left\{ \frac{|u(x)|}{(d(x))^{1-n/p}} \right\} \leq C(n, p) \left(\int_\Omega |\nabla u|^p \mathrm{d}x\right)^{1/p}$$

Hence,

$$\sup_{x \in \Omega} |u(x)| \leq C(n, p) R^{1-n/p} \left(\int_\Omega |\nabla u|^p \mathrm{d}x\right)^{1/p}.$$

□

# Chapter 3

# The distance function

Here we gather some facts for the distance function, a weight in Hardy inequalities. The results will play an essential role in §4. In the first section, among other properties we recall the local semiconcavity of the distance function and the semiconcavity of the squared distance function. We also connect sets of positive reach to these properties and obtain weak upper bounds for the distributional Laplacian of the distance function. In the second section we give a characterization of mean convex sets obtained in [Ps1].

## 3.1 Basic properties and sets with positive reach

**Definition 3.1.** If $\emptyset \neq K \subsetneq \mathbb{R}^n$ is closed, then $d_K : \mathbb{R}^n \to [0, \infty)$ is the *distance function* to $K$, that is

$$d_K(x) := \mathrm{dist}(x, K) = \inf_{y \in K} |x - y|,$$

whenever $x \in \mathbb{R}^n$. Furthermore, $K_1$ is the set of points in $\mathbb{R}^n$ which have a unique closest point on $K$, namely

$$K_1 := \{x \in \mathbb{R}^n : \exists!\, y \in K \text{ such that } d_K(x) = |x - y|\},$$

and the function $\xi_K : K_1 \to K$ associates each $x \in K_1$ with its unique $y \in K$ for which $d_K(x) = |x - y|$.

The following theorem gathers some well known properties of the distance function. It is a minor collection of the properties one can find in [F]-Theorem 4.8.

**Theorem 3.2.** *If $\emptyset \neq K \subsetneq \mathbb{R}^n$ is closed, then*

(i) *for all $x, y \in \mathbb{R}^n$ there holds $|d_K(x) - d_K(y)| \leq |x - y|$.*

(ii) *if $x \in \mathbb{R}^n \setminus K$, then $d_K$ is differentiable at $x$ if and only if $x \in K_1$.*

(iii) *if $d_K$ is differentiable at $x \in \mathbb{R}^n$, then $\nabla d_K(x) = (x - \xi_K(x))/d_K(x)$.*

**Remark 3.3.** Part (i) says that $d_K$ is Lipschitz continuous in $\mathbb{R}^n$. Thus by (ii) and Rademacher's Theorem (see [EvG]-§3.1.2-Theorem 2) we obtain $\mathcal{L}^n(\mathbb{R}^n \setminus K_1) = 0$. By this fact and (iii) we also get $|\nabla d_K| = 1$ a.e. in $\mathbb{R}^n$. □

The next property of $d_K$ can be obtained for example from [CS]-Proposition 2.2.2-(ii) together with Proposition 1.1.3-(c).

**Proposition 3.4.** *The function* $\alpha : \mathbb{R}^n \to \mathbb{R}$ *defined by* $\alpha(x) = C|x|^2/2 - d_K(x)$, *is convex in any open ball* $B \subset\subset \mathbb{R}^n \setminus K$, *for any* $C \geq 1/\operatorname{dist}(B, K)$.

**Proof.** First note that for all $a, b \in \mathbb{R}^n$ with $a \neq 0$ we have

$$|a + b| + |a - b| - 2|a| \leq \frac{|b|^2}{|a|}. \tag{3.1}$$

We choose an open ball $B \subset \mathbb{R}^n \setminus K$ with $r := \operatorname{dist}(B, K) > 0$ and take $x \in B$. Let $y \in K$ be such that $d_K(x) = |x - y|$. For any $z \in \mathbb{R}^n$ such that $x + z, x - z \in B$, we get

$$\begin{aligned}
\alpha(x + z) + \alpha(x - z) - 2\alpha(x) &= C|z|^2 - (d_K(x + z) + d_K(x - z) - 2d_K(x)) \\
&\geq C|z|^2 - (|x + z - y| + |x - z - y| - 2|x - y|) \\
\text{(by (3.1) for } a = x - y \text{ and } b = z) &\geq C|z|^2 - \frac{|z|^2}{|x - y|} \\
&\geq (C - 1/r)|z|^2.
\end{aligned}$$

Since $\alpha(x)$ is also continuous we obtain that $\alpha(x)$ is convex in $B$ for any $C \geq 1/r$ (see [CS]-Proposition A1.2). ■

Next we discuss sets with positive reach introduced by H. Federer in [F].

**Definition 3.5.** Let $\emptyset \neq K \subsetneq \mathbb{R}^n$ be closed. The *reach of a point* $x \in K$ is

$$\operatorname{Reach}(K, x) := \sup\{r \geq 0 : B_r(x) \subset K_1\}.$$

The *reach of the set* $K$ is

$$\operatorname{Reach}(K) := \inf_{x \in K} \operatorname{Reach}(K, x).$$

In general the complement of any domain of class $\mathcal{C}^2$ is a set of positive reach (see Lemma 3.10-(1)). The following example is taken from [KrP]-§6.2 and shows that the complement of a domain that is less smooth than $\mathcal{C}^2$ need not be a set of positive reach.

**Example 3.6.** Let $\Omega_\delta := \{(x_1, x_2) \in \mathbb{R}^2 : x_2 > |x_1|^{2-\delta}\}$, where $0 < \delta < 1$. Then $\Omega_\delta$ is of class $\mathcal{C}^{1, 2-\delta}$ and $\mathbb{R}^n \setminus \Omega_\delta$ is not a set of positive reach. More precisely, there exist $\epsilon_0 = \epsilon_0(\delta) > 0$ such that for any $0 < \epsilon < \epsilon_0$, every point $(0, \epsilon) \in \mathbb{R}^2$ has two closest points in $\mathbb{R}^n \setminus \Omega_\delta$. □

**Remark 3.7** ([F]). Motzkin's Theorem ([Hrm]-Theorem 2.1.30) asserts that a nonempty closed set $K$ in $\mathbb{R}^n$ is convex if and only if every point in $\mathbb{R}^n$ has a unique closest point in $K$. Hence, in terms of the above definition, *a nonempty closed set $K$ in $\mathbb{R}^n$ is convex if and only if* $\mathrm{Reach}(K) = \infty$. □

**Example 3.8.** For an example of a set with zero reach, denote first by $Q_r(x_0)$ an open cube ("square" if $n = 2$) having edge length $r > 0$ and center at $x_0 \in \mathbb{R}^n$. It is then easily verified that $\mathrm{Reach}(\mathbb{R}^n \setminus Q_r(x_0)) = 0$. □

The following folklore lemma illustrates the connection between the reach of the closure of an open set and the distributional Laplacian of the distance function to the complement of the set.

**Lemma 3.9.** *Let $\emptyset \neq \Omega \subsetneq \mathbb{R}^n; n \geq 2$, be open and for any $x \in \mathbb{R}^n$ set*

$$d(x) := \mathrm{dist}(x, \mathbb{R}^n \setminus \Omega) = \inf\{|x - y| : y \in \mathbb{R}^n \setminus \Omega\}.$$

*Then*

(i) *If $\mathrm{Reach}(\overline{\Omega}) = \infty$ then $-\Delta d \geq 0$ in the sense of distributions in $\Omega$.*

(ii) *If $n = 2$ then $\mathrm{Reach}(\overline{\Omega}) = \infty$ if and only if $-\Delta d \geq 0$ in the sense of distributions in $\Omega$.*

(iii) *Let $h := \mathrm{Reach}(\overline{\Omega}) < \infty$. Then*

$$(d + h)(-\Delta d) \geq -(n - 1) \text{ in the sense of distributions in } \Omega. \tag{3.2}$$

**Proof.** Parts (i) and (ii) follow from the basic results of [ArmK] and Remark 3.7. Thus we prove only (iii). In case $h = 0$ we have $\Omega_h \equiv \Omega$ and estimate (3.2) rests on the fact that the function $A : \mathbb{R}^n \to \mathbb{R}$ defined by $A(x) := |x|^2 - d^2(x)$ is convex. To see this we take $x \in \mathbb{R}^n$ and let $y \in \mathbb{R}^n \setminus \Omega$ be such that $d(x) = |x - y|$. For any $z \in \mathbb{R}^n$ we get

$$\begin{aligned} A(x + z) + A(x - z) - 2A(x) &= 2|z|^2 - (d^2(x + z) + d^2(x - z) - 2d^2(x)) \\ &\geq 2|z|^2 - (|x + z - y|^2 + |x - z - y|^2 - 2|x - y|^2) \\ &= 0. \end{aligned}$$

Since $A(x)$ is also continuous we obtain that $A(x)$ is convex (see [CS]-Proposition A1.2). It follows by [EvG]-§6.3-Theorem 2 that the distributional Laplacian of $A$ is a nonnegative Radon measure on $\mathbb{R}^n$. The result follows since we have $\Delta A = 2(n - 1 - d\Delta d)$ in the sense of distributions in $\Omega$.

In case $h > 0$, we set $\Omega_h := \{x \in \mathbb{R}^n : d_{\overline{\Omega}}(x) < h\}$, and as in the previous proof, the continuous function $A_h : \mathbb{R}^n \to \mathbb{R}$ defined by $A_h(x) = |x|^2 - d^2_{\Omega_h^c}(x)$ is convex. Hence the distributional Laplacian of $A_h$ is a nonnegative Radon measure on $\mathbb{R}^n$. The result follows since for every $x \in \Omega$ we have $d_{\Omega_h^c}(x) = d(x) + h$ (see [F]-Corollary 4.9), and thus $\Delta A_h = 2(n - 1 - (d + h)\Delta d) \geq 0$ in the sense of distributions in $\Omega$. ■

## 3.2  $-\triangle d$ and mean convex domains

Let $\emptyset \neq \Omega \subsetneq \mathbb{R}^n$ be open. In this section we consider the choice $K = \mathbb{R}^n \setminus \Omega$ in Definition 3.1. For brevity, as in Lemma 3.9 we write $d$ instead of $d_{\mathbb{R}^n \setminus \Omega}$. Thus

$$d(x) := \operatorname{dist}(x, \mathbb{R}^n \setminus \Omega) = \inf\{|x - y| : y \in \mathbb{R}^n \setminus \Omega\},$$

whenever $x \in \mathbb{R}^n$. Also we will denote by $\Sigma$ the set of points in $\Omega$ which have more than one closest point on $\mathbb{R}^n \setminus \Omega$ (thus in terms of Definition 3.1 we have $K_1 \cap \Omega = \Omega \setminus \Sigma$). If $x \in \Omega \setminus \Sigma$, then $\xi(x)$ will stand for its unique closest point on $\mathbb{R}^n \setminus \Omega$.

Before stating the main result of this subsection we gather some additional properties of the distance function to the boundary that will be in use. Recall that by Theorem 3.2 we have that $d$ is Lipschitz continuous on $\mathbb{R}^n$ and in particular $|\nabla d(x)| = 1$ a.e. in $\Omega$.

From now on $\Omega$ will be a domain, i.e. an open and connected subset of $\mathbb{R}^n$.

The following Lemma follows from Lemmas 14.16 and 14.17 in [GTr].

**Lemma 3.10.** *Let $\Omega \subset \mathbb{R}^n$ be a domain with boundary of class $\mathcal{C}^2$.*

(1) *If in addition $\Omega$ satisfies a uniform interior sphere condition, then there exists $\delta > 0$ such that $\tilde{\Omega}_\delta := \{x \in \overline{\Omega} : d(x) < \delta\} \subset \Omega \setminus \Sigma$ and $d \in C^2(\tilde{\Omega}_\delta)$.*

(2) $d \in C^2(\overline{\Omega} \setminus \overline{\Sigma})$ *and for any $x \in \overline{\Omega} \setminus \overline{\Sigma}$, in terms of a principal coordinate system at $\xi(x) \in \partial\Omega$, there holds*

$$\begin{aligned}
(i) & \quad \nabla d(x) = -\vec{\nu}(\xi(x)) = (0, ..., 0, 1) \\
(ii) & \quad 1 - \kappa_i(\xi(x))d(x) > 0 \text{ for all } i = 1, ..., n-1 \\
(iii) & \quad [D^2 d(x)] = \operatorname{diag}\left[\frac{-\kappa_1(\xi(x))}{1 - \kappa_1(\xi(x))d(x)}, ..., \frac{-\kappa_{n-1}(\xi(x))}{1 - \kappa_{n-1}(\xi(x))d(x)}, 0\right],
\end{aligned}$$

*where $\vec{\nu}(\xi(x))$ is the unit outer normal at $\xi(x) \in \partial\Omega$, and $\kappa_1(\xi(x)), ..., \kappa_{n-1}(\xi(x))$ are the principal curvatures of $\partial\Omega$ at the point $\xi(x) \in \partial\Omega$.*

**Remark 3.11.** Part (2) of the above Lemma is proved in [GTr] only in $\tilde{\Omega}_\delta$. However, it is also true for the largest open set contained in $\Omega \setminus \Sigma$, i.e. $\Omega \setminus \overline{\Sigma}$ (see for instance [CrM], [LN], [CC], [Grg]). □

By Remark 3.3 we have that $\mathcal{L}^n(\Sigma) = 0$. An important fact we will need is that domains with boundary of class $\mathcal{C}^2$ satisfy $\mathcal{L}^n(\overline{\Sigma}) = 0$. This is proved in [Mnn]-Errata-§5.2 (see also [CrM] where however only bounded domains are discussed). In [MM] the authors construct a $\mathcal{C}^{1,1}$ bounded convex domain with $\mathcal{L}^n(\overline{\Sigma}) > 0$.

To state our main theorem we denote by $\mathcal{H}(y) := \frac{1}{n-1}\sum_{i=1}^{n-1} \kappa_i(y)$ the mean curvature of $\partial\Omega$ at the point $y \in \partial\Omega$.

**Theorem 3.12** ([Ps1]). *Let $\Omega \subset \mathbb{R}^n$ be a domain with boundary of class $\mathcal{C}^2$ satisfying a uniform interior sphere condition. Then $\mu := (-\Delta d)\mathrm{d}x$ is a signed Radon measure on $\Omega$. Let $\mu = \mu_{ac} + \mu_s$ be the Lebesgue decomposition of $\mu$ with respect to $\mathcal{L}^n$, i.e. $\mu_{ac} \ll \mathcal{L}^n$ and $\mu_s \perp \mathcal{L}^n$. Then $\mu_s \geq 0$ in $\Omega$, and $\mu_{ac} \geq (n-1)\underline{\mathcal{H}}\mathrm{d}x$ a.e. in $\Omega$, where $\underline{\mathcal{H}} := \inf_{y \in \partial\Omega} \mathcal{H}(y)$.*

**Proof.** Letting $\delta$ be as in Lemma 3.10-(1) we set $\Omega_\delta = \{x \in \Omega : d(x) < \delta\}$. Then $-\Delta d$ is a continuous function on $\Omega_\delta$ and so $\mu^0 := (-\Delta d)\mathrm{d}x$ is a signed Radon measure on $\Omega_\delta$, absolutely continuous with respect to $\mathcal{L}^n$.

Next let $\{B_i\}_{i \geq 1}$ be a cover of the set $\Omega \setminus \Omega_\delta$ comprised of open balls $B_i$ for which $\mathrm{dist}(B_i, \partial\Omega) > \delta/2$ for all $i \geq 1$. According to Proposition 3.4 the function $\alpha(x) := |x|^2/\delta - d(x)$ is convex in each $B_i$. From [EvG]-§6.3-Theorem 2, we deduce that there exist nonnegative Radon measures $\{\nu^i\}_{i \geq 1}$, respectively on $\{B_i\}_{i \geq 1}$, such that

$$\int_{B_i} \phi \Delta \tilde{A} \mathrm{d}x = \int_{B_i} \phi \mathrm{d}\nu^i,$$

for all $\phi \in C_c^\infty(B_i)$. Since $\Delta \alpha = 2n/\delta - \Delta d$ in the sense of distributions, we get

$$\int_{B_i} \phi(-\Delta d)\mathrm{d}x = \int_{B_i} \phi \mathrm{d}\nu^i - \frac{2n}{\delta} \int_{B_i} \phi \mathrm{d}x, \tag{3.3}$$

for all $\phi \in C_c^\infty(B_i)$. Hence $\mu^i := (-\Delta d)\mathrm{d}x = \nu^i - \frac{2n}{\delta}\mathrm{d}x$ is a signed Radon measure on $B_i$. Let $\{\eta_i\}_{i \geq 1}$ be a $C^\infty$ partition of unity subordinated to the open covering $\{B_i\}_{i \geq 1}$ of $\Omega \setminus \Omega_\delta$, i.e.

$$\eta_i \in C_c^\infty(B_i), \quad 0 \leq \eta_i(x) \leq 1 \text{ in } B_i \quad \text{and} \quad \sum_{i=1}^\infty \eta_i(x) = 1 \text{ in } \Omega \setminus \Omega_\delta.$$

Further, for $x \in \Omega$ define $\eta_0(x) = 1 - \sum_{i=1}^\infty \eta_i(x)$. We then have

$$\mathrm{sprt}\, \eta_0 \subset \Omega_\delta, \quad \eta_0(x) = 1 \text{ in } \Omega_{\delta/2} \quad \text{and} \quad \sum_{i=0}^\infty \eta_i(x) = 1 \text{ in } \Omega.$$

We will now show that $\mu := \sum_{i=0}^\infty \eta_i \mu^i$ is a well defined signed Radon measure on $\Omega$, and $\mu = (-\Delta d)\mathrm{d}x$. To this end, for any $\phi \in C_c^\infty(\Omega)$ we have

$$\begin{aligned}
\int_\Omega \phi(-\Delta d)\mathrm{d}x &= \sum_{i=0}^\infty \int_\Omega \phi \eta_i(-\Delta d)\mathrm{d}x \\
\text{(by (3.3))} &= \int_\Omega \phi \eta_0 \mathrm{d}\mu^0 + \sum_{i=1}^\infty \left( \int_\Omega \phi \eta_i \mathrm{d}\nu^i - \frac{2n}{\delta} \int_\Omega \phi \eta_i \mathrm{d}x \right) \\
&= \int_\Omega \phi \eta_0 \mathrm{d}\mu^0 + \int_\Omega \phi \sum_{i=1}^\infty \eta_i \mathrm{d}\nu^i - \frac{2n}{\delta} \int_\Omega \phi \sum_{i=1}^\infty \eta_i \mathrm{d}x \\
&= \int_\Omega \phi \eta_0 \mathrm{d}\mu^0 + \int_\Omega \phi \sum_{i=1}^\infty \eta_i \mathrm{d}\mu^i \\
&= \int_\Omega \phi \mathrm{d}\mu,
\end{aligned}$$

where the middle equality follows since $\nu^i$ are positive Radon measures and thus $\sum_{i=0}^m \eta_i \nu^i$ is increasing in $m$ (see [EvG]-Section 1.9).

Next, by the Lebesgue Decomposition Theorem ([EvG]-§1.3-Theorem 3), $\mu = \mu_{ac} + \mu_s$ where

$$\mu_s = \sum_{i=0}^\infty \eta_i \mu_s^i = \sum_{i=1}^\infty \eta_i \mu_s^i = \sum_{i=1}^\infty \eta_i \nu_s^i \geq 0,$$

since $\mu^i = \nu^i - \frac{2n}{\delta} \mathrm{d}x$ and $\nu^i$ are nonnegative. Finally, from Lemma 3.10-(2) we get

$$\begin{aligned}
-\Delta d(x) &= \sum_{i=1}^{n-1} \frac{\kappa_i(\xi(x))}{1 - \kappa_i(\xi(x))d(x)} \\
&\geq \sum_{i=1}^{n-1} \kappa_i(\xi(x)) \\
&= (n-1)\mathcal{H}(\xi(x)) \\
&\geq (n-1)\underline{\mathcal{H}}, \quad \forall x \in \Omega \setminus \overline{\Sigma}.
\end{aligned}$$

Now by Lemma 3.10-(2), $-\Delta d$ is a continuous function on $\Omega \setminus \overline{\Sigma}$ and so

$$\mu_{ac} = (-\Delta d)\mathrm{d}x \geq (n-1)\underline{\mathcal{H}}\mathrm{d}x \text{ in } \Omega \setminus \overline{\Sigma}.$$

Recalling that $\mathcal{L}^n(\overline{\Sigma}) = 0$ when $\partial\Omega \in \mathcal{C}^2$ and since $\Omega = (\Omega \setminus \overline{\Sigma}) \cup \overline{\Sigma}$, we conclude $\mu_{ac} \geq (n-1)\underline{\mathcal{H}}\mathrm{d}x$ a.e. in $\Omega$. ∎

**Definition 3.13.** A domain $\Omega$ with boundary of class $\mathcal{C}^2$ is said to be *mean convex* if $\mathcal{H}(y) \geq 0$ for all $y \in \partial\Omega$.

Theorem 3.12 along with Lemma 3.10 provides us a characterization of mean convexity in terms of the distance function for sufficiently smooth domains. More precisely

**Corollary 3.14.** *Let $\Omega$ be a domain with boundary of class $\mathcal{C}^2$ satisfying a uniform interior sphere condition. Then $\Omega$ is mean convex if and only if $-\Delta d \geq 0$ in the sense of distributions in $\Omega$.*

**Remark 3.15.** The resulting lower bound $-(\Delta d)\mathrm{d}x \geq (n-1)\underline{\mathcal{H}}\mathrm{d}x$ is optimal. To see this assume first that $\Omega$ is bounded and choose a point $y_0 \in \partial\Omega$ such that $\mathcal{H}(y_0) = \underline{\mathcal{H}}$. Pick $0 \leq \phi_\delta \in C_c^\infty(\Omega)$, such that $\mathrm{sprt}\{\phi_\delta\} \subset B_\delta(y_0) \cap \Omega_\delta$; $\delta > 0$. For sufficiently small $\delta$, as $x \in \Omega_\delta$ approaches $y_0$ we have $-\Delta d(x) = (n-1)\underline{\mathcal{H}} + O(d(x))$. Thus, as $\delta \downarrow 0$ we have $-\Delta d(x) = (n-1)\underline{\mathcal{H}} + o_\delta(1)$ for all $x \in B_\delta(y_0) \cap \Omega_\delta$, and so

$$\begin{aligned}
\inf_{0 \leq \phi \in C_c^\infty(\Omega) \setminus \{0\}} \frac{\int_\Omega \phi(-\Delta d)\mathrm{d}x}{\int_\Omega \phi \mathrm{d}x} &\leq \frac{\int_{B_\delta(y_0) \cap \Omega_\delta} \phi_\delta(-\Delta d)\mathrm{d}x}{\int_{B_\delta(y_0) \cap \Omega_\delta} \phi_\delta \mathrm{d}x} \\
&= (n-1)\underline{\mathcal{H}} + o_\delta(1).
\end{aligned}$$

If $\Omega$ is unbounded we may consider a sequence $\{y_k\} \subset \partial\Omega$ converging to $y_0$, and repeat the above argument for any such point, to obtain

$$\inf_{0 \leq \phi \in C_c^\infty(\Omega)\setminus\{0\}} \frac{\int_\Omega \phi(-\Delta d)\mathrm{d}x}{\int_\Omega \phi \mathrm{d}x} \leq (n-1)\mathcal{H}(y_k) + o_\delta(1).$$

Since $\mathcal{H}(y)$ is a continuous function on $\partial\Omega$, we end up by letting $k \to \infty$. □

**Remark 3.16.** Assume that $\Omega$ is a domain with $\partial\Omega \in \mathcal{C}^2$ and such that it satisfies a uniform interior sphere condition. By Theorem 3.12 and (2.8) for $s > 1$, we obtain

$$\int_\Omega \frac{|\nabla u|^q}{d^{s-q}}\mathrm{d}x \geq \left(\frac{s-1}{q}\right)^q \int_\Omega \frac{|u|^q}{d^s}\mathrm{d}x + (n-1)\underline{\mathcal{H}}\left(\frac{s-1}{q}\right)^{q-1} \int_\Omega \frac{|u|^q}{d^{s-1}}\mathrm{d}x; \quad q \geq 1. \quad (3.4)$$

**Remark 3.17.** Corollary 3.14 was also noted in [LL] without proof. After [Ps1] was submitted, a proof not based on Theorem 3.12) was given in [LLL] where also (3.4) was obtained for $s = q > 1$. □

# Chapter 4

# $L^1$ Hardy inequalities with weights

In this chapter we prove Theorems I, III and IV, stated in the introduction. These results are the content of the work [Ps1]. We will also give improvements to the $L^1$-weighted Hardy inequality involving the $L^1$-norm of the length of the gradient. These improvements are sharp in the sense that extremal sets exist, i.e. sets in which the remainder term obtained is the optimal one.

## 4.1 Inequalities in sets without regularity assumptions on the boundary

As in §2.1.1 all inequalities of this chapter will follow by the integration by parts formula which we formalize as follows: let $\Omega$ be an open set in $\mathbb{R}^n$ and $\vec{T}$ be a vector field on $\Omega$. Integrating by parts and using elementary inequalities we get

$$\int_\Omega |\vec{T}||\nabla u|\mathrm{d}x \geq \int_\Omega \mathrm{div}(\vec{T})|u|\mathrm{d}x, \tag{4.1}$$

for all $u \in C_c^\infty(\Omega)$, where we have also used the fact that $|\nabla |u|| = |\nabla u|$ a.e. in $\Omega$. For example, setting $\vec{T}(x) = -(d(x))^{1-s}\nabla d(x)$ for a.e. $x \in \Omega$ we deduce (2.8) for $q = 1$:

**Lemma 4.1.** *Let $\Omega \subsetneq \mathbb{R}^n$ be open. For all $u \in C_c^\infty(\Omega)$ and all $s \geq 1$*

$$\int_\Omega \frac{|\nabla u|}{d^{s-1}}\mathrm{d}x \geq (s-1)\int_\Omega \frac{|u|}{d^s}\mathrm{d}x + \int_\Omega \frac{|u|}{d^{s-1}}(-\Delta d)\mathrm{d}x, \tag{4.2}$$

*where $-\Delta d$ is meant in the distributional sense. If $\Omega$ has finite volume then equality holds for $u_\varepsilon(x) = (d(x))^{s-1+\varepsilon} \in W_0^{1,1}(\Omega; d^{-(s-1)})$, where $\varepsilon > 0$.*

### 4.1.1 General sets

A covering of $\Omega$ by cubes was used in [Avkh] to prove the next Theorem. We present an elementary proof.

**Theorem 4.2.** *Let $\Omega \subsetneq \mathbb{R}^n$ be open. For all $u \in C_c^\infty(\Omega)$ and all $s > n$, there holds*

$$\int_\Omega \frac{|\nabla u|}{d^{s-1}} \mathrm{d}x \geq (s-n) \int_\Omega \frac{|u|}{d^s} \mathrm{d}x. \tag{4.3}$$

**Proof.** Coupling (4.2) and Lemma 3.9-(iii) with $h = 0$, we get

$$\int_\Omega \frac{|\nabla u|}{d^{s-1}} \mathrm{d}x \geq (s-1) \int_\Omega \frac{|u|}{d^s} \mathrm{d}x - (n-1) \int_\Omega \frac{|u|}{d^s} \mathrm{d}x = (s-n) \int_\Omega \frac{|u|}{d^s} \mathrm{d}x. \qquad\blacksquare$$

**Remark 4.3.** The constant obtained is just a lower bound for the best constant. The best constant differs from one open set to another. However $\mathbb{R}^n \setminus \{0\}$ serves as an extremal domain for Theorem 4.2. More precisely, letting $\Omega = \mathbb{R}^n \setminus \{0\}$ we have $d(x) = |x|$ and (4.3) reads as follows

$$\int_{\mathbb{R}^n} \frac{|\nabla u|}{|x|^{s-1}} \mathrm{d}x \geq (s-n) \int_{\mathbb{R}^n} \frac{|u|}{|x|^s} \mathrm{d}x; \quad s > n, \tag{4.4}$$

for all $u \in C_c^\infty(\mathbb{R}^n \setminus \{0\})$. To illustrate the optimality of the constant on the right hand side of (4.4), we define the following function

$$u_\delta(x) := \chi_{B_\eta(0) \setminus B_\delta(0)}(x); \quad x \in \mathbb{R}^n, \tag{4.5}$$

where $0 < \delta < \eta$ and $\eta$ is fixed. The distributional gradient of $u_\delta$ is $\nabla u_\delta = \vec{\nu}_{\partial B_\delta(0)} \delta_{\partial B_\delta(0)} - \vec{\nu}_{\partial B_\eta(0)} \delta_{\partial B_\eta(0)}$ where, for any $r > 0$, $\vec{\nu}_{\partial B_r(0)}$ stands for the outward pointing unit normal vector field along $\partial B_r(0) = \{x \in \mathbb{R}^n : |x| = r\}$, and by $\delta_{\partial B_r(0)}$ we denote the Dirac measure on $\partial B_r(0)$. Moreover the total variation of $\nabla u_\delta$ is $|\nabla u_\delta| = \delta_{\partial B_\delta(0)} + \delta_{\partial B_\eta(0)}$. Using the co-area formula, we get

$$\begin{aligned}
\frac{\int_{\mathbb{R}^n} \frac{|\nabla u_\delta|}{|x|^{s-1}} \mathrm{d}x}{\int_{\mathbb{R}^n} \frac{|u_\delta|}{|x|^s} \mathrm{d}x} &= \frac{\delta^{1-s} |\partial B_\delta(0)| + \eta^{1-s} |\partial B_\eta(0)|}{\int_\delta^\eta r^{-s} |\partial B_r(0)| \mathrm{d}r} \\
&= \frac{\delta^{n-s} + \eta^{n-s}}{\int_\delta^\eta r^{n-s-1} \mathrm{d}r} \\
&= (s-n) \frac{\delta^{n-s} + \eta^{n-s}}{\delta^{n-s} - \eta^{n-s}} \\
&\to s-n, \quad \text{as } \delta \downarrow 0. \qquad \square
\end{aligned} \tag{4.6}$$

**Remark 4.4.** Although not smooth, functions like $u_\delta$ defined in (4.5) belong to $BV(\mathbb{R}^n)$ (the space of functions of bounded variation in $\mathbb{R}^n$), and thus we can use a $C_c^\infty$ approximation so that the calculation above to hold in the limit (see for instance [EvG]-§5.2). We illustrate bellow the steps of such an approximation for this particular case. Similar to the bellow mollification arguments can be used in any such case that appears throughout this chapter, so that all the calculations we present using $BV$ functions are legitimate. In this Remark, whenever $r > 0$ we write $B_r$ instead of $B_r(0)$ to denote an open ball of radius $r$ with center at the origin. To begin with we set

$$\varphi(x) = \begin{cases} c\exp\left(\frac{1}{|x|^2-1}\right) & \text{if } |x| < 1 \\ 0 & \text{if } |x| \geq 1, \end{cases}$$

where the constant $c > 0$ is selected so that $\int_{\mathbb{R}^n} \varphi(x)\mathrm{d}x = 1$. Define the mollifier $\varphi_\varepsilon(x) = (1/\varepsilon^n)\varphi(x/\varepsilon)$, $0 < \varepsilon < \delta < \eta/3$. Note that $\mathrm{sprt}\{\varphi_\varepsilon\} = \overline{B}_\varepsilon$. We then set

$$u_{\delta,\varepsilon}(x) = (\varphi_\varepsilon * u_\delta)(x) := \int_{B_\eta \setminus B_\delta} \varphi_\varepsilon(x-y)\mathrm{d}y, \quad x \in \mathbb{R}^n. \tag{4.7}$$

The differentiation under the integral sign is permitted and so $u_{\delta,\varepsilon} \in C_c^\infty(\mathbb{R}^n)$, with $\mathrm{sprt}\{u_{\delta,\varepsilon}\} = \overline{B}_{\eta+\varepsilon} \setminus B_{\delta-\varepsilon}$. In addition

$$\begin{cases} u_{\delta,\varepsilon}(x) \equiv 1 & \text{in } B_{\eta-\varepsilon} \setminus B_{\delta+\varepsilon} \\ 0 \leq u_{\delta,\varepsilon}(x) \leq 1 & \text{in } (B_{\eta+\varepsilon} \setminus B_{\eta-\varepsilon}) \cup (B_{\delta+\varepsilon} \setminus B_{\delta-\varepsilon}). \end{cases}$$

As $\varepsilon \downarrow 0$, the function $u_{\delta,\varepsilon}$ approaches the characteristic function of the set $B_\eta \setminus B_\delta$. Thus we have

$$\int_{\mathbb{R}^n} \frac{|u_{\delta,\varepsilon}(x)|}{|x|^s}\mathrm{d}x = \int_{B_{\eta-\varepsilon}\setminus B_{\delta+\varepsilon}} \frac{1}{|x|^s}\mathrm{d}x + \int_{B_{\eta+\varepsilon}\setminus B_{\eta-\varepsilon}} \frac{u_{\delta,\varepsilon}(x)}{|x|^s}\mathrm{d}x + \int_{B_{\delta+\varepsilon}\setminus B_{\delta-\varepsilon}} \frac{u_{\delta,\varepsilon}(x)}{|x|^s}\mathrm{d}x$$

$$= \int_\delta^\eta r^{-s}|\partial B_r|\mathrm{d}r + o_\varepsilon(1). \tag{4.8}$$

On the other hand, changing variables and using the Gauss-Green Theorem

$$\nabla u_{\delta,\varepsilon}(x) = -\int_{B_\eta \setminus B_\delta} \nabla_y \varphi_\varepsilon(x-y)\mathrm{d}y$$

$$= \int_{\partial B_\delta} \varphi_\varepsilon(x-y)\vec{\nu}_{\partial B_\delta}(y)\mathrm{d}S_y - \int_{\partial B_\eta} \varphi_\varepsilon(x-y)\vec{\nu}_{\partial B_\eta}(y)\mathrm{d}S_y,$$

and so,

$$\int_{\mathbb{R}^n} \frac{|\nabla u_{\delta,\varepsilon}(x)|}{|x|^{s-1}}\mathrm{d}x \leq \int_{\partial B_\delta}\int_{B_{\eta+\varepsilon}\setminus B_{\delta-\varepsilon}} \frac{\varphi_\varepsilon(x-y)}{|x|^{s-1}}\mathrm{d}x\mathrm{d}S_y + \int_{\partial B_\eta}\int_{B_{\eta+\varepsilon}\setminus B_{\delta-\varepsilon}} \frac{\varphi_\varepsilon(x-y)}{|x|^{s-1}}\mathrm{d}x\mathrm{d}S_y$$

$$=: \mathcal{I}_\varepsilon + \mathcal{J}_\varepsilon,$$

where we have also reversed the order of integration. Since

$$\lim_{\varepsilon\downarrow 0}\int_{B_{\eta+\varepsilon}\setminus B_{\delta-\varepsilon}} \frac{\varphi_\varepsilon(x-y)}{|x|^{s-1}}\mathrm{d}x = \frac{1}{|y|^{s-1}},$$

we conclude $\mathcal{I}_\varepsilon \to \delta^{1-s}|\partial B_\delta|$ as $\varepsilon \downarrow 0$. Similarly $\mathcal{J}_\varepsilon \to \eta^{1-s}|\partial B_\eta|$ as $\varepsilon \downarrow 0$. From these last two facts we may write

$$\int_{\mathbb{R}^n} \frac{|\nabla u_{\delta,\varepsilon}(x)|}{|x|^{s-1}}\mathrm{d}x \leq \delta^{1-s}|\partial B_\delta| + \eta^{1-s}|\partial B_\eta| + o_\varepsilon(1). \tag{4.9}$$

Thus, by (4.8) and (4.9)

$$\frac{\int_{\mathbb{R}^n} \frac{|\nabla u_{\delta,\varepsilon}|}{|x|^{s-1}} \mathrm{d}x}{\int_{\mathbb{R}^n} \frac{|u_{\delta,\varepsilon}|}{|x|^s} \mathrm{d}x} \leq \frac{\delta^{1-s}|\partial B_\delta| + \eta^{1-s}|\partial B_\eta| + o_\varepsilon(1)}{\int_\delta^\eta r^{-s}|\partial B_r|\mathrm{d}r + o_\varepsilon(1)},$$

and letting $\varepsilon \downarrow 0$, the calculation following (4.6) can be repeated. □

**General sets with finite inner radius.** If we consider sets having finite inner radius we can improve (4.3) for $s \geq n$. First we prove the following

**Theorem 4.5.** *Let* $\Omega \subsetneq \mathbb{R}^n$ *be open and such that* $R := \sup_{x \in \Omega} d(x) < \infty$. *For all* $u \in C_c^\infty(\Omega)$, *all* $s \geq n, \gamma > 1$, *there holds*

$$\int_\Omega \frac{|\nabla u|}{d^{s-1}}\mathrm{d}x \geq (s-n)\int_\Omega \frac{|u|}{d^s}\mathrm{d}x + \frac{C}{R^{s-n}}\int_\Omega \frac{|u|}{d^n} X^\gamma\Big(\frac{d}{R}\Big)\mathrm{d}x, \qquad (4.10)$$

*where* $C \geq \gamma - 1$.

**Proof.** We set $\vec{T}(x) = -(d(x))^{1-s}[1 - (d(x)/R)^{s-n} X^{\gamma-1}(d(x)/R)]\nabla d(x)$ for a.e. $x \in \Omega$. Since $|1 - (d(x)/R)^{s-n} X^{\gamma-1}(d(x)/R)| \leq 1$ for all $x \in \Omega$, we have

$$\int_\Omega |\vec{T}||\nabla u|\mathrm{d}x \leq \int_\Omega \frac{|\nabla u|}{d^{s-1}}\mathrm{d}x.$$

Using the rule $\nabla X^{\gamma-1}(d(x)/R) = (\gamma-1)X^\gamma(d(x)/R)\frac{\nabla d(x)}{d(x)}$ for a.e. $x \in \Omega$, we compute

$$\begin{aligned}\mathrm{div}(\vec{T}) &= (s-1)d^{-s}[1-(d/R)^{s-n}X^{\gamma-1}(d/R)] + \frac{s-n}{R^{s-n}}d^{-n}X^{\gamma-1}(d/R) \\ &\quad + \frac{\gamma-1}{R^{s-n}}d^{-n}X^\gamma(d/R) + d^{1-s}[1-(d/R)^{s-n}X^{\gamma-1}(d/R)](-\Delta d).\end{aligned}$$

Since $1 - (d(x)/R)^{s-n}X^{\gamma-1}(d(x)/R) \geq 0$ for all $x \in \Omega$, we use Lemma 3.9-(iii) with $h=0$ on the last term of the above equality. After a straightforward computation

$$\mathrm{div}(\vec{T}) \geq (s-n)d^{-s} + \frac{\gamma-1}{R^{s-n}}d^{-n}X^\gamma(d/R).$$

This means that

$$\int_\Omega \mathrm{div}(\vec{T})|u|\mathrm{d}x \geq (s-n)\int_\Omega \frac{|u|}{d^s}\mathrm{d}x + \frac{\gamma-1}{R^{s-n}}\int_\Omega \frac{|u|}{d^n}X^\gamma(d/R)\mathrm{d}x.$$

and the result follows from (4.1). ■

**Remark 4.6.** A punctured domain serves as an extremal domain for Theorem 4.5. More precisely, let $\Omega = U \setminus \{0\}$, where $U$ is an open, connected subset of $\mathbb{R}^n$ containing the origin and satisfying $R := \sup_{x \in U} d(x) < \infty$. We define $u_\delta$ as in (4.5) where $\eta$ is fixed and sufficiently small such that $d(x) = |x|$ in $B_\eta$. For any $s \geq n$, we have

$$\frac{\int_\Omega \frac{|\nabla u_\delta|}{|x|^{s-1}}dx - (s-n)\int_\Omega \frac{|u_\delta|}{|x|^s}dx}{\int_\Omega \frac{|u_\delta|}{|x|^n}X(|x|/R)dx} = \frac{\delta^{1-s}|\partial B_\delta| + \eta^{1-s}|\partial B_\eta| - (s-n)\int_\delta^\eta r^{-s}|\partial B_r|dr}{\int_\delta^\eta r^{-n}X(r/R)|\partial B_r|dr}$$

$$= \frac{\delta^{n-s} + \eta^{n-s} - (s-n)\int_\delta^\eta r^{n-s-1}dr}{\int_\delta^\eta r^{-1}X(r/R)dr}$$

$$= \frac{2\eta^{n-s}}{\log\left(\frac{X(\eta/R)}{X(\delta/R)}\right)}$$

$$= o_\delta(1).$$

Thus, for a punctured domain, inequality (4.10) does not hold when $\gamma = 1$, as well as the exponent $n$ on the second term of the right hand side in (4.10) cannot be increased. □

The second improvement for sets having finite inner radius appears here for the first time. Further extension is under investigation.

**Theorem 4.7.** *Let $\Omega \subsetneq \mathbb{R}^n$ be open and such that $R := \sup_{x \in \Omega} d(x) < \infty$. For all $u \in C_c^\infty(\Omega)$ and all $s > n$, there holds*

$$\int_\Omega \frac{|\nabla u|}{d^{s-1}}dx - (s-n)\int_\Omega \frac{|u|}{d^s}dx \geq \frac{1}{R^{s-n}}\int_\Omega \frac{|\nabla u|}{d^{n-1}}dx. \tag{4.11}$$

**Proof.** We set $\vec{T}(x) = -(d(x))^{1-s}[1 - (d(x)/R)^{s-n}]\nabla d(x)$ for a.e. $x \in \Omega$. Since

$$|\vec{T}(x)| = (d(x))^{1-s}\left[1 - \left(\frac{d(x)}{R}\right)^{s-n}\right] \quad \text{a.e. } x \in \Omega,$$

we have

$$\int_\Omega |\vec{T}||\nabla u|dx = \int_\Omega \frac{|\nabla u|}{d^{s-1}}dx - \frac{1}{R^{s-n}}\int_\Omega \frac{|\nabla u|}{d^{n-1}}dx.$$

We also calculate

$$\operatorname{div}(\vec{T}) = (s-1)d^{-s}[1 - (d/R)^{s-n}] + \frac{s-n}{R^{s-n}}d^{-n} + d^{1-s}[1 - (d/R)^{s-n}](-\Delta d), \quad \text{in } \Omega,$$

in the distributional sense. Since $1 - (d(x)/R)^{s-n} \geq 0$ for all $x \in \Omega$, we may use Lemma 3.9-(iii) with $h = 0$ on the last term of the above equality and after a straightforward computation to obtain

$$\int_\Omega \operatorname{div}(\vec{T})|u|dx \geq (s-n)\int_\Omega \frac{|u|}{d^s}dx.$$

The result follows from (4.1). ∎

**Remark 4.8.** A punctured domain serves also as an extremal domain for Theorem 4.7. As before, letting $\Omega = U \setminus \{0\}$, where $U$ is an open, connected subset of $\mathbb{R}^n$ containing the origin and satisfying $R := \sup_{x \in U} d(x) < \infty$, we define $u_\delta$ as in (4.5) where $\eta$ is fixed and sufficiently small such that $d(x) = |x|$ in $B_\eta$. By the co-area formula, for any $\varepsilon \geq 0$ we have

$$\begin{aligned}
\frac{\int_\Omega \frac{|\nabla u_\delta|}{|x|^{s-1}} \mathrm{d}x - (s-n) \int_\Omega \frac{|u_\delta|}{|x|^s} \mathrm{d}x}{\int_\Omega \frac{|\nabla u_\delta|}{|x|^{n-1+\varepsilon}} \mathrm{d}x} &= \frac{\delta^{1-s}|\partial B_\delta| + \eta^{1-s}|\partial B_\eta| - (s-n) \int_\delta^\eta r^{-s}|\partial B_r| \mathrm{d}r}{\delta^{1-n-\varepsilon}|\partial B_\delta| + \eta^{1-n-\varepsilon}|\partial B_\eta|} \\
&= \frac{\delta^{n-s} + \eta^{n-s} + \int_\delta^\eta (r^{n-s})' \mathrm{d}r}{\delta^{-\varepsilon} + \eta^{-\varepsilon}} \\
&= \frac{2\eta^{n-s}}{\delta^{-\varepsilon} + \eta^{-\varepsilon}} \\
&= \begin{cases} o_\delta(1) & \text{if } \varepsilon > 0 \\ \eta^{n-s} & \text{if } \varepsilon = 0. \end{cases}
\end{aligned}$$

Note that if $\varepsilon = 0$ then $\eta^{n-s} \downarrow R^{n-s}$ as $\eta \uparrow R$. □

### 4.1.2 Sets satisfying property $(\mathfrak{C})$

In this subsection we assume that

$$-\Delta d \geq 0 \text{ in } \Omega, \text{ in the sense of distributions.} \tag{$\mathfrak{C}$}$$

This condition was first used in the context of Hardy inequalities in [BFT1-3] and has been used intensively in [FMzT1-3] and [FMschT]. As we have showed in §3, domains with sufficiently smooth boundary carrying condition $(\mathfrak{C})$ are characterized as domains with nonnegative mean curvature of their boundary. In this section we do not impose regularity on the boundary.

**Theorem 4.9.** *Let $\Omega \subsetneq \mathbb{R}^n$ be open and such that condition $(\mathfrak{C})$ holds. For all $u \in C_c^\infty(\Omega)$ and all $s > 1$, there holds*

$$\int_\Omega \frac{|\nabla u|}{d^{s-1}} \mathrm{d}x \geq (s-1) \int_\Omega \frac{|u|}{d^s} \mathrm{d}x. \tag{4.12}$$

*Moreover the constant appearing on the right hand side of (4.12) is sharp.*

**Proof.** Since $(\mathfrak{C})$ holds we may cancel the last term in (4.2) and (4.12) follows. To prove the sharpness of the constant we pick $y \in \partial\Omega$ and define the family of $W_0^{1,1}(\Omega; d^{-(s-1)})$ functions by $u_\varepsilon(x) := \phi(x)(d(x))^{s-1+\varepsilon}$, $\varepsilon > 0$, where $\phi \in C_c^\infty(B_\delta(y))$, $0 \leq \phi \leq 1$ and

$\phi \equiv 1$ in $B_{\delta/2}(y)$ for some small but fixed $\delta$. We have

$$\begin{aligned}\frac{\int_\Omega \frac{|\nabla u_\varepsilon|}{d^{s-1}}\mathrm{d}x}{\int_\Omega \frac{|u_\varepsilon|}{d^s}\mathrm{d}x} &\leq s-1+\varepsilon + \frac{\int_\Omega |\nabla\phi| d^\varepsilon \mathrm{d}x}{\int_\Omega \phi d^{-1+\varepsilon}\mathrm{d}x}\\ &\leq s-1+\varepsilon + \frac{C}{\int_{\Omega\cap B_{\delta/2}(y)} d^{-1+\varepsilon}\mathrm{d}x}\\ &\leq s-1+o_\varepsilon(1),\end{aligned}$$

where $C$ is some universal constant (not depending on $\varepsilon$). ∎

**Remark 4.10.** In view of Theorem 4.9 and Lemma 4.1 we see that if $\Omega$ is bounded and condition ($\mathfrak{C}$) holds, then all constants appearing in (4.2) are optimal. □

**Sets satisfying property** ($\mathfrak{C}$) **and having finite inner radius.** The counterpart to Theorem 4.5 result for sets with finite inner radius reads as follows

**Theorem 4.11.** *Let $\Omega \subsetneq \mathbb{R}^n$ be open and such that condition ($\mathfrak{C}$) holds. Suppose in addition that $R := \sup_{x\in\Omega} d(x) < \infty$. For all $u \in C_c^\infty(\Omega)$, all $s \geq 1, \gamma > 1$, there holds*

$$\int_\Omega \frac{|\nabla u|}{d^{s-1}}\mathrm{d}x \geq (s-1)\int_\Omega \frac{|u|}{d^s}\mathrm{d}x + \frac{C}{R^{s-1}}\int_\Omega \frac{|u|}{d}X^\gamma\Big(\frac{d}{R}\Big)\mathrm{d}x, \qquad (4.13)$$

*where $C \geq \gamma - 1$.*

**Proof.** We set $\vec{T}(x) = -(d(x))^{1-s}[1 - (d(x)/R)^{s-1}X^{\gamma-1}(d(x)/R)]\nabla d(x)$ for a.e. $x \in \Omega$. Since $|1 - (d(x)/R)^{s-1}X^{\gamma-1}(d(x)/R)| \leq 1$ for all $x \in \Omega$, we have

$$\int_\Omega |\vec{T}||\nabla u|\mathrm{d}x \leq \int_\Omega \frac{|\nabla u|}{d^{s-1}}\mathrm{d}x.$$

Using the rule $\nabla X^{\gamma-1}(d(x)/R) = (\gamma-1)X^\gamma(d(x)/R)\frac{\nabla d(x)}{d(x)}$ for a.e. $x \in \Omega$, by a straightforward calculation we arrive at

$$\begin{aligned}\int_\Omega \mathrm{div}(\vec{T})|u|\mathrm{d}x &= (s-1)\int_\Omega \frac{|u|}{d^s}\mathrm{d}x + \frac{\gamma-1}{R^{s-1}}\int_\Omega \frac{|u|}{d}X^\gamma(d/R)\mathrm{d}x\\ &+ \int_\Omega \frac{|u|}{d^{s-1}}[1-(d/R)^{s-1}X^{\gamma-1}(d/R)](-\Delta d)\mathrm{d}x.\end{aligned}$$

Since $1 - (d(x)/R)^{s-1}X^{\gamma-1}(d(x)/R) \geq 0$ for all $x \in \Omega$ and also ($\mathfrak{C}$) holds, we may cancel the last term and the result follows by (4.1). ∎

**Remark 4.12.** We prove in §4.2.1-Example 4.25 that an infinite strip is an extremal domain for Theorem 4.11. More precisely, if $\Omega = \{x = (x', x_n) : x' \in \mathbb{R}^{n-1}, 0 < x_n < 2R\}$ for some $R > 0$, then (4.13) fails for $\gamma = 1$ and thus the exponent 1 on the distance to the boundary in the remainder term of (4.13) cannot be increased. □

The counterpart of Theorem 4.7 reads as follows

**Theorem 4.13.** *Let* $\Omega \subsetneq \mathbb{R}^n$ *be open, satisfies condition* $(\mathfrak{C})$ *and* $R := \sup_{x \in \Omega} d(x) < \infty$. *Then for all* $u \in C_c^\infty(\Omega)$ *and all* $s > 1$

$$\int_\Omega \frac{|\nabla u|}{d^{s-1}} \mathrm{d}x \geq (s-1) \int_\Omega \frac{|u|}{d^s} \mathrm{d}x + \frac{1}{R^{s-1}} \int_\Omega |\nabla u| \mathrm{d}x. \tag{4.14}$$

**Proof.** We insert $\vec{T}(x) = -(d(x))^{1-s}[1 - (d(x)/R)^{s-1}]\nabla d(x)$; a.e. $x \in \Omega$, in (4.1). Since $(d(x)/R)^{s-1} \leq 1$ for all $x \in \Omega$ we have

$$\int_\Omega |\vec{T}||\nabla u| \mathrm{d}x = \int_\Omega \frac{|\nabla u|}{d^{s-1}} \mathrm{d}x - \frac{1}{R^{s-1}} \int_\Omega |\nabla u| \mathrm{d}x.$$

On the other hand

$$\begin{aligned}\int_\Omega \mathrm{div}(\vec{T})|u| \mathrm{d}x &= (s-1) \int_\Omega \frac{|u|}{d^s} \mathrm{d}x + \int_\Omega \frac{|u|}{d^{s-1}}(1 - (d/R)^{s-1})(-\Delta d)\mathrm{d}x \\ &\geq (s-1) \int_\Omega \frac{|u|}{d^s} \mathrm{d}x,\end{aligned}$$

where now we have used again the fact that $(d(x)/R)^{s-1} \leq 1$ for all $x \in \Omega$ and also $(\mathfrak{C})$. The result follows. ∎

An infinite strip is an extremal domain for (4.14), in the following sense

**Lemma 4.14.** *For fixed* $R > 0$, *set* $S := \{x = (x', x_n) : x' \in \mathbb{R}^n \setminus \{0\}, 0 < x_n < 2R\}$. *Suppose that for some nonnegative* $\alpha$ *and* $s \geq 1$, *there holds*

$$\mathcal{C} := \inf_{u \in C_c^\infty(S) \setminus \{0\}} \tilde{Q}[u] \geq C_0 > 0,$$

*where*

$$\tilde{Q}[u] := \frac{\int_S \frac{|\nabla u|}{d^{s-1}} \mathrm{d}x - (s-1) \int_S \frac{|u|}{d^s} \mathrm{d}x}{\int_S \frac{|\nabla u|}{d^\alpha} \mathrm{d}x}.$$

*Then* $\alpha = 0$.

**Proof.** For $s = 1$ it is obvious. Note also that it is enough to assume that $0 < \alpha < 1$. Let $s > 1$. Pick any $\phi \equiv \phi(x') \in C_c^1(\mathbb{R}^n \setminus \{0\})$ such that $\mathrm{sprt}\{\phi\} \subset B_1$, where $B_1$ is the $n-1$ dimensional open ball with radius $1$ centered at $0'$. Let $\delta > 0$ and set $\phi_\delta \equiv \phi_\delta(x') := \phi(\delta x')$. Let also $0 < \varepsilon < \eta \leq R$. We test $\mathcal{C}$ with $u_{\varepsilon,\delta}(x) := \chi_{(\varepsilon,\eta)}(x_n)\phi_\delta(x')$. First note that

$$\nabla u_{\varepsilon,\delta}(x) = (\chi_{(\varepsilon,\eta)}(x_n)\nabla_{x'}\phi_\delta(x'), (\delta(x_n - \varepsilon) - \delta(x_n - \eta))\phi_\delta(x')),$$

where $\nabla_{x'} = (\frac{\partial}{\partial x_1}, \frac{\partial}{\partial x_2}, ..., \frac{\partial}{\partial x_{n-1}})$. Thus

$$|\nabla u_{\varepsilon,\delta}(x)| = \chi_{(\varepsilon,\eta)}(x_n)|\nabla_{x'}\phi_\delta(x')| + (\delta(x_n - \varepsilon) + \delta(x_n - \eta))|\phi_\delta(x')|.$$

Since $\eta \leq R/2$ we may substitute $d(x)$ by $x_n$ in $\tilde{Q}[u]$, and so

$$\tilde{Q}[u_{\varepsilon,\delta}] = \frac{\int_S \frac{|\nabla u_{\varepsilon,\delta}|}{x_n^{s-1}}\mathrm{d}x - (s-1)\int_S \frac{|u_{\varepsilon,\delta}|}{x_n^s}\mathrm{d}x}{\int_S \frac{|\nabla u_{\varepsilon,\delta}|}{x_n^\alpha}\mathrm{d}x}$$

$$= \frac{\int_\varepsilon^\eta \int_{B_{1/\delta}} \frac{|\nabla_{x'}\phi_\delta|}{x_n^{s-1}}\mathrm{d}x'\mathrm{d}x_n + (\frac{1}{\varepsilon^{s-1}} + \frac{1}{\eta^{s-1}})\int_{B_{1/\delta}}|\phi_\delta|\mathrm{d}x' - (s-1)\int_\varepsilon^\eta \int_{B_{1/\delta}} \frac{|\phi_\delta|}{x_n^s}\mathrm{d}x'\mathrm{d}x_n}{\int_\varepsilon^\eta \int_{B_{1/\delta}} \frac{|\nabla_{x'}\phi_\delta|}{x_n^\alpha}\mathrm{d}x'\mathrm{d}x_n + (\frac{1}{\varepsilon^\alpha} + \frac{1}{\eta^\alpha})\int_{B_{1/\delta}}|\phi_\delta|\mathrm{d}x'}$$

$$= \frac{K_\delta \int_\varepsilon^\eta x_n^{1-s}\mathrm{d}x_n + M_\delta(\varepsilon^{1-s} + \eta^{1-s}) + M_\delta \int_\varepsilon^\eta (x_n^{1-s})'\mathrm{d}x_n}{K_\delta \int_\varepsilon^\eta x_n^{-\alpha}\mathrm{d}x_n + M_\delta(\varepsilon^{-\alpha} + \eta^{-\alpha})},$$

where we have set $K_\delta := \int_{B_{1/\delta}} |\nabla_{x'}\phi_\delta(x')|\mathrm{d}x'$ and $M_\delta := \int_{B_{1/\delta}} |\phi_\delta(x')|\mathrm{d}x'$. Performing the integration appeared in the last term of the numerator we arrive at

$$\tilde{Q}[u_{\varepsilon,\delta}] = \frac{K_\delta \int_\varepsilon^\eta x_n^{1-s}\mathrm{d}x_n + 2M_\delta \eta^{1-s}}{K_\delta \int_\varepsilon^\eta x_n^{-\alpha}\mathrm{d}x_n + M_\delta(\varepsilon^{-\alpha} + \eta^{-\alpha})}.$$

By the change of variables $y' = \delta x'$ we obtain $K_\delta = \delta^{2-n}K_1$ and $M_\delta = \delta^{1-n}M_1$. Thus

$$\tilde{Q}[u_{\varepsilon,\delta}] = \frac{\delta^{2-n}K_1 \int_\varepsilon^\eta x_n^{1-s}\mathrm{d}x_n + 2\delta^{1-n}M_1\eta^{1-s}}{\delta^{2-n}K_1 \int_\varepsilon^\eta x_n^{-\alpha}\mathrm{d}x_n + \delta^{1-n}M_1(\varepsilon^{-\alpha} + \eta^{-\alpha})}$$

$$= \frac{\delta K_1 \int_\varepsilon^\eta x_n^{1-s}\mathrm{d}x_n + 2M_1\eta^{1-s}}{\frac{\delta K_1}{1-\alpha}(\eta^{1-\alpha} - \varepsilon^{1-\alpha}) + M_1(\varepsilon^{-\alpha} + \eta^{-\alpha})}.$$

To proceed we distinguish cases:
  • Let $1 < s < 2$. Then

$$\tilde{Q}[u_{\varepsilon,\delta}] = \frac{\frac{\delta K_1}{2-s}(\eta^{2-s} - \varepsilon^{2-s}) + 2M_1\eta^{1-s}}{\frac{\delta K_1}{1-\alpha}(\eta^{1-\alpha} - \varepsilon^{1-\alpha}) + M_1(\varepsilon^{-\alpha} + \eta^{-\alpha})}$$

$$= o_\varepsilon(1).$$

  • Now let $s = 2$. Then

$$\tilde{Q}[u_{\varepsilon,\delta}] = \frac{\delta K_1 \log(\eta/\varepsilon) + 2M_1\eta^{-1}}{\frac{\delta K_1}{1-\alpha}(\eta^{1-\alpha} - \varepsilon^{1-\alpha}) + M_1(\varepsilon^{-\alpha} + \eta^{-\alpha})}$$

$$= o_\varepsilon(1).$$

  • Finally let $s > 2$. Then

$$\tilde{Q}[u_{\varepsilon,\delta}] = \frac{\frac{\delta K_1}{s-2}(\varepsilon^{2-s} - \eta^{2-s}) + 2M_1\eta^{1-s}}{\frac{\delta K_1}{1-\alpha}(\eta^{1-\alpha} - \varepsilon^{1-\alpha}) + M_1(\varepsilon^{-\alpha} + \eta^{-\alpha})}.$$

We may set $\delta = \varepsilon^{s-2}$ so that $\tilde{Q}[u_{\varepsilon,\delta}] = o_\varepsilon(1)$. ∎

### 4.1.3 Sets with positive reach

In this subsection we obtain an interpolation inequality between (4.3) and (4.12) via sets with positive reach.

**Theorem 4.15.** *Let $\Omega \subsetneq \mathbb{R}^n$ be open and set $h := \text{Reach}(\overline{\Omega})$. Suppose in addition that $R := \sup_{x \in \Omega} d(x) < \infty$. For all $u \in C_c^\infty(\Omega)$ and all $s > \frac{h+nR}{h+R}$, there holds*

$$\int_\Omega \frac{|\nabla u|}{d^{s-1}} \mathrm{d}x \geq \left((s-1)\frac{h}{h+R} + (s-n)\frac{R}{h+R}\right) \int_\Omega \frac{|u|}{d^s} \mathrm{d}x. \tag{4.15}$$

**Proof.** Inserting (3.2) to (4.2) we obtain

$$\begin{aligned}
\int_\Omega \frac{|\nabla u|}{d^{s-1}} \mathrm{d}x &\geq (s-1) \int_\Omega \frac{|u|}{d^s} \mathrm{d}x - (n-1) \int_\Omega \frac{|u|}{d^s} \frac{d}{h+d} \mathrm{d}x. \\
&= \int_\Omega \frac{(s-1)h + (s-n)d}{h+d} \frac{|u|}{d^s} \mathrm{d}x \\
&\geq \frac{(s-1)h + (s-n)R}{h+R} \int_\Omega \frac{|u|}{d^s} \mathrm{d}x,
\end{aligned}$$

where the last inequality follows since $R < \infty$ and $\frac{(s-1)h+(s-n)d}{h+d}$ is decreasing in $d$. ∎

Note that this inequality interpolates between the case of a general open set $\Omega \subsetneq \mathbb{R}^n$ where we have $h = 0$ and the constant becomes $s - n$, and the case of a convex set $\Omega$ where $h = \infty$ and the constant becomes $s - 1$.

## 4.2 Proof of Theorem I

Let $\Omega$ be a domain satisfying property $(\mathfrak{C})$. We define the quotient

$$Q_\beta[u] := \frac{\int_\Omega \frac{|\nabla u|}{d^{s-1}} \mathrm{d}x - (s-1) \int_\Omega \frac{|u|}{d^s} \mathrm{d}x}{\int_\Omega \frac{|u|}{d^{s-\beta}} \mathrm{d}x}; \quad s > 1, \tag{4.16}$$

and we consider the following minimization problem

$$\mathcal{B}_\beta(\Omega) := \inf\{Q_\beta[u] : u \in C_c^\infty(\Omega) \setminus \{0\}\}; \quad 0 < \beta \leq s - 1.$$

The next Proposition shows that the essential range for $\beta$ is smaller.

**Proposition 4.16.** *Let $\Omega$ be a domain with boundary of class $\mathcal{C}^2$ satisfying property $(\mathfrak{C})$. If $s \geq 2$ then $\mathcal{B}_\beta(\Omega) = 0$ for all $0 < \beta < 1$. If $1 < s < 2$ then $\mathcal{B}_\beta(\Omega) = 0$ for all $0 < \beta \leq s-1$.*

**Proof.** For small $\delta > 0$, let $\Omega_\delta := \{x \in \Omega : d(x) < \delta\}$ and $\Omega_\delta^c = \Omega \setminus \Omega_\delta$. We test (4.16) with $u_\delta(x) = \chi_{\Omega_\delta^c}(x)\phi(x)$, where $\phi \in C_c^\infty(B_\varepsilon(y_0))$ for a fixed $y_0 \in \partial\Omega$ and sufficiently small $\varepsilon$, satisfying $\varepsilon > 3\delta$. We may suppose in addition that $0 \leq \phi \leq 1$ in $B_\varepsilon(y_0)$, $\phi \equiv 1$ in $B_{\varepsilon/2}(y_0)$ and $|\nabla\phi| \leq 1/\varepsilon$. This function is not in $C_c^\infty(\Omega)$ but since it is in $BV(\Omega)$ we can mollify the characteristic function so that the calculations below to hold in the limit. The distributional gradient of $u_\delta$ is $\nabla u_\delta = \chi_{\Omega_\delta^c}\nabla\phi - \vec{\nu}\phi\delta_{\partial\Omega_\delta^c}$, where $\vec{\nu}$ is the outward pointing unit normal vector field along $\partial\Omega_\delta^c$, and $\delta_{\partial\Omega_\delta^c}$ is the Dirac measure on $\partial\Omega_\delta^c$. Moreover the total variation of $\nabla u_\delta$ is $|\nabla u_\delta| = \chi_{\Omega_\delta^c}|\nabla\phi| + \phi\delta_{\partial\Omega_\delta^c}$. Since $\partial\Omega_\delta^c = \{x \in \Omega : d(x) = \delta\}$, we obtain

$$Q_\beta[u_\delta] = \frac{\int_{\Omega_\delta^c} |\nabla\phi| d^{1-s} \mathrm{d}x + \delta^{1-s} \int_{\partial\Omega_\delta^c} \phi \mathrm{d}S_x - (s-1) \int_{\Omega_\delta^c} \phi d^{-s} \mathrm{d}x}{\int_{\Omega_\delta^c} \phi d^{\beta-s} \mathrm{d}x}. \tag{4.17}$$

Using the fact that $|\nabla d(x)| = 1$ for a.e. $x \in \Omega$, we may perform an integration by parts in the last term of the numerator as follows

$$\begin{aligned}(s-1)\int_{\Omega_\delta^c} \phi d^{-s} \mathrm{d}x &= -\int_{\Omega_\delta^c} \phi \nabla d \cdot \nabla d^{1-s} \mathrm{d}x \\ &= \int_{\Omega_\delta^c} [\nabla\phi \cdot \nabla d] d^{1-s} \mathrm{d}x + \int_{\Omega_\delta^c} \phi d^{1-s} \Delta d \mathrm{d}x - \delta^{1-s} \int_{\partial\Omega_\delta^c} \phi \nabla d \cdot \vec{\nu} \mathrm{d}S_x.\end{aligned}$$

Since $\nabla d$ is the inner unit normal to $\partial\Omega$ we have $\nabla d \cdot \vec{\nu} = -1$, and substituting the above equality in (4.17) the surface integrals will be canceled to get

$$Q_\beta[u_\delta] = \frac{\int_{\Omega_\delta^c} [|\nabla\phi| - \nabla\phi \cdot \nabla d] d^{1-s} \mathrm{d}x + \int_{\Omega_\delta^c} \phi d^{1-s}(-\Delta d) \mathrm{d}x}{\int_{\Omega_\delta^c} \phi d^{\beta-s} \mathrm{d}x}.$$

By the fact that $-\Delta d(x) \leq c$ for all $x \in \Omega_\delta^c \cap B_\varepsilon$ and by the properties we imposed on $\phi$, we get

$$\begin{aligned}Q_\beta[u_\delta] &\leq \frac{\frac{2}{\varepsilon}\int_{\Omega_\delta^c \cap B_\varepsilon} d^{1-s} \mathrm{d}x + c\int_{\Omega_\delta^c \cap B_\varepsilon} d^{1-s} \mathrm{d}x}{\int_{\Omega_\delta^c \cap B_{\varepsilon/2}} d^{\beta-s} \mathrm{d}x} \\ &= c(\varepsilon)\frac{\int_{\Omega_\delta^c \cap B_\varepsilon} d^{1-s} \mathrm{d}x}{\int_{\Omega_\delta^c \cap B_{\varepsilon/2}} d^{\beta-s} \mathrm{d}x} \\ &=: c(\varepsilon)\frac{N(\delta)}{D(\delta)}.\end{aligned}$$

Using now the co-area formula we compute

$$\begin{aligned}N(\delta) &= \int_\delta^\varepsilon r^{1-s} \int_{\{x \in \Omega_\delta^c \cap B_\varepsilon : d(x) = r\}} \mathrm{d}S_x \mathrm{d}r \\ &\leq c_1(\varepsilon) \int_\delta^\varepsilon r^{1-s} \mathrm{d}r,\end{aligned}$$

where $c_1(\varepsilon) = \max_{r \in [0,\varepsilon]} |\{x \in \Omega_\delta^c \cap B_\varepsilon : d(x) = r\}|$. Also,

$$\begin{aligned} D(\delta) &= \int_\delta^{\varepsilon/2} r^{\beta-s} \int_{\{x \in \Omega_\delta^c \cap B_{\varepsilon/2} : d(x) = r\}} \mathrm{d}S_x \mathrm{d}r \\ &\geq \int_\delta^{\varepsilon/3} r^{\beta-s} \int_{\{x \in \Omega_\delta^c \cap B_{\varepsilon/2} : d(x) = r\}} \mathrm{d}S_x \mathrm{d}r \\ &\geq c_2(\varepsilon) \int_\delta^{\varepsilon/3} r^{\beta-s} \mathrm{d}r, \end{aligned}$$

where $c_2(\varepsilon) = \min_{r \in [0,\varepsilon/3]} |\{x \in \Omega_\delta^c \cap B_{\varepsilon/2} : d(x) = r\}|$. A direct computation reveals that if $s \geq 2$ then $Q_\beta[u_\delta] \leq o_\delta(1)$ for all $0 < \beta < 1$, and also if $1 < s < 2$ then $Q_\beta[u_\delta] \leq o_\delta(1)$ for all $0 < \beta \leq s - 1$. ∎

### 4.2.1 Lower and upper estimates for $\mathcal{B}_1(\Omega)$

In this subsection we obtain upper and lower estimates for $\mathcal{B}_1(\Omega)$. In particular we prove Theorem I and the optimality in Theorem IV of the introduction.

**Theorem 4.17** (**Lower estimate**). *Let $\Omega$ be a domain with boundary of class $\mathcal{C}^2$ satisfying a uniform interior sphere condition. If $s \geq 1$ then*

$$\mathcal{B}_1(\Omega) \geq (n-1)\underline{\mathcal{H}}, \tag{4.18}$$

*where $\underline{\mathcal{H}}$ is the infimum of the mean curvature of $\partial\Omega$.*

**Proof.** The estimate follows directly from (4.2) using Theorem 3.12. ∎

**Remark 4.18.** By Theorem 4.9, if condition ($\mathfrak{C}$) is satisfied then the first term in (4.2) is sharp. The passage from (4.2) to inequality (4.18) via Theorem 3.12 is also sharp, i.e. the constant $(n-1)\underline{\mathcal{H}}$ in the inequality

$$\int_\Omega \frac{|u|}{d^{s-1}}(-\Delta d)\mathrm{d}x \geq (n-1)\underline{\mathcal{H}} \int_\Omega \frac{|u|}{d^{s-1}}\mathrm{d}x, \quad \forall u \in C_c^\infty(\Omega),$$

is optimal. To see this set $v = d^{1-s}|u|$, to get

$$\begin{aligned} \inf_{u \in C_c^\infty(\Omega) \setminus \{0\}} \frac{\int_\Omega \frac{|u|}{d^{s-1}}(-\Delta d)\mathrm{d}x}{\int_\Omega \frac{|u|}{d^{s-1}}\mathrm{d}x} &\leq \inf_{0 \leq v \in C_c^\infty(\Omega) \setminus \{0\}} \frac{\int_\Omega v(-\Delta d)\mathrm{d}x}{\int_\Omega v \mathrm{d}x} \\ &\leq (n-1)\underline{\mathcal{H}} + o_\delta(1), \end{aligned}$$

by Remark 3.15. □

We next present upper bounds. We begin with an upper bound which although not sharp enough for our problem it is of independent interest.

**Definition 4.19.** The *Cheeger constant* $h(\Omega)$ of a bounded domain $\Omega$ with piecewise $\mathcal{C}^1$ boundary, is defined by $h(\Omega) := \inf_\omega \frac{|\partial \omega|}{|\omega|}$, where the infimum is taken over all subdomains $\omega \subset\subset \Omega$ with piecewise $\mathcal{C}^1$ boundary.

For existence of minimizers, uniqueness and regularity results concerning the Cheeger constant, we refer to [FrK] and references therein (especially in [StrZ]). See also [P] for an up to date survey.

**Proposition 4.20.** *Let $\Omega$ be a bounded domain with piecewise $\mathcal{C}^1$ boundary such that condition ($\mathfrak{C}$) holds. For all $s \geq 1$ we have $\mathcal{B}_1(\Omega) \leq h(\Omega)$.*

**Proof.** Take $\omega \subset\subset \Omega$ with piecewise $\mathcal{C}^1$ boundary and let $u_\omega(x) = (d(x))^{s-1}\chi_\omega(x)$. The distributional gradient and the total variation of this $BV(\Omega)$ function are respectively $\nabla u_\omega = (s-1)d^{s-2}\chi_\omega \nabla d - \vec{\nu} d^{s-1}\delta_{\partial \omega}$ and $|\nabla u_\omega| = (s-1)d^{s-2}\chi_\omega + d^{s-1}\delta_{\partial\omega}$, where $\vec{\nu}$ is the outward pointing unit normal vector field along $\partial\omega$, and $\delta_{\partial\omega}$ is the uniform Dirac measure on $\partial\omega$. We test (4.16) with $u_\omega$ to get

$$Q_1[u_\omega] = \frac{(s-1)\int_\omega d^{-1}\mathrm{d}x + \int_{\partial\omega}\mathrm{d}S_x - (s-1)\int_\omega d^{-1}\mathrm{d}x}{\int_\omega \mathrm{d}x} = \frac{|\partial\omega|}{|\omega|}.$$

In particular $h(\Omega) = \inf_\omega Q_1[u_\omega]$. By the standard $C_c^\infty$ approximation of the characteristic function of the domain $\Omega$ we obtain $\mathcal{B}_1(\Omega) \leq \frac{|\partial\omega|}{|\omega|}$ and thus $\mathcal{B}_1(\Omega) \leq h(\Omega)$. ∎

From Theorem 4.17 and Proposition 4.20 for $s = 1$, we conclude

**Corollary 4.21.** *If $\Omega$ is a strictly mean convex, bounded domain with boundary of class $\mathcal{C}^2$, there holds $h(\Omega) \geq (n-1)\underline{\mathcal{H}}$.*

Note that in [AltC] it is proved that if a bounded convex domain $\Omega$ is a self-minimizer of $h(\Omega)$, then it belongs to the class $\mathcal{C}^{1,1}$ and also the stronger estimate $h(\Omega) \geq (n-1)\overline{\mathcal{H}}$ holds. Here $\overline{\mathcal{H}}$ is the essential supremum of the mean curvature of the boundary (the last being defined in the almost everywhere sense since $\partial\Omega \in \mathcal{C}^{1,1}$).

The following result states a more useful upper bound for $\mathcal{B}_1(\Omega)$. It will be combined with Theorem 4.17 to give the best possible constant for special geometries.

**Theorem 4.22.** *Let $\Omega$ be a domain with boundary of class $\mathcal{C}^2$ satisfying a uniform interior sphere condition. If $s \geq 2$ then for all $\phi \in C_c^1(\partial\Omega)$,*

$$\mathcal{B}_1(\Omega) \leq (n-1)\frac{\int_{\partial\Omega}|\phi(y)|\mathcal{H}(y)\mathrm{d}S}{\int_{\partial\Omega}|\phi(y)|\mathrm{d}S} + \frac{\int_{\partial\Omega}|\nabla\phi(y)|\mathrm{d}S}{\int_{\partial\Omega}|\phi(y)|\mathrm{d}S},$$

*where $\mathcal{H}(y)$ is the mean curvature at the point $y \in \partial\Omega$.*

**Proof.** Let $\delta > 0$ such that for all $x \in \tilde{\Omega}_\delta := \{x \in \overline{\Omega} : d(x) < \delta\}$ there exists a unique point

$$\xi \equiv \xi(x) = x - d(x)\nabla d(x) \in \partial\Omega, \tag{4.19}$$

with $d(x) = |x - \xi|$. For any $t \in [0, \delta]$ the surface area element of $\partial\Omega_t^c = \{x \in \Omega : d(x) = t\}$ is given by

$$dS_t = (1 - \kappa_1 t)...(1 - \kappa_{n-1} t)dS = (1 - (n-1)t\mathcal{H} + O(t^2))dS, \tag{4.20}$$

where $\kappa_1, ..., \kappa_{n-1}$, are the principal curvatures of $\partial\Omega$, $dS$ is the surface area element of $\partial\Omega$ and $\mathcal{H}$ is the mean curvature of $\partial\Omega$ (see [S]-§13.5 & 13.6). Now let $0 < \varepsilon < \delta$ and chose $\phi \in C_c^1(\partial\Omega)$. We test (4.16) with $u_\varepsilon(x) = \chi_{\Omega_\varepsilon^c \setminus \Omega_\delta^c}(x)\phi(\xi(x))$, $\xi(x)$ as in (4.19), and then we will check the limit as $\varepsilon \downarrow 0$. The distributional gradient of $u_\varepsilon$, is $\nabla u_\varepsilon = (\vec{\nu}_\delta \delta_{\partial\Omega_\delta^c} - \vec{\nu}_\varepsilon \delta_{\partial\Omega_\varepsilon^c})\phi(\xi) + \chi_{\Omega_\varepsilon^c \setminus \Omega_\delta^c}\nabla_x \phi(\xi)$, where $\vec{\nu}_\delta, \vec{\nu}_\varepsilon$ are respectively, the outward pointing unit normal vector fields along $\partial\Omega_\delta^c, \partial\Omega_\varepsilon^c$. Its total variation is $|\nabla u_\varepsilon| = (\delta_{\partial\Omega_\delta^c} + \delta_{\partial\Omega_\varepsilon^c})|\phi(\xi)| + \chi_{\Omega_\varepsilon^c \setminus \Omega_\delta^c}|\nabla_x\phi(\xi)|$. Thus

$$\int_\Omega \frac{|\nabla u_\varepsilon|}{d^{s-1}}dx = \delta^{1-s}\int_{\partial\Omega_\delta^c}|\phi(\xi)|dS_\delta + \varepsilon^{1-s}\int_{\partial\Omega_\varepsilon^c}|\phi(\xi)|dS_\varepsilon + \int_{\Omega_\varepsilon^c \setminus \Omega_\delta^c}\frac{|\nabla_x\phi(\xi)|}{d^{s-1}}dx. \tag{4.21}$$

The first integral on the right-hand side of (4.21) is a constant since we will keep $\delta$ fixed. We perform the change of variables $y = \xi(x)$ in the second integral. Using (4.20) we have

$$\begin{aligned}\varepsilon^{1-s}\int_{\partial\Omega_\varepsilon^c}|\phi(\xi)|dS_\varepsilon &= \varepsilon^{1-s}\int_{\partial\Omega}|\phi(y)|(1-(n-1)\varepsilon\mathcal{H}(y) + O(\varepsilon^2))dS \\ &= \varepsilon^{1-s}M - (n-1)\varepsilon^{2-s}M_\mathcal{H} + O(\varepsilon^{3-s}),\end{aligned} \tag{4.22}$$

where $M := \int_{\partial\Omega}|\phi|dS$ and $M_\mathcal{H} := \int_{\partial\Omega}|\phi|\mathcal{H}dS$. Using the co-area formula the third term on the right-hand side of (4.21) is written as follows

$$\int_{\Omega_\varepsilon^c \setminus \Omega_\delta^c}\frac{|\nabla_x\phi(\xi)|}{d^{s-1}}dx = \int_\varepsilon^\delta t^{1-s}\int_{\partial\Omega_t^c}|\nabla_x\phi(\xi)|dS_t dt. \tag{4.23}$$

From (4.19) we have $\xi_i(x) = x_i - d(x)\frac{\partial}{\partial x_i}(d(x))$ and thus by Lemma 3.10-(c) we compute

$$\begin{aligned}\nabla_x\phi(\xi) &= \Big(\sum_{i=1}^n \phi_{\xi_i}(\xi)\frac{\partial\xi_i}{\partial x_1}, ..., \sum_{i=1}^n \phi_{\xi_i}(\xi)\frac{\partial\xi_i}{\partial x_n}\Big) \\ &= \Big(\frac{\phi_{\xi_1}(\xi)}{1-\kappa_1 d}, ..., \frac{\phi_{\xi_{n-1}}(\xi)}{1-\kappa_{n-1}d}, 0\Big).\end{aligned}$$

Thus (4.23) becomes

$$\int_{\Omega_\varepsilon^c \setminus \Omega_\delta^c} \frac{|\nabla_x \phi(\xi)|}{d^{s-1}} \mathrm{d}x = \int_\varepsilon^\delta t^{1-s} \int_{\partial\Omega} \left( \sum_{i=1}^{n-1} \left( \frac{\phi_{y_i}}{1-\kappa_i t} \right)^2 \right)^{1/2} \mathrm{d}S_t \mathrm{d}t$$

$$= \int_\varepsilon^\delta t^{1-s} \int_{\partial\Omega} \left( \sum_{i=1}^{n-1} \left( \phi_{y_i} \prod_{j=1, j\neq i}^{n-1} (1-\kappa_j t) \right)^2 \right)^{1/2} \mathrm{d}S \mathrm{d}t,$$

where we have changed variables by $y = \xi(x)$ in the last inequality. Expanding the product as in (4.20) we get

$$\int_{\Omega_\varepsilon^c \setminus \Omega_\delta^c} \frac{|\nabla_x \phi(\xi)|}{d^{s-1}} \mathrm{d}x \leq \int_\varepsilon^\delta t^{1-s} \int_{\partial\Omega} \left( \sum_{i=1}^{n-1} \phi_{y_i}^2 \left( 1 - [(n-1)\mathcal{H} - \kappa_i]t + c_1 t^2 \right)^2 \right)^{1/2} \mathrm{d}S \mathrm{d}t$$

$$\leq K \int_\varepsilon^\delta t^{1-s} \mathrm{d}t + c_2 \int_\varepsilon^\delta t^{2-s} \mathrm{d}t, \qquad (4.24)$$

for some $c_1, c_2 \geq 0$, where $K := \int_{\partial\Omega} |\nabla\phi| \mathrm{d}S$. Next, using co-area formula and the same change of variables we get

$$(s-1) \int_\Omega \frac{|u_\varepsilon|}{d^s} \mathrm{d}x = (s-1) \int_\varepsilon^\delta t^{-s} \int_{\partial\Omega_t^c} |\phi(\xi)| \mathrm{d}S_t \mathrm{d}t$$

$$\geq (s-1) \int_\varepsilon^\delta t^{-s} \int_{\partial\Omega} |\phi(y)|[1 - (n-1)t\mathcal{H}(y) + c_3 t^2] \mathrm{d}S \mathrm{d}t$$

$$= M\varepsilon^{1-s} - (s-1)(n-1)M_\mathcal{H} \int_\varepsilon^\delta t^{1-s} \mathrm{d}t + c_4 \int_\varepsilon^\delta t^{2-s} \mathrm{d}t, \quad (4.25)$$

for some $c_3, c_4 \in \mathbb{R}$, and similarly

$$\int_\Omega \frac{|u_\varepsilon|}{d^{s-\beta}} \mathrm{d}x \geq M \int_\varepsilon^\delta t^{\beta-s} \mathrm{d}t - (n-1)M_\mathcal{H} \int_\varepsilon^\delta t^{1+\beta-s} \mathrm{d}t + c_5 \int_\varepsilon^\delta t^{2+\beta-s} \mathrm{d}t, \qquad (4.26)$$

for some $c_5 \in \mathbb{R}$. Thus inserting (4.22), (4.24), (4.25) into (4.21), and by (4.26) for $\beta = 1$, we get

$$Q_\beta[u_\varepsilon] \leq \frac{(n-1)M_\mathcal{H}[(s-1)\int_\varepsilon^\delta t^{1-s}\mathrm{d}t - \varepsilon^{2-s}] + K\int_\varepsilon^\delta t^{1-s}\mathrm{d}t + c_6 \int_\varepsilon^\delta t^{2-s}\mathrm{d}t}{M\int_\varepsilon^\delta t^{\beta-s}\mathrm{d}t - (n-1)M_\mathcal{H}\int_\varepsilon^\delta t^{1+\beta-s}\mathrm{d}t + c_5\int_\varepsilon^\delta t^{2+\beta-s}\mathrm{d}t}, \qquad (4.27)$$

for some $c_6 \in \mathbb{R}$. If $s = 2$ then

$$Q_1[u_\varepsilon] \leq \frac{((n-1)M_\mathcal{H} + K)\log(\delta/\varepsilon) + O_\varepsilon(1)}{M\log(\delta/\varepsilon) + O_\varepsilon(1)},$$

while if $s > 2$ then

$$Q_1[u_\varepsilon] \leq \frac{\frac{1}{s-2}((n-1)M_\mathcal{H} + K)\varepsilon^{2-s} + c_7 \int_\varepsilon^\delta t^{2-s}\mathrm{d}t}{\frac{1}{s-2}M\varepsilon^{2-s} - (n-1)M_\mathcal{H}\int_\varepsilon^\delta t^{2-s}\mathrm{d}t + c_8 \int_\varepsilon^\delta t^{3-s}\mathrm{d}t},$$

for some $c_7, c_8 \in \mathbb{R}$. In any case, letting $\varepsilon \downarrow 0$ we deduce $\mathcal{B}_1(\Omega) \leq \frac{(n-1)M_\mathcal{H}+K}{M}$. ∎

An immediate consequence is

**Corollary 4.23** (**Upper estimate**). *Let $\Omega$ be a bounded domain with boundary of class $C^2$. If $s \geq 2$ then*

$$\mathcal{B}_1(\Omega) \leq \frac{n-1}{|\partial\Omega|} \int_{\partial\Omega} \mathcal{H}(y)\mathrm{d}S$$

*where $\mathcal{H}(y)$ is the mean curvature at the point $y \in \partial\Omega$.*

**Proof.** Since $\Omega$ is bounded we can chose $\varphi \equiv 1$ in the above Theorem. ∎

The proof of Theorem I follows from Proposition 4.16, Theorem 4.17 and Corollary 4.23.

**Example 4.24** (**Ball**). Let $B_R$ be a ball of radius $R$. By Theorem 4.17 we have $\mathcal{B}_1(B_R) \geq \frac{n-1}{R}$ and by Corollary 4.23, $\mathcal{B}_1(B_R) \leq \frac{n-1}{R}$. We conclude that if $s \geq 2$ then $\mathcal{B}_1(B_R) = \frac{n-1}{R}$. See §4.3.

**Example 4.25** (**Infinite strip: proof of the optimality in Theorem IV**). Let $S_R = \{x = (x', x_n) : x' \in \mathbb{R}^{n-1}, 0 < x_n < 2R\}$. If $s \geq 2$, then combining Theorem 4.17 and Theorem 4.22 we can prove that $\mathcal{B}_1(S_R) = 0$. In fact we have $\mathcal{B}_\beta(S_R) = 0$ for any $1 < \beta \leq s - 1$ and in particular we will prove that if $\gamma = 1$, there is not positive constant $C$ such that (4.13) holds for $\gamma = 1$. To see this pick any $\phi \equiv \phi(x') \in C_c^1(\mathbb{R}^{n-1})$ such that $\mathrm{sprt}\{\phi\} \subset B_1 \subset \mathbb{R}^{n-1}$, where $B_1$ is the open ball in $\mathbb{R}^{n-1}$ with radius 1 centered at $0'$. Let $\eta > 0$ and set $\phi_\eta \equiv \phi_\eta(x') := \phi(\eta x')$. Note that $\mathrm{sprt}\{\phi_\eta\} \subset B_{1/\eta}$. Let also $0 < \varepsilon < \delta$ for some fixed $\delta \leq R$ (so that $d(x) = x_n$). The quotient in correspondence with (4.13) is

$$Q_\gamma[u] = \frac{\int_{S_R} \frac{|\nabla u|}{d^{s-1}}\mathrm{d}x - (s-1)\int_{S_R} \frac{|u|}{d^s}\mathrm{d}x}{\int_{S_R} \frac{|u|}{d}X^\gamma(\frac{d}{R})\mathrm{d}x} \tag{4.28}$$

As in the proof of Theorem 4.22, we test (4.28) with $u_{\varepsilon,\eta}(x) := \chi_{(\varepsilon,\delta)}(x_n)\phi_\eta(x')$ to arrive at

$$Q_\gamma[u_{\varepsilon,\eta}] = \frac{K_\eta \int_\varepsilon^\delta x_n^{1-s}\mathrm{d}x_n + 2M_\eta \delta^{1-s}}{M_\eta \int_\varepsilon^\delta x_n^{-1}X^\gamma(x_n/R)\mathrm{d}x_n},$$

where we have set $M_\eta := \int_{B_{1/\eta}} |\phi_\eta(x')|\mathrm{d}x'$ and $K_\eta := \int_{B_{1/\eta}} |\nabla_{x'}\phi_\eta(x')|\mathrm{d}x'$. Changing variables by $y' = \delta x'$ we obtain

$$\frac{K_\eta}{M_\eta} = \frac{K_1 \eta^{-(n-2)}}{M_1 \eta^{-(n-1)}} = \frac{K_1}{M_1}\eta,$$

where $M_1 = \int_{B_1} |\phi(y')| \mathrm{d}y'$ and $K_1 = \int_{B_1} |\nabla_{y'} \phi(y')| \mathrm{d}y'$. Thus

$$Q_\gamma[u_{\varepsilon,\eta}] = \frac{\frac{K_1}{M_1}\eta \int_\varepsilon^\delta x_n^{1-s} \mathrm{d}x_n + 2\delta^{1-s}}{\int_\varepsilon^\delta x_n^{-1} X^\gamma(x_n/R) \mathrm{d}x_n}.$$

Now we select $\eta = \varepsilon^{s-2+\epsilon}$ for some fixed $\epsilon > 0$. We deduce

$$Q_1[u_{\varepsilon,\eta}] = \frac{\frac{K_1}{M_1}\varepsilon^{s-2+\epsilon} \int_\varepsilon^\delta x_n^{1-s} \mathrm{d}x_n + 2\delta^{1-s}}{\log(\frac{X(\delta/R)}{X(\varepsilon/R)})}.$$

It follows that $Q_1[u_{\varepsilon,\eta}] \to 0$ as $\varepsilon \downarrow 0$. Thus for $\Omega = S_R$ inequality (4.13) does not hold when $\gamma = 1$ and the exponent 1 on the distance function in the remainder term in (4.13) cannot be increased.

## 4.3 The case of a ball

In this section we assume $\Omega$ is a ball of radius $R$. Without loss of generality we assume it is centered at the origin and denote it by $B_R$. The distance function to the boundary is then $d(x) = R - r$ where $r := |x|$. Moreover,

$$-\Delta d(x) = \frac{n-1}{R - d(x)}, \quad x \in B_R \setminus \{0\}. \tag{4.29}$$

This section is devoted to the proof of the following fact (Theorem III of the introduction)

**Theorem 4.26.** (1) *For all* $u \in C_c^\infty(B_R)$, $s \geq 2$ *and* $\gamma > 1$, *there holds*

$$\int_{B_R} \frac{|\nabla u|}{d^{s-1}} \mathrm{d}x \geq (s-1) \int_{B_R} \frac{|u|}{d^s} \mathrm{d}x + \sum_{k=1}^{[s]-1} \frac{n-1}{R^k} \int_{B_R} \frac{|u|}{d^{s-k}} \mathrm{d}x + \frac{C}{R^{s-1}} \int_{B_R} \frac{|u|}{d} X^\gamma\left(\frac{d}{R}\right) \mathrm{d}x, \quad (4.30)$$

*where* $C \geq \gamma - 1$. *The exponents* $s - k$; $k = 1, 2, ..., [s] - 1$, *on the distance function as well as the constants* $(n-1)/R^k$; $k = 1, 2, ..., [s] - 1$, *in the summation terms are optimal. If* $\gamma = 1$ *the above inequality fails in the sense of* (4.33).

(2) *For all* $u \in C_c^\infty(B_R)$, $1 \leq s < 2$ *and* $\gamma > 1$, *there holds*

$$\int_{B_R} \frac{|\nabla u|}{d^{s-1}} \mathrm{d}x \geq (s-1) \int_{B_R} \frac{|u|}{d^s} \mathrm{d}x + \frac{C}{R^{s-1}} \int_{B_R} \frac{|u|}{d} X^\gamma\left(\frac{d}{R}\right) \mathrm{d}x, \tag{4.31}$$

*where* $C \geq \gamma - 1$. *If* $\gamma = 1$ *the above inequality fails in the sense of* (4.33).

**Remark 4.27.** The optimality of the exponents and the constants stated in the above Theorem is meant in the following sense: for any $s \geq 1$ set

$$I_0[u] := \int_{B_R} \frac{|\nabla u|}{d^{s-1}} \mathrm{d}x - (s-1) \int_{B_R} \frac{|u|}{d^s} \mathrm{d}x,$$

and also for any $s \geq 2$ set

$$I_m[u] := I_0[u] - \sum_{k=1}^{m} \frac{n-1}{R^k} \int_{B_R} \frac{|u|}{d^{s-k}} \mathrm{d}x, \quad m = 1, ..., [s]-1.$$

Then for any $s \geq 2$

$$\inf_{u \in C_c^\infty(B_R) \setminus \{0\}} \frac{I_m[u]}{\int_{B_R} \frac{|u|}{d^\beta} \mathrm{d}x} = \begin{cases} (n-1)/R^{m+1}, & \text{if } \beta = s - m - 1 \\ 0, & \text{if } \beta > s - m - 1, \end{cases} \quad (4.32)$$

for all $m \in \{0, ..., [s]-2\}$. Further, for any $s \geq 1$

$$\inf_{u \in C_c^\infty(B_R) \setminus \{0\}} \frac{I_{[s]-1}[u]}{\int_{B_R} \frac{|u|}{d} X(d/R) \mathrm{d}x} = 0. \quad (4.33)$$

**Proof.** Inequality (4.31) is evident by Theorem 4.11. Let $s \geq 2$ and $\gamma > 1$. Since inequality (4.30) is scale invariant it suffices to prove it for $R = 1$. Testing (4.1) with

$$T(x) = -(d(x))^{1-s}[1 - (d(x))^{s-1} X^{\gamma-1}(d(x))] \nabla d(x); \quad x \in B_1 \setminus \{0\}.$$

we arrive at

$$\int_{B_1} \mathrm{div}(T)|u| \mathrm{d}x = (s-1) \int_{B_1} \frac{|u|}{d^s} \mathrm{d}x + \int_{B_1} \frac{|u|}{d^{s-1}} (1 - d^{s-1} X^{\gamma-1}(d))(-\Delta d) \mathrm{d}x$$
$$+ (\gamma - 1) \int_{B_1} \frac{|u|}{d} X^\gamma(d) \mathrm{d}x.$$

Thus, using (4.29) for $R = 1$ we obtain

$$\int_{B_1} \mathrm{div}(T)|u| \mathrm{d}x = (s-1) \int_{B_1} \frac{|u|}{d^s} \mathrm{d}x + (n-1) \int_{B_1} \frac{|u|}{d^{s-1}} \frac{1 - d^{s-1} X^{\gamma-1}(d)}{1-d} \mathrm{d}x$$
$$+ (\gamma - 1) \int_{B_1} \frac{|u|}{d} X^\gamma(d) \mathrm{d}x. \quad (4.34)$$

Since $s \geq 2$ we take into account in (4.34) the fact that

$$\frac{1 - d^{s-1} X^{\gamma-1}(d)}{1-d} \geq \frac{1 - d^{s-1}}{1-d} \geq \frac{1 - d^{[s]-1}}{1-d} = \sum_{k=1}^{[s]-1} d^{k-1}, \quad x \in B_1 \setminus \{0\},$$

and finally arrive at

$$I_0[u] \geq (s-1)\int_{B_1} \frac{|u|}{d^s}\mathrm{d}x + (n-1)\sum_{k=1}^{[s]-1}\int_{B_1} \frac{|u|}{d^{s-k}}\mathrm{d}x + (\gamma-1)\int_{B_1} \frac{|u|}{d}X^\gamma(d)\mathrm{d}x,$$

which is (4.30) for $R=1$.

We next prove (4.32). Suppose first that $2 \leq s < 3$. In this case all we have to prove is that

$$\inf_{u \in C_0^1(B_1)\setminus\{0\}} \frac{I_0[u]}{\int_{B_1}\frac{|u|}{d^\beta}\mathrm{d}x} = \begin{cases} n-1, & \text{if } \beta = s-1 \\ 0, & \text{if } \beta > s-1. \end{cases} \qquad (4.35)$$

To this end we pick $u_\delta(x) = \chi_{B_{1-\delta}}(x)$ where $x \in B_1$ and $0 < \delta < 1$. This function is in $BV(B_1)$ and we can take a $C_c^\infty$ approximation of it, so that the calculations bellow to hold in the limit. The distributional gradient of $u_\delta$ is $\nabla u_\delta = -\vec{\nu}_{\partial B_{1-\delta}}\delta_{\partial B_{1-\delta}}$ and the total variation of $\nabla u_\delta$ is $|\nabla u_\delta| = \delta_{\partial B_{1-\delta}}$. Using co-area formula we get

$$\frac{I_0[u_\delta]}{\int_{B_1}\frac{|u_\delta|}{d^\beta}\mathrm{d}x} = \frac{\delta^{1-s}|\partial B_{1-\delta}| - (s-1)\int_0^{1-\delta}(1-r)^{-s}|\partial B_r|\mathrm{d}r}{\int_0^{1-\delta}(1-r)^{-\beta}|\partial B_r|\mathrm{d}r}$$

$$= \frac{\delta^{1-s}(1-\delta)^{n-1} - \int_0^{1-\delta}((1-r)^{1-s})'r^{n-1}\mathrm{d}r}{\int_0^{1-\delta}(1-r)^{-\beta}r^{n-1}\mathrm{d}r}$$

$$= (n-1)\frac{\int_0^{1-\delta}(1-r)^{1-s}r^{n-2}\mathrm{d}r}{\int_0^{1-\delta}(1-r)^{-\beta}r^{n-1}\mathrm{d}r}.$$

Thus

$$\frac{I_0[u_\delta]}{\int_{B_1}\frac{|u_\delta|}{d^\beta}\mathrm{d}x} \to \begin{cases} n-1, & \text{if } \beta = s-1 \\ 0, & \text{if } \beta > s-1 \end{cases}, \quad \text{as } \delta \downarrow 0.$$

Assume next that $3 \leq s < 4$. This time we have besides (4.35) to prove that

$$\inf_{u \in C_c^\infty(B_1)\setminus\{0\}} \frac{I_1[u]}{\int_{B_1}\frac{|u|}{d^\beta}\mathrm{d}x} = \begin{cases} n-1, & \text{if } \beta = s-2 \\ 0, & \text{if } \beta > s-2. \end{cases}$$

Picking the same $u_\delta$ as before and performing the same integration by parts in the second term of the numerator, we conclude

$$\frac{I_1[u_\delta]}{\int_{B_1}\frac{|u_\delta|}{d^\beta}\mathrm{d}x} = \frac{(n-1)\int_0^{1-\delta}(1-r)^{1-s}r^{n-2}\mathrm{d}r - (n-1)\int_0^{1-\delta}(1-r)^{1-s}r^{n-1}\mathrm{d}r}{\int_0^{1-\delta}(1-r)^{-\beta}r^{n-1}\mathrm{d}r}$$

$$= (n-1)\frac{\int_0^{1-\delta}(1-r)^{2-s}r^{n-2}\mathrm{d}r}{\int_0^{1-\delta}(1-r)^{-\beta}r^{n-1}\mathrm{d}r}.$$

Thus

$$\frac{I_1[u_\delta]}{\int_{B_1} \frac{|u_\delta|}{d^\beta} dx} \to \begin{cases} n-1, & \text{if } \beta = s-2 \\ 0, & \text{if } \beta > s-2, \end{cases} \quad \text{as } \delta \downarrow 0.$$

We continue in the same fashion for $4 \leq s < 5$, then $5 \leq s < 6$ and so on.

Next we prove (4.33). We pick $u_\delta$ as before and perform the same integration by parts, to get

$$\begin{aligned}
\frac{I_{[s]-1}[u_\delta]}{\int_{B_1} \frac{|u_\delta|}{d} X(d) dx} &= \frac{(n-1)\int_0^{1-\delta}(1-r)^{1-s}r^{n-2}dr - (n-1)\sum_{k=1}^{[s]-1}\int_0^{1-\delta}(1-r)^{k-s}r^{n-1}dr}{\int_0^{1-\delta}(1-r)^{-1}r^{n-1}X(1-r)dr} \\
&= (n-1)\frac{\int_0^{1-\delta}(1-r)^{[s]-s}r^{n-2}dr}{\int_1^{1-\log\delta} t^{-1}(1-e^{1-t})^{n-1}dt} \\
&=: (n-1)\frac{N_\delta}{D_\delta}.
\end{aligned}$$

Since $[s] - s > -1$ we have $N_\delta = O_\delta(1)$ as $\delta \downarrow 0$. Also, $D_\delta \geq \int_1^{1-\log\delta} t^{-1} dt + O_\delta(1) \to \infty$ as $\delta \downarrow 0$. ∎

# Chapter 5

# Hardy-Sobolev type inequalities for $p > n$ - Distance to the origin

In this chapter we will prove Theorems V, VI and VII stated in the introduction. These theorems are the content of the work [Ps2]. In this chapter we set

$$I[u] := \int_\Omega |\nabla u|^p \mathrm{d}x - \Big(\frac{p-n}{p}\Big)^p \int_\Omega \frac{|u|^p}{|x|^p} \mathrm{d}x; \quad p > n,$$

whenever $u \in C_c^\infty(\Omega \setminus \{0\})$.

## 5.1 The Hardy-Sobolev inequality for $p > n \geq 1$

In this section we prove Theorem V. For $n = 1$, by (2.26) we have

$$\begin{aligned}
u(x) &\leq \frac{1}{2}\int_0^R |u'(t)|\mathrm{d}t \\
(\text{setting } u(t) = t^{1-1/p}v(t)) &\leq \frac{1}{2}\int_0^R t^{1-1/p}|v'(t)|\mathrm{d}t + \frac{p-1}{2p}\int_0^R t^{-1/p}|v(t)|\mathrm{d}t \\
&\leq \int_0^R t^{1-1/p}|v'(t)|\mathrm{d}t,
\end{aligned}$$

by (2.1) for $q = 1$ and $s = 1/p$. The proof follows applying Hölder's inequality and using (2.10) for $p \geq 2$ or (2.12) if $1 < p < 2$. We omit further details.

Assume next that $n \geq 2$. From Lemma 2.16 we get

$$u(x) \leq \frac{1}{n\omega_n}\int_\Omega \frac{|\nabla u(z)|}{|x-z|^{n-1}} \mathrm{d}z,$$

for $x \in \Omega$. Setting $u(z) = |z|^{1-n/p}v(z)$, we arrive at

$$n\omega_n|u(x)| \leq \underbrace{\int_\Omega \frac{|z|^{1-n/p}|\nabla v(z)|}{|x-z|^{n-1}}\mathrm{d}z}_{=:K(x)} + \frac{p-n}{p}\underbrace{\int_\Omega \frac{|v(z)|}{|z|^{n/p}|x-z|^{n-1}}\mathrm{d}z}_{=:\Lambda(x)}. \tag{5.1}$$

We first estimate $K(x)$. Using Hölder's inequality we get

$$K(x) \leq \left(\underbrace{\int_\Omega \frac{1}{|x-z|^{(n-1)p/(p-1)}} \mathrm{d}z}_{=:M(x)}\right)^{1-1/p} \left(\int_\Omega |z|^{p-n} |\nabla v(z)|^p \mathrm{d}z\right)^{1/p}. \tag{5.2}$$

Note that $M(x)$ is finite since $(n-1)p/(p-1) < n$ if and only if $p > n$. To estimate $M(x)$ we set $R := [\mathcal{L}^n(\Omega)/\omega_n]^{1/n}$, so that the volume of a ball with radius $R$ to be equal to the volume of $\Omega$. Then $M(x)$ increases if we change the domain of integration from $\Omega$ to $B_R(x)$. Therefore

$$\begin{aligned} M(x) &\leq \int_{B_R(x)} \frac{1}{|x-z|^{(n-1)p/(p-1)}} \mathrm{d}z \\ &= n\omega_n \frac{p-1}{p-n} [\mathcal{L}^n(\Omega)/\omega_n]^{(p-n)/n(p-1)}, \end{aligned} \tag{5.3}$$

and using (2.10) in the second factor of (5.2), we get

$$K(x) \leq C_1(n,p) [\mathcal{L}^n(\Omega)]^{1/n - 1/p} (I[u])^{1/p}. \tag{5.4}$$

Next we bound $\Lambda(x)$. Using Hölder's inequality with conjugate exponents $p/(p-1-\varepsilon)$ and $p/(1+\varepsilon)$, where $0 < \varepsilon < (p-n)/n$ is fixed and depending only on $n, p$, we get

$$\Lambda(x) \leq \left(\underbrace{\int_\Omega \frac{1}{|x-z|^{(n-1)p/(p-1-\varepsilon)}} \mathrm{d}z}_{=:M(\varepsilon,x)}\right)^{1-(1+\varepsilon)/p} \left(\int_\Omega \frac{|v(z)|^{p/(1+\varepsilon)}}{|z|^{n/(1+\varepsilon)}} \mathrm{d}z\right)^{(1+\varepsilon)/p}. \tag{5.5}$$

Note that $M(\varepsilon, x)$ is finite since $(n-1)p/(p-1-\varepsilon) < n$ if and only if $\varepsilon < (p-n)/n$. As before, we set $R := [\mathcal{L}^n(\Omega)/\omega_n]^{1/n}$ and $M(\varepsilon, x)$ increases if we change the domain of integration from $\Omega$ to $B_R(x)$. Therefore,

$$\begin{aligned} M(\varepsilon,x) &\leq \int_{B_R(x)} \frac{1}{|x-z|^{(n-1)p/(p-1-\varepsilon)}} \mathrm{d}z \\ &= n\omega_n \frac{p-1-\varepsilon}{p-n-n\varepsilon} [\mathcal{L}^n(\Omega)/\omega_n]^{(p-n-n\varepsilon)/n(p-1-\varepsilon)}, \end{aligned} \tag{5.6}$$

and using inequality (2.5) of Remark 2.4 with $s = n/(1+\varepsilon)$ and $q = p/(1+\varepsilon)$ in the second factor of (5.5), we obtain

$$\Lambda(x) \leq C_2(n,p) [\mathcal{L}^n(\Omega)]^{1/n - 1/p - \varepsilon/p} \left(\int_\Omega |z|^{(p-n)/(1+\varepsilon)} |\nabla v(z)|^{p/(1+\varepsilon)} \mathrm{d}z\right)^{(1+\varepsilon)/p}.$$

Using once more Hölder's inequality with conjugate exponents $1 + 1/\varepsilon$ and $1 + \varepsilon$, we get

$$\begin{aligned}
\Lambda(x) &\leq C_2(n,p)[\mathcal{L}^n(\Omega)]^{1/n-1/p-\varepsilon/p}\left[[\mathcal{L}^n(\Omega)]^{\varepsilon/(1+\varepsilon)}\left(\int_\Omega |z|^{p-n}|\nabla v(z)|^p \mathrm{d}z\right)^{1/(1+\varepsilon)}\right]^{(1+\varepsilon)/p} \\
&= C_2(n,p)[\mathcal{L}^n(\Omega)]^{1/n-1/p}\left(\int_\Omega |z|^{p-n}|\nabla v(z)|^p \mathrm{d}z\right)^{1/p}. \\
&\leq C_3(n,p)[\mathcal{L}^n(\Omega)]^{1/n-1/p}(I[u])^{1/p},
\end{aligned} \qquad (5.7)$$

by (2.10). The proof follows inserting (5.7) and (5.4) in (5.1). ∎

## 5.2 An optimal Hardy-Morrey inequality for $n = 1$

Theorem VI in the one-dimensional case has an easier proof. We present it in this separate section. Firstly note that it suffices to restrict ourselves in the case $\Omega = (0,1)$ We start with the following lemma

**Lemma 5.1.** *Let $q > 1$, $\beta > 1 - q$. There exists a positive constant $c = c(q, \beta)$ such that, for all $v \in C_c^\infty(0,1)$ and any $D \geq 1$*

$$\sup_{x \in (0,1)} \left\{|v(x)|X^{(\beta+q-1)/q}(x/D)\right\} \leq c\left(\int_0^1 t^{q-1}|v'(t)|^q X^\beta(t/D)\mathrm{d}t\right)^{1/q}.$$

**Proof.** Letting $D \geq 1$, we have

$$\begin{aligned}
v(x) &= -\int_x^D v'(t)\mathrm{d}t \\
&\leq \left(\int_x^D t^{-1}X^{-\beta/(q-1)}(t/D)\mathrm{d}t\right)^{1-1/q}\left(\int_x^D t^{q-1}|v'(t)|^q X^\beta(t/D)\mathrm{d}t\right)^{1/q},
\end{aligned}$$

where in the last inequality we have used Hölder's inequality with conjugate exponents $q$ and $q/(q-1)$. Since $v \in C_c^\infty(0,1)$ the second integral is actually over $(x,1)$. In addition

$$[X^{-1-\beta/(q-1)}(t/D)]' = \left(-1 - \frac{\beta}{q-1}\right)t^{-1}X^{-\beta/(q-1)}(t/D),$$

so that the left integral can be computed. We find

$$|v(x)| \leq c\left(X^{-1-\beta/(q-1)}(x/D) - 1\right)^{1-1/q}\left(\int_0^1 t^{q-1}|v'(t)|^q X^\beta(t/D)\mathrm{d}t\right)^{1/q},$$

where $c = [(q-1)/(q-1+\beta)]^{1-1/q}$, or

$$|v(x)|X^{(\beta+q-1)/q}(x/D) \leq c\left(1-X^{(\beta+q-1)/(q-1)}(x/D)\right)^{1-1/q}\left(\int_0^1 t^{q-1}|v'(t)|^q X^\beta(t/D)dt\right)^{1/q}.$$

The result follows since $(1 - X^{(\beta+q-1)/(q-1)}(x/D))^{1-1/q} \leq 1$, for all $x \in (0,1]$. ∎

We proceed by proving Theorem VI with one point in the Hölder seminorm taken to be the origin. In particular, we have

**Proposition 5.2.** *Let $p > 1$. There exists positive constant $c_p$ depending only on $p$, such that for all $u \in C_c^\infty(0,1)$ and any $D \geq 1$*

$$\sup_{x \in (0,1)}\left\{\frac{|u(x)|}{x^{1-1/p}}X^{1/p}(x/D)\right\} \leq c_p(I[u])^{1/p}. \tag{5.8}$$

**Proof.** We set $v(x) = x^{-1+1/p}u(x)$. If $1 < p < 2$, by Lemma 5.1 for $q = p$ and $\beta = 2 - p$, we have

$$|v(x)|X^{1/p}(x/D) \leq c_p\left(\int_0^1 t^{p-1}|v'(t)|^p X^{2-p}(t/D)dt\right)^{1/p},$$

for any $D \geq 1$. The result follows by (2.12). If $p \geq 2$, by Lemma 5.1 for $q = 2$ and $\beta = 0$, we have

$$|w(x)|X^{1/2}(x/D) \leq c_p\left(\int_0^1 t|w'(t)|^2 dt\right)^{1/2},$$

for any $w \in C_c^\infty(0,1)$ and any $D \geq 1$. For $w(x) = |v(x)|^{p/2}$, we obtain

$$|v(x)|X^{1/p}(x/D) \leq c_p\left(\int_0^1 t|v(t)|^{p-2}|v'(t)|^2 dt\right)^{1/p}.$$

The result follows by (2.11). ∎

Now we use (5.8) to obtain its counterpart inequality with the exact Hölder seminorm.

**Proof of Theorem VI in case $n = 1$.** For $0 < y < x < 1$ we have

$$|u(x) - u(y)| = \left|\int_y^x u'(t)dt\right|$$

$$(\text{setting } u(t) = t^{1-1/p}v(t)) \leq \underbrace{\int_y^x t^{1-1/p}|v'(t)|dt}_{:=K(x,y)} + \frac{p-1}{p}\underbrace{\int_y^x \frac{|v(t)|}{t^{1/p}}dt}_{:=\Lambda(x,y)}. \tag{5.9}$$

To estimate $\mathrm{K}(x,y)$ we use Hölder's inequality to get

$$\mathrm{K}(x,y) \leq (x-y)^{1-1/p}\left(\int_y^x t^{p-1}|v'(t)|^p \mathrm{d}t\right)^{1/p}$$

$$(\text{by (2.10)}) \leq c_1(p)(x-y)^{1-1/p}(I[u])^{1/p}$$

$$\leq c_1(p)(x-y)^{1-1/p}X^{-1/p}((x-y)/D)(I[u])^{1/p}, \quad (5.10)$$

for any $D \geq 1$, since $0 \leq X(t) \leq 1$ for all $t \in (0,1]$. To estimate $\Lambda(x,y)$ we return to the original variable by $v(t) = u(t)/t^{1-1/p}$, thus

$$\Lambda(x,y) = \int_y^x \frac{|u(t)|}{t}\mathrm{d}t.$$

Inserting (5.8) in $\Lambda(x,y)$ we obtain

$$\Lambda(x,y) \leq c_2(p)(I[u])^{1/p}\int_y^x t^{-1/p}X^{-1/p}(t/D)\mathrm{d}t$$

$$\leq c_3(p)(x-y)^{1-1/p}X^{-1/p}((x-y)/D)(I[u])^{1/p}, \quad (5.11)$$

for any $D \geq e^\eta$, where $\eta$ depends only on $p$, due to Lemma A.1-(ii) for $\alpha = -1/p$ and $\beta = 1/p$. Coupling (5.10) and (5.11) with (5.9), we end up with

$$|u(x) - u(y)| \leq c_4(p)(x-y)^{1-1/p}X^{-1/p}((x-y)/D)(I[u])^{1/p},$$

for all $0 < y < x < 1$ and any $D \geq e^\eta$, which is the desired estimate with $B = e^\eta$. ∎

## 5.3 A nonhomogeneous local remainder

Here we prove Theorem VII in case $0 \in \overline{B}_r$. To emphasize that under this assumption we will prove (1.29) for general $p, q > 1$ satisfying $1 \leq q < p$, $p \neq n$, we present it as a separate theorem. The proof in case $0 \notin \overline{B}_r$ and $p > n$ is given after the proof of Proposition 5.4.

**Theorem 5.3.** *Suppose $\Omega$ is a bounded domain in $\mathbb{R}^n$; $n \geq 2$, containing the origin and let $p \neq n$ and $1 \leq q < p$. There exist constants $\Theta = \Theta(n,p,q) \geq 0$ and $C = C(n,p,q) > 0$ such that for all $u \in C_c^\infty(\Omega \setminus \{0\})$, any open ball $B_r$ with $r \in (0, \mathrm{diam}(\Omega))$ and $0 \in \overline{B}_r$, and any $D \geq e^\Theta \mathrm{diam}(\Omega)$*

$$r^{n/p}X^{1/p}(r/D)\left(\frac{1}{|B_r|}\int_{B_r}\frac{|u|^q}{|x|^q}\mathrm{d}x\right)^{1/q} \leq C\left(\int_\Omega |\nabla u|^p\mathrm{d}x - \left|\frac{p-n}{p}\right|^p\int_\Omega \frac{|u|^p}{|x|^p}\mathrm{d}x\right)^{1/p}, \quad (5.12)$$

*where $X(t) = (1 - \log t)^{-1}$; $t \in (0, 1]$.*

**Proof.** Let $\Omega$ be a bounded domain in $\mathbb{R}^n$ containing the origin and let $1 \leq q < p$, $p \neq n$. Suppose $u \in C_c^\infty(\Omega \setminus \{0\})$ and let also $B_r$ be any ball containing zero in its closure with $r \in (0, \text{diam}(\Omega))$. Setting $u(x) = |x|^{1-n/p} v(x)$ we get

$$\int_{B_r} \frac{|u|^q}{|x|^q} dx = \int_{B_r} \frac{|v|^q}{|x|^{nq/p}} dx$$
$$\leq \left[\frac{pq}{n(p-q)}\right]^q \underbrace{\int_{B_r} |x|^{q(p-n)/p} |\nabla v|^q dx}_{=:N_r}$$
$$+ \frac{pq}{n(p-q)} \underbrace{\int_{\partial B_r} \frac{|v|^q}{|x|^{nq/p}} x \cdot \vec{\nu} dS_x}_{=:P_r}, \quad (5.13)$$

where we have used Lemma 2.3 for $V = B_r$ and $s = nq/p$. We use Hölder's inequality with conjugate exponents $p/(p-q)$ and $p/q$, to get

$$N_r \leq (\omega_n r^n)^{(p-q)/p} \left(\int_{B_r} |x|^{p-n} |\nabla v|^p dx\right)^{q/p}$$
$$\text{(by (2.10))} \leq C_1(n,p,q) r^{n(p-q)/p} (I[u])^{q/p}$$
$$\leq C_1(n,p,q) r^{n(p-q)/p} X^{-q/p}(r/D) (I[u])^{q/p}, \quad (5.14)$$

for any $D \geq \text{diam}(\Omega)$ since $0 \leq X(t) \leq 1$ for all $t \in (0,1]$. For $P_r$ we write first

$$P_r = \int_{\partial B_r} \left\{X^{-q/p}(|x|/D)\right\} \left\{\frac{|v|^q}{|x|^{nq/p}} X^{q/p}(|x|/D)\right\} x \cdot \vec{\nu} dS_x$$
$$\leq \underbrace{\left(\int_{\partial B_r} X^{-q/(p-q)}(|x|/D) x \cdot \vec{\nu} dS_x\right)^{(p-q)/p}}_{=:S_r} \underbrace{\left(\int_{\partial B_r} \frac{|v|^p}{|x|^n} X(|x|/D) x \cdot \vec{\nu} dS_x\right)^{q/p}}_{=:T_r} (5.15)$$

where we have used once more Hölder's inequality with exponents $p/(p-q)$ and $p/q$. Note that $\nu$ is the outward pointing unit normal vector field along $\partial B_r$ and that since $0 \in \overline{B}_r$ we have $x \cdot \nu \geq 0$ for all $x \in \partial B_r$. By the divergence theorem we have

$$S_r = \int_{B_r} \text{div}[X^{-q/(p-q)}(|x|/D) \, x] dx$$
$$= n \int_{B_r} X^{-q/(p-q)}(|x|/D) dx - \frac{q}{p-q} \int_{B_r} X^{1-q/(p-q)}(|x|/D) dx$$
$$\leq n \int_{B_r(0)} X^{-q/(p-q)}(|x|/D) dx,$$

since the integral increases if we change the domain of integration from $B_r$ to $B_r(0)$. Thus

$$\begin{aligned} S_r &\leq n^2 \omega_n \int_0^r t^{n-1} X^{-q/(p-q)}(t/D) \mathrm{d}t \\ &\leq C_2(n) r^n X^{-q/(p-q)}(r/D), \end{aligned} \quad (5.16)$$

for any $D \geq e^\eta \operatorname{diam}(\Omega)$, where $\eta \geq 0$ depends only on $n, p, q$, due to Lemma A.1-(i) for $\alpha = n-1$ and $\beta = q/(p-q)$ in the Appendix. $T_r$ will be estimated after an integration by parts. More precisely

$$T_r = \int_{B_r} \operatorname{div}\left\{ \frac{X(|x|/D)}{|x|^n} x \right\} |v|^p \mathrm{d}x + \int_{B_r} \frac{X(|x|/D)}{|x|^n} x \cdot \nabla(|v|^p) \mathrm{d}x.$$

A simple calculation shows that $\operatorname{div}\{\frac{X(|x|/D)}{|x|^n} x\} = \frac{X^2(|x|/D)}{|x|^n}$ for any $x \in \Omega \setminus \{0\}$. In the second integral we compute the gradient and note that $x \cdot \nabla |v(x)| \leq |x||\nabla v(x)|$ for a.e. $x \in \Omega$. Thus,

$$T_r \leq \int_\Omega \frac{|v|^p}{|x|^n} X^2(|x|/D) \mathrm{d}x + p \int_\Omega |x|^{1-n} |v|^{p-1} |\nabla v| X(|x|/D) \mathrm{d}x.$$

We rearrange the integrand in the second integral above as follows

$$\begin{aligned} T_r &\leq \int_\Omega \frac{|v|^p}{|x|^n} X^2(|x|/D) \mathrm{d}x + p \int_\Omega \left\{ |x|^{1-n/2} |v|^{p/2-1} |\nabla v| \right\} \left\{ \frac{|v|^{p/2}}{|x|^{n/2}} X(|x|/D) \right\} \mathrm{d}x \\ &\leq \int_\Omega \frac{|v|^p}{|x|^n} X^2(|x|/D) \mathrm{d}x + p \left( \int_\Omega |x|^{2-n} |v|^{p-2} |\nabla v|^2 \mathrm{d}x \right)^{1/2} \left( \int_\Omega \frac{|v|^p}{|x|^n} X^2(|x|/D) \mathrm{d}x \right)^{1/2}, \end{aligned}$$

by the Cauchy-Schwarz inequality. According to Theorem 2.11-(i), there exist constants $\theta \geq 0$ and $b > 0$, both depending only on $n, p$, such that for any $D \geq e^\theta \sup_{x \in \Omega} |x|$, the first term and the second radicand on the right hand side (when returned to the original function by $v(x) = |x|^{n/p-1} u(x)$) are bounded above by $bI[u]$. Due to (2.11), the first radicand is also bounded above by $C(n, p)I[u]$. It follows that

$$T_r \leq C_3(n, p) I[u], \quad (5.17)$$

for any $D \geq e^\theta \sup_{x \in \Omega} |x|$. Setting $\Theta = \max\{\theta, \eta\}$ and noting that $0 \in \Omega$ implies $\sup_{x \in \Omega} |x| \leq \operatorname{diam}(\Omega)$, we insert (5.16) and (5.17) into estimate (5.15) to obtain

$$P_r \leq C_4(n, p, q) r^{n(p-q)/p} X^{-q/p}(r/D) (I[u])^{q/p},$$

for any $D \geq e^\Theta \operatorname{diam}(\Omega)$. The last inequality together with (5.14), when applied to estimate (5.13), give

$$\int_{B_r} \frac{|u|^q}{|x|^q} \mathrm{d}x \leq C_5(n, p, q) r^{n(p-q)/p} X^{-q/p}(r/D) (I[u])^{q/p}.$$

for any $D \geq e^\Theta \operatorname{diam}(\Omega)$. Rearranging, raising in $1/q$ and taking the supremum over all $B_r$ containing zero with $r \in (0, \operatorname{diam}(\Omega))$, the result follows. ∎

## 5.4 An optimal Hardy-Morrey inequality for $n \geq 2$

In this section we prove Theorem VI when $n \geq 2$. We first obtain (1.28) with one point in the Hölder seminorm taken to be the origin. More precisely, we prove

**Proposition 5.4.** *Let $p > n \geq 2$ and suppose $\Omega$ is a bounded domain in $\mathbb{R}^n$ containing the origin. There exist constants $\Theta \geq 0$ and $C > 0$ both depending only on $n, p$, such that for all $u \in C_c^\infty(\Omega \setminus \{0\})$ and any $D \geq e^\Theta \operatorname{diam}(\Omega)$*

$$\sup_{x \in \Omega} \left\{ \frac{|u(x)|}{|x|^{1-n/p}} X^{1/p}(|x|/D) \right\} \leq C (I[u])^{1/p}. \tag{5.18}$$

**Proof.** Let $B_r$ be a ball containing zero with $r \in (0, \operatorname{diam}(\Omega))$ and set $u_{B_r} = |B_r|^{-1} \int_{B_r} u \, dz$. Letting $x \in B_r$, the local version of the representation formula (Lemma 2.19) asserts

$$|u(x) - u_{B_r}| \leq \frac{2^n}{n \omega_n} \int_{B_r} \frac{|\nabla u(z)|}{|x-z|^{n-1}} dz.$$

Setting $u(z) = |z|^{1-n/p} v(z)$, we arrive at

$$\frac{n \omega_n}{2^n} |u(x) - u_{B_r}| \leq \underbrace{\int_{B_r} \frac{|z|^{1-n/p} |\nabla v(z)|}{|x-z|^{n-1}} dz}_{=: K_r(x)} + \frac{p-n}{p} \underbrace{\int_{B_r} \frac{|v(z)|}{|z|^{n/p} |x-z|^{n-1}} dz}_{=: \Lambda_r(x)}. \tag{5.19}$$

We will derive suitable bounds for $K_r(x), \Lambda_r(x)$. For $K_r(x)$ we use Hölder's inequality

$$K_r(x) \leq \left( \int_{B_r} \frac{1}{|x-z|^{(n-1)p/(p-1)}} dz \right)^{1-1/p} \left( \int_{B_r} |z|^{p-n} |\nabla v|^p dz \right)^{1/p}.$$

Both integrals increase if we integrate over $B_r(x)$ and $\Omega$ respectively. Hence

$$K_r(x) \leq \left( \int_{B_r(x)} \frac{1}{|x-z|^{(n-1)p/(p-1)}} dz \right)^{1-1/p} \left( \int_\Omega |z|^{p-n} |\nabla v|^p dz \right)^{1/p}.$$

Computing the first integral and using (2.10) for the second, we arrive at

$$\begin{aligned} K_r(x) &\leq C_1(n,p) r^{1-n/p} (I[u])^{1/p} \\ &\leq C_1(n,p) r^{1-n/p} X^{-1/p}(r/D) (I[u])^{1/p}, \end{aligned} \tag{5.20}$$

for any $D \geq \operatorname{diam}(\Omega)$, where the last inequality follows since $0 < X(t) \leq 1$ for all $t \in (0, 1]$. Next we bound $\Lambda_r(x)$. Using Hölder's inequality with conjugate exponents $p/(p-1-\varepsilon)$ and $p/(1+\varepsilon)$, where $0 < \varepsilon < (p-n)/n$ is fixed but depending only on $n, p$, we obtain

$$\Lambda_r(x) \leq \underbrace{\left(\int_{B_r} \frac{1}{|x-z|^{(n-1)p/(p-1-\varepsilon)}} dz\right)^{1-(1+\varepsilon)/p}}_{=:M_r(x)} \left(\int_{B_r} \frac{|v|^{p/(1+\varepsilon)}}{|z|^{n/(1+\varepsilon)}} dz\right)^{(1+\varepsilon)/p}. \quad (5.21)$$

By (5.6) and returning to the original function on the second integral on the right of (5.21), we obtain

$$\Lambda_r(x) \leq C_2(n,p) r^{1-n/p-n\varepsilon/p} \left(\int_{B_r} \frac{|u|^{p/(1+\varepsilon)}}{|z|^{p/(1+\varepsilon)}} dz\right)^{(1+\varepsilon)/p}.$$

Using Theorem VII with $q = p/(1+\varepsilon)$,

$$\Lambda_r(x) \leq C_3(n,p) r^{1-n/p-n\varepsilon/p} \left[r^{n\varepsilon/(1+\varepsilon)} X^{-1/(1+\varepsilon)}(r/D)(I[u])^{1/(1+\varepsilon)}\right]^{(1+\varepsilon)/p}$$
$$= C_3(n,p) r^{1-n/p} X^{-1/p}(r/D)(I[u])^{1/p}, \quad (5.22)$$

for any $D \geq e^\Theta \operatorname{diam}(\Omega)$, where $\Theta$ depends only on $n, p, \varepsilon$ (and thus only on $n, p$).

Applying estimates (5.22) and (5.20) to estimate (5.19), we finally conclude

$$|u(x) - u_{B_r}| \leq C_4(n,p) r^{1-n/p} X^{-1/p}(r/D)(I[u])^{1/p}, \quad (5.23)$$

for all $x \in B_r$ and any $D \geq e^\Theta \operatorname{diam}(\Omega)$. Since $0 \in B_r$, it follows from (6.21), that

$$|u_{B_r}| \leq C_4(n,p) r^{1-n/p} X^{-1/p}(r/D)(I[u])^{1/p}.$$

Hence

$$|u(x)| \leq |u(x) - u_{B_r}| + |u_{B_r}|$$
$$\leq 2C_4(n,p) r^{1-n/p} X^{-1/p}(r/D)(I[u])^{1/p},$$

for all $x \in B_r$ and any $D \geq e^\Theta \operatorname{diam}(\Omega)$. Now if $x \in \Omega$, we consider a ball $B_r$ of radius $r = |x|$, containing $x$. Then the previous inequality yields

$$|u(x)| \leq C(n,p) |x|^{1-n/p} X^{-1/p}(|x|/D)(I[u])^{1/p},$$

for any $D \geq e^\Theta \operatorname{diam}(\Omega)$. Rearranging and taking the supremum over all $x \in \Omega$, the result follows. ∎

**Completion of proof of Theorem VII.** Let $r \in (0, \operatorname{diam}(\Omega))$ and $p > n$. Using (5.18) we obtain

$$\int_{B_r} \frac{|u|^q}{|x|^q} dx \leq C^q (I[u])^{q/p} \int_{B_r} \frac{1}{|x|^{nq/p} X^{q/p}(|x|/D)} dx$$
$$\leq C^q (I[u])^{q/p} \int_{B_r(0)} \frac{1}{|x|^{nq/p} X^{q/p}(|x|/D)} dx$$
$$= C^q n \omega_n (I[u])^{q/p} \int_0^r t^{n-1-nq/p} X^{-q/p}(t/D) dt$$
$$\leq C(n,q,p) r^{n-nq/p} X^{-q/p}(r/D)(I[u])^{q/p},$$

for any $D \geq \max\{e^\eta, e^\Theta\} \operatorname{diam}(\Omega)$, where $\eta = \eta(n, p, q)$, by Lemma A.1-(i) (with $\alpha = n - 1 - nq/p$ and $\beta = q/p$). Rearranging, raising in $1/q$ and taking the supremum over all $B_r$ with $r \in (0, \operatorname{diam}(\Omega))$, we arrive at (1.29) without having assumed $0 \in \overline{B}_r$. ∎

Now we utilize (5.18) in order to obtain its counterpart inequality with the exact Hölder seminorm, i.e. inequality (1.28).

**Proof of Theorem VI in case $n \geq 2$.** Letting $x, y \in \Omega$ with $x \neq y$, we consider a ball $B_r$ of radius $r := |x - y|$ containing $x, y$. Note that $r \in (0, \operatorname{diam}(\Omega))$. We have

$$
\begin{aligned}
|u(x) - u(y)| &\leq |u(x) - u_{B_r}| + |u(y) - u_{B_r}| \\
&\leq \frac{2^n}{n\omega_n}\Big\{ \underbrace{\int_{B_r} \frac{|\nabla u(z)|}{|x-z|^{n-1}} \mathrm{d}z}_{=:J(x)} + \underbrace{\int_{B_r} \frac{|\nabla u(z)|}{|y-z|^{n-1}} \mathrm{d}z}_{=:J(y)} \Big\},
\end{aligned}
\qquad (5.24)
$$

where we have used Lemma 2.19 twice. We will bound $J(x)$ independently on $x$ so that the same estimate holds also for $J(y)$. We start with the substitution $u(z) = |z|^{1-n/p} v(z)$, to get

$$
J(x) \leq \underbrace{\int_{B_r} \frac{|z|^{1-n/p}|\nabla v(z)|}{|x-z|^{n-1}} \mathrm{d}z}_{=:\mathrm{K}_r(x)} + \frac{p-n}{p} \underbrace{\int_{B_r} \frac{|v(z)|}{|z|^{n/p}|x-z|^{n-1}} \mathrm{d}z}_{=:\Lambda_r(x)}.
\qquad (5.25)
$$

We estimate $\mathrm{K}_r(x)$ in the same manner as in (5.20). So

$$
\mathrm{K}_r(x) \leq C_1(n, p) r^{1-n/p} X^{-1/p}(r/D)(I[u])^{1/p},
\qquad (5.26)
$$

for any $D \geq \operatorname{diam}(\Omega)$. To estimate $\Lambda_r(x)$ we return to the original function by $v(z) = |z|^{n/p-1} u(z)$, thus

$$
\Lambda_r(x) = \int_{B_r} \frac{|u(z)|}{|z||x-z|^{n-1}} \mathrm{d}z.
$$

Inserting (5.18) in $\Lambda_r(x)$, we obtain

$$
\Lambda_r(x) \leq C_2(n, p)(I[u])^{1/p} \int_{B_r} \frac{X^{-1/p}(|z|/D)}{|z|^{n/p}|x-z|^{n-1}} \mathrm{d}z,
$$

for any $D \geq e^\Theta \operatorname{diam}(\Omega)$. Using Hölder's inequality with conjugate exponents $Q$ and $Q' = Q/(Q - 1)$, where $n < Q < p$ is fixed but depending only on $n, p$, we deduce

$$
\Lambda_r(x) \leq C_2(n, p)(I[u])^{1/p} \left( \int_{B_r} \frac{X^{-Q/p}(|z|/D)}{|z|^{nQ/p}} \mathrm{d}z \right)^{1/Q} \left( \int_{B_r} \frac{1}{|x-z|^{(n-1)Q'}} \mathrm{d}z \right)^{1/Q'}.
$$

Note that both integrals above are finite since $nQ/p < n$ if and only if $Q < p$, and $(n-1)Q' < n$ if and only if $n < Q$. Further, both integrals increase if we integrate over $B_r(0)$ and $B_r(x)$ respectively. Therefore

$$\Lambda_r(x) \leq C_2(n,p)(I[u])^{1/p} \left( \int_{B_r(0)} \frac{X^{-Q/p}(|z|/D)}{|z|^{nQ/p}} \mathrm{d}z \right)^{1/Q} \left( \int_{B_r(x)} \frac{1}{|x-z|^{(n-1)Q'}} \mathrm{d}z \right)^{1/Q'}$$

$$= C_3(n,p)(I[u])^{1/p} \left( \int_0^r t^{n-1-nQ/p} X^{-Q/p}(t/D) \mathrm{d}t \right)^{1/Q} r^{n/Q'-n+1}, \qquad (5.27)$$

for any $D \geq e^\Theta \operatorname{diam}(\Omega)$. Lemma A.1-(i) for $\alpha = n-1-nQ/p$ and $\beta = Q/p$ ensures the existence of constants $\eta \geq 0$ and $c > 0$ both depending only on $n, p, Q$ (and thus only on $n, p$), such that

$$\int_0^r t^{n-1-nQ/p} X^{-Q/p}(t/D) \mathrm{d}t \leq c r^{n-nQ/p} X^{-Q/p}(r/D),$$

for any $D \geq e^\eta \operatorname{diam}(\Omega)$. Thus (5.27) becomes

$$\Lambda_r(x) \leq C_4(n,p) r^{1-n/p} X^{-1/p}(r/D) (I[u])^{1/p}, \qquad (5.28)$$

for any $D \geq e^{\Theta'} \operatorname{diam}(\Omega)$ where $\Theta' = \max\{\Theta, \eta\}$.

Summarizing, in view of (5.26) and (5.28), estimate (5.25) becomes

$$J(x) \leq C_5(n,p) r^{1-n/p} X^{-1/p}(r/D) (I[u])^{1/p},$$

for any $D \geq e^{\Theta'} \operatorname{diam}(\Omega)$. The same estimate holds for $J(y)$ and thus (6.22) becomes

$$|u(x) - u(y)| \leq C_6(n,p) r^{1-n/p} X^{-1/p}(r/D) (I[u])^{1/p},$$

for any $D \geq e^{\Theta'} \operatorname{diam}(\Omega)$. The proof of (1.28) in case $n \geq 2$ is completed with $B = e^{\Theta'}$. ∎

## 5.5 Optimality of the logarithmic correction

In this section we prove the optimality of the exponent $1/p$ on $X$, in the Hölder semi-norm inequality (1.28) of Theorem VI. Note that we can pick one point in (1.28) to be the origin, and therefore it is enough to prove the alleged optimality in (5.18) and (5.8).

We consider the radially symmetric, Lipschitz continuous function
$$u_\delta(x) = \begin{cases} (\delta^2|x|)^H(6 - \frac{\log|x|}{\log\delta}), & \delta^6 \leq |x| < \delta^5 \\ (\delta^{-3}|x|^2)^H, & \delta^5 \leq |x| < \delta^4 \\ (\delta|x|)^H(1 + H\log(|x|/\delta^4)), & \delta^4 \leq |x| < \delta^3 \\ (\delta|x|)^H(1 - H\log(|x|/\delta^2)), & \delta^3 \leq |x| < \delta^2 \\ \delta^{3H}, & \delta^2 \leq |x| < \delta \\ (\delta^2|x|)^H \frac{\log|x|}{\log\delta}, & \delta \leq |x| \leq 1 \end{cases}$$
where $0 < \delta < 1$ and $H := (p-n)/p$ with $p > n \geq 1$. With $u_\delta$ we associate the quotient
$$Q_\epsilon[u_\delta; x] := \frac{(I[u_\delta])^{1/p}}{|u_\delta(x)||x|^{-1+n/p}X^{1/p-\epsilon}(|x|/D)}, \quad 0 \leq \epsilon < 1/p, \quad \delta^6 < |x| < 1.$$
Note that due to (5.18) and (5.8) (and after an approximation of $u_\delta$ by smooth functions), we have $Q_0[u_\delta; x] \geq C$, for some positive constant $C = C(n,p)$. To prove that the exponent $1/p$ on the correction weight $X$ cannot be decreased, we fix $0 < \epsilon < 1/p$ and taking $x$ such that $|x| = \delta^3$ we will prove that $Q_\epsilon[u_\delta; \delta^3] \to 0$ as $\delta \downarrow 0$.

We begin by computing $I[u_\delta]$. Setting $A_k := \{x \in \mathbb{R}^n : \delta^k < |x| < \delta^{k-1}\}$, $k = 1, ..., 6$, we have $I[u_\delta] = \sum_{k=1}^{6} I[u_\delta](A_k)$, where
$$I[u_\delta](A_k) := \int_{A_k} |\nabla u_\delta(x)|^p dx - H^p \int_{A_k} \frac{|u_\delta(x)|^p}{|x|^p} dx, \quad k = 1, ..., 6.$$
By the fact that $u_\delta$ is radially symmetric, we may use polar coordinates to get
$$I[u_\delta](A_k) = n\omega_n \left[ \int_{\delta^k}^{\delta^{k-1}} |\tilde{u}'_\delta(t)|^p t^{n-1} dt - H^p \int_{\delta^k}^{\delta^{k-1}} |\tilde{u}_\delta(t)|^p t^{n-1-p} dt \right], \quad k = 1, ..., 6,$$
where $\tilde{u}_\delta(t) = u_\delta(x)$ with $t = |x|$. We then have
$$I[u_\delta](A_1) = n\omega_n \frac{\delta^{2pH}}{\log^p(1/\delta)} \left[ \int_\delta^1 t^{-1}|1 - H\log(1/t)|^p dt - \int_\delta^1 t^{-1}(H\log(1/t))^p dt \right]$$
$$= n\omega_n \frac{\delta^{2pH}}{\log^p(1/\delta)} \left[ \int_\delta^{e^{-1/H}} + \int_{e^{-1/H}}^1 t^{-1}|1 - H\log(1/t)|^p dt - \int_\delta^1 t^{-1}(H\log(1/t))^p dt \right]$$
$$= \frac{n\omega_n}{(p+1)H} \delta^{2pH} \log(1/\delta) \left[ \left(H - \frac{1}{\log(1/\delta)}\right)^{p+1} + \left(\frac{1}{\log(1/\delta)}\right)^{p+1} - H^{p+1} \right],$$
where (since we will let $\delta \downarrow 0$) we have assumed $\delta < e^{-1/H}$ in order to get rid of the absolute value. Now we compute
$$I[u_\delta](A_6) = n\omega_n \frac{\delta^{2pH}}{\log^p(1/\delta)} \left[ \int_{\delta^6}^{\delta^5} t^{-1}(1 + H\log(t/\delta^6))^p dt - \int_{\delta^6}^{\delta^5} t^{-1}(H\log(t/\delta^6))^p dt \right]$$
$$= \frac{n\omega_n}{(p+1)H} \delta^{2pH} \log(1/\delta) \left[ \left(H + \frac{1}{\log(1/\delta)}\right)^{p+1} - \left(\frac{1}{\log(1/\delta)}\right)^{p+1} - H^{p+1} \right].$$

Thus $I[u_\delta](A_1) + I[u_\delta](A_6) =$

$$\frac{2n\omega_n}{(p+1)H}\delta^{2pH}\log(1/\delta)\left[\left(H + \frac{1}{\log(1/\delta)}\right)^{p+1} + \left(H - \frac{1}{\log(1/\delta)}\right)^{p+1} - 2H^{p+1}\right].$$

The factor in the square brackets is of order $o(1)$, as $\delta \downarrow 0$. Since $H = (p-n)/p$, we get

$$I[u_\delta](A_1) + I[u_\delta](A_6) = o(\delta^{2(p-n)}\log(1/\delta)), \quad \text{as } \delta \downarrow 0. \tag{5.29}$$

Similarly,

$$\begin{aligned} I[u_\delta](A_2) &= -n\omega_n H^p \delta^{3pH} \int_{\delta^2}^{\delta} t^{-pH-1} dt \\ &= -\frac{n\omega_n}{p} H^{p-1} \delta^{pH}(1 - \delta^{pH}), \end{aligned}$$

and

$$\begin{aligned} I[u_\delta](A_5) &= n\omega_n H^p \delta^{-3pH}\left[2^p \int_{\delta^5}^{\delta^4} t^{pH-1} dt - \int_{\delta^5}^{\delta^4} t^{pH-1} dt\right] \\ &= \frac{n\omega_n}{p} H^{p-1}(2^p - 1)\delta^{pH}(1 - \delta^{pH}). \end{aligned}$$

Hence

$$\begin{aligned} I[u_\delta](A_2) + I[u_\delta](A_5) &= \frac{n\omega_n}{p} H^{p-1}(2^p - 2)\delta^{pH}(1 - \delta^{pH}) \\ &= O(\delta^{p-n}), \quad \text{as } \delta \downarrow 0. \end{aligned} \tag{5.30}$$

Finally, the first summand of the last pair is

$$\begin{aligned} I[u_\delta](A_3) &= n\omega_n H^p \delta^{pH}\left[\int_{\delta^3}^{\delta^2}\left(-H\log\frac{t}{\delta^2}\right)^p t^{-1} dt - \int_{\delta^3}^{\delta^2}\left(1 - H\log\frac{t}{\delta^2}\right)^p t^{-1} dt\right] \\ &= \frac{n\omega_n}{p+1} H^{2p}\delta^{pH}\left(\log\frac{1}{\delta}\right)^{p+1}\left[1 + \frac{1}{(H\log\frac{1}{\delta})^{p+1}} - \left(1 + \frac{1}{H\log\frac{1}{\delta}}\right)^{p+1}\right], \end{aligned}$$

and the second one

$$I[u_\delta](A_4) = n\omega_n H^p \delta^{pH}\left[\int_{\delta^4}^{\delta^3}\left(2 + H\log\frac{t}{\delta^4}\right)^p t^{-1} dt - \int_{\delta^4}^{\delta^3}\left(1 + H\log\frac{t}{\delta^4}\right)^p t^{-1} dt\right]$$

$$= \frac{n\omega_n}{p+1} H^{2p}\delta^{pH}\left(\log\frac{1}{\delta}\right)^{p+1}\left[\left(1 + \frac{2}{H\log\frac{1}{\delta}}\right)^{p+1} - \left(1 + \frac{1}{H\log\frac{1}{\delta}}\right)^{p+1} + \frac{1 - 2^{p+1}}{(H\log\frac{1}{\delta})^{p+1}}\right].$$

Adding, we find $I[u_\delta](A_3) + I[u_\delta](A_4) =$

$$\frac{n\omega_n}{p+1} H^{2p}\delta^{pH}\left(\log\frac{1}{\delta}\right)^{p+1}\left[1 + \left(1 + \frac{2}{H\log\frac{1}{\delta}}\right)^{p+1} - 2\left(1 + \frac{1}{H\log\frac{1}{\delta}}\right)^{p+1} - \frac{2^{p+1} - 2}{(H\log\frac{1}{\delta})^{p+1}}\right].$$

The factor in square brackets is of order $\frac{p(p+1)}{H^2}o(\frac{1}{\log^2(1/\delta)})$, as $\delta \downarrow 0$, and we get

$$I[u_\delta](A_3) + I[u_\delta](A_4) = pn\omega_n H^{2(p-1)}\delta^{p-n}(\log(1/\delta))^{p-1} + o(\delta^{p-n}(\log(1/\delta))^{p-1}), \quad (5.31)$$

as $\delta \downarrow 0$. From (5.29), (5.30) and (5.31), we see that the leading term in $I[u_\delta]$, comes from $I[u_\delta](A_3) + I[u_\delta](A_4)$. More precisely

$$I[u_\delta] = pn\omega_n H^{2(p-1)}\delta^{p-n}(\log(1/\delta))^{p-1} + o(\delta^{p-n}(\log(1/\delta))^{p-1}),$$

as $\delta \downarrow 0$. Finally, we compute the denominator of $Q_\epsilon[u_\delta; \delta^3]$

$$\begin{aligned}|u_\delta(\delta^3)|\delta^{3(-1+n/p)}X^{1/p-\epsilon}(\delta^3/D) &= \delta^{1-n/p}(1+H\log(1/\delta))X^{1/p-\epsilon}(\delta^3/D)\\ &= \delta^{1-n/p}(1+H\log(1/\delta))(1-\log(\delta^3/D))^{-1/p+\epsilon}\\ &= 3H\delta^{1-n/p}(\log(1/\delta))^{1-1/p+\epsilon} + o(\delta^{1-n/p}(\log(1/\delta))^{1-1/p}),\end{aligned}$$

as $\delta \downarrow 0$. Dividing the two endmost relations we conclude

$$\begin{aligned}Q_\epsilon[u_\delta; \delta^3] &= \frac{1}{3}(pn\omega_n)^{1/p}H^{1-2/p}\frac{[\delta^{p-n}(\log 1/\delta)^{p-1} + o(\delta^{p-n}(\log(1/\delta))^{p-1})]^{1/p}}{\delta^{1-n/p}(\log(1/\delta))^{1-1/p+\epsilon} + o(\delta^{1-n/p}(\log(1/\delta))^{1-1/p})}\\ &= \frac{1}{3}(pn\omega_n)^{1/p}H^{1-2/p}\frac{[1+o(1)]^{1/p}}{(\log(1/\delta))^\epsilon + o(1)}\\ &\to 0, \quad \text{as } \delta \downarrow 0. \quad \blacksquare\end{aligned}$$

# Chapter 6

# Hardy-Sobolev type inequalities for $p > n$ - Distance to the boundary

In this chapter we will prove Theorems VIII, IX and X stated in the introduction. The proofs appear here for the first time. In this chapter we set

$$I[u] := \int_\Omega |\nabla u|^p \mathrm{d}x - \Big(\frac{p-1}{p}\Big)^p \int_\Omega \frac{|u|^p}{d^p}\mathrm{d}x; \quad p > 1,$$

whenever $u \in C_c^\infty(\Omega)$. Recall that $d \equiv d(x) := \mathrm{dist}(x, \mathbb{R}^n \setminus \Omega); x \in \Omega$.

## 6.1 The Hardy-Sobolev inequality for $p > n \geq 1$

In this section we prove of Theorem VIII. For $n = 1$, by (2.26) we have

$$u(x) \leq \frac{1}{2}\int_0^R |u'|\mathrm{d}t$$

$$(\text{setting } u(t) = (d(t))^{1-1/p}v(t)) \leq \frac{1}{2}\int_0^R d^{1-1/p}|v'|\mathrm{d}t + \frac{p-1}{2p}\int_0^R d^{-1/p}|v|\mathrm{d}t$$

$$\leq \int_0^R d^{1-1/p}|v'|\mathrm{d}t,$$

by (2.2) for $q = 1$ and $s = 1/p$. The proof follows applying Hölder's inequality and using (2.13) for $p \geq 2$ or (2.15) if $1 < p < 2$. We omit further details.

Assume next that $n \geq 2$. From Lemma 2.16 we get

$$u(x) \leq \frac{1}{n\omega_n}\int_\Omega \frac{|\nabla u(z)|}{|x-z|^{n-1}}\mathrm{d}z.$$

for $x \in \Omega$. Setting $u(z) = (d(z))^{1-1/p}v(z)$ we arrive at

$$n\omega_n|u(x)| \leq \underbrace{\int_\Omega \frac{(d(z))^{1-1/p}|\nabla v(z)|}{|x-z|^{n-1}}\mathrm{d}z}_{=:K(x)} + \frac{p-1}{p}\underbrace{\int_\Omega \frac{|v(z)|}{(d(z))^{1/p}|x-z|^{n-1}}\mathrm{d}z}_{=:\Lambda(x)}. \qquad (6.1)$$

We first estimate $K(x)$. Using Hölder's inequality we get

$$K(x) \leq \left(\underbrace{\int_\Omega \frac{1}{|x-z|^{(n-1)p/(p-1)}} \mathrm{d}z}_{=:M(x)}\right)^{1-1/p} \left(\int_\Omega d^{p-1}|\nabla v|^p \mathrm{d}z\right)^{1/p}. \tag{6.2}$$

By (5.3) and using (2.13) in the second factor of (6.2) we get

$$K(x) \leq C_1(n,p)[\mathcal{L}^n(\Omega)]^{1/n-1/p}(I[u])^{1/p}. \tag{6.3}$$

Next we bound $\Lambda(x)$. Using Hölder's inequality with conjugate exponents $p/(p-1-\varepsilon)$ and $p/(1+\varepsilon)$, where $0 < \varepsilon < (p-n)/p$ is fixed and depending only on $n,p$, we get

$$\Lambda(x) \leq \left(\underbrace{\int_\Omega \frac{1}{|x-z|^{(n-1)p/(p-1-\varepsilon)}} \mathrm{d}z}_{=:M(\varepsilon,x)}\right)^{1-(1+\varepsilon)/p} \left(\int_\Omega \frac{|v|^{p/(1+\varepsilon)}}{d^{1/(1+\varepsilon)}} \mathrm{d}z\right)^{(1+\varepsilon)/p}. \tag{6.4}$$

By (5.6) and using inequality (2.8) with $s = 1/(1+\varepsilon)$ and $q = p/(1+\varepsilon)$ in the second factor of (6.4), we obtain

$$\Lambda(x) \leq C_2(n,p)[\mathcal{L}^n(\Omega)]^{1/n-1/p-\varepsilon/p} \bigg[ \left(\frac{p}{\varepsilon}\right)^{p/(1+\varepsilon)} \int_\Omega d^{(p-1)/(1+\varepsilon)}|\nabla v|^{p/(1+\varepsilon)} \mathrm{d}z$$

$$+ \frac{p}{\varepsilon} \int_\Omega d^{\varepsilon/(1+\varepsilon)}|v|^{p/(1+\varepsilon)}(-\Delta d)\mathrm{d}z \bigg]^{(1+\varepsilon)/p}.$$

Using once more Hölder's inequality with conjugate exponents $1 + 1/\varepsilon$ and $1 + \varepsilon$ in both terms inside brackets, we get

$$\Lambda(x) \leq C_3(n,p)[\mathcal{L}^n(\Omega)]^{1/n-1/p-\varepsilon/p} \bigg[ [\mathcal{L}^n(\Omega)]^{\varepsilon/(1+\varepsilon)} \left(\int_\Omega d^{p-1}|\nabla v|^p \mathrm{d}z\right)^{1(1+\varepsilon)}$$

$$+ \left(\int_\Omega d(-\Delta d)\mathrm{d}z\right)^{\varepsilon/(1+\varepsilon)} \left(\int_\Omega |v|^p(-\Delta d)\mathrm{d}z\right)^{1/(1+\varepsilon)} \bigg]^{(1+\varepsilon)/p}.$$

Integration by parts shows that $\int_\Omega d(-\Delta d)\mathrm{d}z = \mathcal{L}^n(\Omega)$. Thus

$$\Lambda(x) \leq C_3(n,p)[\mathcal{L}^n(\Omega)]^{1/n-1/p} \bigg[ \left(\int_\Omega d^{p-1}|\nabla v|^p \mathrm{d}z\right)^{\frac{1}{1+\varepsilon}} + \left(\int_\Omega |v|^p(-\Delta d)\mathrm{d}z\right)^{\frac{1}{1+\varepsilon}} \bigg]^{\frac{1+\varepsilon}{p}}$$

$$\leq C_4(n,p)[\mathcal{L}^n(\Omega)]^{1/n-1/p} \left(\int_\Omega d^{p-1}|\nabla v|^p \mathrm{d}z + \int_\Omega |v|^p(-\Delta d)\mathrm{d}z\right)^{1/p},$$

where the last inequality follows by the fact that $a^{1/q} + b^{1/q} \leq 2^{1-1/q}(a+b)^{1/q}$ for all $q > 1$ and $a, b$ nonnegative. By (2.13) we conclude

$$\Lambda(x) \leq C_5(n,p)[\mathcal{L}^n(\Omega)]^{1/n-1/p}(I[u])^{1/p}. \tag{6.5}$$

The proof follows inserting (6.5) and (6.3) in (6.1). ∎

## 6.2 An optimal Hardy-Morrey inequality for $n = 1$

We present it in this separate section the proof of Theorem X. Firstly note that it suffices to restrict ourselves to the case $\Omega = (0, 1)$. We start with the following lemma

**Lemma 6.1.** *Let $q > 1$, $\beta > 1 - q$. There exists positive constant $c = c(q, \beta)$ such that, for all $v \in C_c^\infty(0, 1)$ and any $D \geq 1/2$*

$$\sup_{x \in (0,1)} \left\{ |v(x)| X^{(\beta+q-1)/q}(d(x)/D) \right\} \leq c \left( \int_0^1 d^{q-1} |v'|^q X^\beta(d/D) \mathrm{d}t + |v(1/2)|^q \right)^{1/q}.$$

**Proof.** Letting $D \geq 1/2$ we have for any $x \in (0, 1)$

$$\begin{aligned} v(x) &= -\int_x^{1/2} v' \mathrm{d}t + v(1/2) \\ &\leq \left| \int_x^{1/2} d^{-1} X^{-\beta/(q-1)}(d/D) \mathrm{d}t \right|^{1-1/q} \left| \int_x^{1/2} d^{q-1} |v'|^q X^\beta(d/D) \mathrm{d}t \right|^{1/q} + |v(1/2)|, \end{aligned}$$

by Hölder's inequality with conjugate exponents $q$ and $q/(q-1)$. Since

$$[X^{-1-\beta/(q-1)}(d(t)/D)]' = \left( -1 - \frac{\beta}{q-1} \right) (d(t))^{-1} X^{-\beta/(q-1)}(d(t)/D) d'(t),$$

recalling that $d'(t) = 1$ if $t \in (0, 1/2)$ and $d'(t) = -1$ if $t \in (1/2, 1)$, the left integral above can be computed. We find

$$\left| \int_x^{1/2} d^{-1} X^{-\beta/(q-1)}(d/D) \mathrm{d}t \right| = \frac{q-1}{\beta+q-1} |X^{-1-\beta/(q-1)}(1/2D) - X^{-1-\beta/(q-1)}(d(x)/D)|,$$

for all $x \in (0, 1)$. In addition, the second integral satisfies

$$\left| \int_x^{1/2} d^{q-1} |v'|^q X^\beta(d/D) \mathrm{d}t \right| \leq \int_0^1 d^{q-1} |v'|^q X^\beta(d/D) \mathrm{d}t,$$

for all $x \in (0, 1)$. Coupling the above estimates we arrive at

$$\begin{aligned} |v(x)| \\ \leq c \Big| X^{-1-\beta/(q-1)}(1/2D) - X^{-1-\beta/(q-1)}(d(x)/D) \Big|^{1-1/q} & \left( \int_0^1 d^{q-1} |v'|^q X^\beta(d/D) \mathrm{d}t \right)^{1/q} \\ & + |v(1/2)|, \end{aligned}$$

where $c \equiv c(q, \beta) = [(q-1)/(\beta + q - 1)]^{1-1/q}$, or

$$|v(x)|X^{(\beta+q-1)/q}(d(x)/D)$$
$$\leq c \Big|\Big(\frac{X(d(x)/D)}{X(1/2D)}\Big)^{1+\beta/(q-1)} - 1\Big|^{1-1/q} \Big(\int_0^1 d^{q-1}|v'|^q X^\beta(d/D)\mathrm{d}t\Big)^{1/q}$$
$$+ |v(1/2)|,$$

where we have also used the fact that $X^{(\beta+q-1)/q}(d(x)/D) \leq 1$ for all $x \in (0,1)$, since $\beta > 1 - q$. Next, since $d(x) \leq 1/2$ in $(0, 1)$ we have

$$\Big|\Big(\frac{X(d(x)/D)}{X(1/2D)}\Big)^{1+\beta/(q-1)} - 1\Big| \leq 1.$$

Hence

$$|v(x)|X^{(\beta+q-1)/q}(d(x)/D) \leq c\Big(\int_0^1 d^{q-1}|v'|^q X^\beta(d/D)\mathrm{d}t\Big)^{1/q} + |v(1/2)|.$$

The desired inequality follows using $a^{1/q} + b^{1/q} \leq 2^{1-1/q}(a+b)^{1/q}$, for all $q > 1$ and $a, b \geq 0$. ∎

Next we prove the counterpart of Proposition 5.2.

**Proposition 6.2.** *Let $p > 1$. There exists positive constant $c$, depending only on $p$, such that for all $u \in C_c^\infty(0,1)$ and any $D \geq 1/2$*

$$\sup_{x \in (0,1)} \Big\{\frac{|u(x)|}{(d(x))^{1-1/p}} X^{1/p}(d(x)/D)\Big\} \leq c(I[u])^{1/p}. \tag{6.6}$$

**Proof.** We set $u(x) = d(x)^{1-1/p}v(x)$. If $1 < p < 2$, by Lemma 6.1 for $q = p$ and $\beta = 2 - p$ we have

$$|v(x)|X^{1/p}(d(x)/D) \leq c\Big(\int_0^1 d^{p-1}|v'|^p X^{2-p}(d/D)\mathrm{d}t + |v(1/2)|^p\Big)^{1/p},$$

for any $D \geq 1/2$. The result follows by (2.15). If $p \geq 2$, by Lemma 6.1 for $q = 2$ and $\beta = 0$ we have

$$|w(x)|X^{1/2}(d(x)/D) \leq c\Big(\int_0^1 d|w'|^2 \mathrm{d}t + |w(1/2)|^2\Big)^{1/2},$$

for any $w \in C_c^\infty(0,1)$ and any $D \geq 1/2$. For $w(x) = |v(x)|^{p/2}$ we obtain

$$|v(x)|X^{1/p}(d(x)/D) \leq c\Big(\int_0^1 d|v|^{p-2}|v'|^2 \mathrm{d}t + |v(1/2)|^p\Big)^{1/p}.$$

The result follows by (2.14). ∎

Now we use (6.6) to obtain its counterpart inequality with the exact Hölder seminorm.

**Proof of Theorem X.** For $0 < y < x < 1$ we have

$$|u(x) - u(y)| = \left| \int_y^x u' \, dt \right|$$

$$\text{(setting } u(t) = (d(t))^{1-1/p} v(t)) \leq \underbrace{\int_y^x d^{1-1/p} |v'| \, dt}_{:= K(x,y)} + \frac{p-1}{p} \underbrace{\int_y^x d^{1/p} |v| \, dt}_{:= \Lambda(x,y)}. \quad (6.7)$$

where we have used the fact that $|d'(t)| = 1$ a.e. in $(0, 1)$. To estimate $K(x, y)$ we use Hölder's inequality to get

$$K(x, y) \leq (x-y)^{1-1/p} \left( \int_y^x d^{p-1} |v'|^p \, dt \right)^{1/p}$$

$$\text{(by (2.13))} \leq c_1(p)(x-y)^{1-1/p} (I[u])^{1/p}$$

$$\leq c_1(p)(x-y)^{1-1/p} X^{-1/p}((x-y)/D)(I[u])^{1/p}, \quad (6.8)$$

for any $D \geq 1$ since $0 \leq X(t) \leq 1$ for all $t \in (0, 1]$. To estimate $\Lambda(x, y)$ we return to the original variable by $v(t) = (d(t))^{1/p-1} u(t)$. Thus

$$\Lambda(x, y) = \int_y^x \frac{|u|}{d} \, dt.$$

Inserting (6.6) in $\Lambda(x, y)$ we obtain

$$\Lambda(x, y) \leq c_2(p)(I[u])^{1/p} \int_y^x d^{-1/p} X^{-1/p}(d/D) \, dt$$

$$\leq c_3(p)(x-y)^{1-1/p} X^{-1/p}((x-y)/D)(I[u])^{1/p}, \quad (6.9)$$

by virtue of Lemma A.1-(ii) in the Appendix. Coupling (6.8) and (6.9) with (6.7) we end up with

$$|u(x) - u(y)| \leq c_4(p)(x-y)^{1-1/p} X^{-1/p}((x-y)/D)(I[u])^{1/p},$$

for all $0 < y < x < 1$ and any $D \geq 1$, which is the desired estimate.

The proof of the optimality of the exponent $1/p$ on $X$ follows by considering the function $u_\delta$ as defined in §5.5 for sufficiently small $\delta$, so that the distance in (6.6) is taken only from $0$. The computations then are identical to §5.5. ∎

## 6.3   A Hardy-Morrey inequality for $n \geq 2$ and $\Omega = B_R$

In this section we prove Theorem IX. We start with

**Proposition 6.3.** *Letting $p > n \geq 2$, there exist constants $\Theta \geq 0$ and $C > 0$ both depending only on $n, p$, such that for all $u \in C_c^\infty(B_R)$ and any $D \geq e^\Theta R$*

$$\sup_{x \in B_R} \left\{ \frac{|u(x)|}{(d(x))^{1-n/p}} X^{1/p}\left(\frac{d(x)}{D}\right) \right\} \leq C \Big[ \int_{B_R} |\nabla u|^p \mathrm{d}x - \Big(\frac{p-1}{p}\Big)^p \int_{B_R} \frac{|u|^p}{d^p} \mathrm{d}x \Big]^{1/p}, \quad (6.10)$$

*where $X(t) = (1 - \log t)^{-1}$, $t \in (0, 1]$.*

**Proof.** Letting $x \in B_R$ we pick $\xi \in \partial B_R$ such that $r := d(x) = |x - \xi|$ and we consider the ball $B_r(\xi)$. For $u \in C_c^\infty(B_R)$ we set

$$u_{B_r(\xi)} = \frac{1}{\omega_n r^n} \int_{B_r(\xi)} u(z) \mathrm{d}z.$$

Letting $y \in B_r(\xi) \cap B_R$ we get from Lemma 2.19

$$|u(y) - u_{B_r(\xi)}| \leq \frac{2^n}{n \omega_n} \int_{B_r(\xi)} \frac{|\nabla u(z)|}{|y - z|^{n-1}} \mathrm{d}z.$$

Setting $u(z) = (d(z))^{1-1/p} v(z)$ we arrive at

$$\frac{n \omega_n}{2^n} |u(y) - u_{B_r}| \leq \underbrace{\int_{B_r(\xi)} \frac{(d(z))^{1-1/p} |\nabla v(z)|}{|y - z|^{n-1}} \mathrm{d}z}_{=: K_r(y)}$$
$$+ \frac{p-1}{p} \underbrace{\int_{B_r(\xi)} \frac{|v(z)|}{(d(z))^{1/p} |y - z|^{n-1}} \mathrm{d}z}_{=: \Lambda_r(y)}. \quad (6.11)$$

We will derive suitable bounds for $K_r(y), \Lambda_r(y)$. For $K_r(y)$ we use Hölder's inequality

$$K_r(y) \leq \left( \int_{B_r(\xi)} \frac{1}{|y - z|^{(n-1)p/(p-1)}} \mathrm{d}z \right)^{1-1/p} \left( \int_{B_r(\xi)} d^{p-n} |\nabla v|^p \mathrm{d}z \right)^{1/p}.$$

Both integrals increase if we integrate over $B_r(y)$ and $B_R$ respectively. Thus

$$K_r(y) \leq \left( \int_{B_r(y)} \frac{1}{|y - z|^{(n-1)p/(p-1)}} \mathrm{d}z \right)^{1-1/p} \left( \int_{B_R} d^{p-n} |\nabla v|^p \mathrm{d}z \right)^{1/p}.$$

Computing the first integral and using (2.13) for the second we arrive at

$$\begin{aligned} K_r(x) &\leq C_1(n,p) r^{1-n/p} (I[u])^{1/p} \\ &\leq C_1(n,p) r^{1-n/p} X^{-1/p}(r/D) (I[u])^{1/p}, \end{aligned} \qquad (6.12)$$

for any $D \geq R$ since $X(t) \leq 1; t \in (0,1]$.

Next we bound $\Lambda_r(y)$. Using Hölder's inequality with conjugate exponents $p/(p-1-\varepsilon)$ and $p/(1+\varepsilon)$, where $0 < \varepsilon < (p-n)/n$ is fixed but depending only on $n, p$, we obtain

$$\Lambda_r(y) \leq \left( \underbrace{\int_{B_r(\xi)} \frac{1}{|y-z|^{(n-1)p/(p-1-\varepsilon)}} dz}_{=:M_r(y)} \right)^{1-(1+\varepsilon)/p} \left( \int_{B_r(\xi)} \frac{|v|^{p/(1+\varepsilon)}}{d^{1/(1+\varepsilon)}} dz \right)^{(1+\varepsilon)/p} \qquad (6.13)$$

By (5.6) and using (2.6) for $V = B_r(\xi)$, $s = 1/(1+\varepsilon)$ and $q = p/(1+\varepsilon)$, we get from (6.13)

$$\begin{aligned} \Lambda_r(x) \leq\ & C_1(n,p) r^{(p-n-n\varepsilon)/p} \Bigg( \underbrace{\int_{B_r(\xi)} d^{(p-1)/(1+\varepsilon)} |\nabla v|^{p/(1+\varepsilon)} dz}_{=:N_r} \\ & + \underbrace{\int_{\partial B_r(\xi)} d^{\varepsilon/(1+\varepsilon)} |v|^{p/(1+\varepsilon)} \nabla d \cdot \nu \, dS_z}_{=:P_r} \\ & + \underbrace{\int_{B_r(\xi)} d^{\varepsilon/(1+\varepsilon)} |v|^{p/(1+\varepsilon)} (-\Delta d) dz}_{=:Q_r} \Bigg)^{(1+\varepsilon)/p}, \end{aligned} \qquad (6.14)$$

where $\nu$ is the outward pointing unit normal vector field along $\partial B_r$. Now we need to estimate $N_r, P_r$ and $Q_r$. For $N_r$, we use Hölder's inequality with conjugate exponents $1 + 1/\varepsilon$ and $1 + \varepsilon$ to get

$$\begin{aligned} N_r &\leq (\omega_n r^n)^{\varepsilon/(1+\varepsilon)} \left( \int_{B_r} d^{p-1} |\nabla v|^p dz \right)^{1/(1+\varepsilon)} \\ \text{(by (2.13))} &\leq C_3(n,p) r^{n\varepsilon/(1+\varepsilon)} (I[u])^{1/(1+\varepsilon)} \\ &\leq C_3(n,p) r^{n\varepsilon/(1+\varepsilon)} X^{-1/(1+\varepsilon)}(r/D) (I[u])^{1/(1+\varepsilon)}, \end{aligned} \qquad (6.15)$$

for any $D \geq R$, since $X(t) \leq 1$ for all $t \in (0,1]$. Now set

$$\partial B_r^+(\xi) := \{ z \in B_R : |z - \xi| = r \text{ and } \nabla d(z) \cdot \nu(z) \geq 0 \}.$$

Then $P_r$ increases if we change the domain of integration from $\partial B_r(\xi)$ to $\partial B_r^+(\xi)$. Hence for any $D \geq R$

$$P_r \leq \int_{\partial B_r^+(\xi)} \left\{ d^{\varepsilon/(1+\varepsilon)} X^{-1/(1+\varepsilon)}(d/D) \right\} \left\{ |v|^{p/(1+\varepsilon)} X^{1/(1+\varepsilon)}(d/D) \right\} \nabla d \cdot \nu \mathrm{d} S_z$$

$$\leq \underbrace{\left( \int_{\partial B_r^+(\xi)} dX^{-1/\varepsilon}(d/D) \nabla d \cdot \nu \mathrm{d} S_z \right)}_{=:S_r}^{\varepsilon/(1+\varepsilon)}$$

$$\times \underbrace{\left( \int_{\partial B_r^+(\xi)} |v|^p X(d/D) \nabla d \cdot \nu \mathrm{d} S_z \right)}_{=:T_r}^{1/(1+\varepsilon)}, \quad (6.16)$$

where we have used once more Hölder's inequality with exponents $1 + 1/\varepsilon$ and $1 + \varepsilon$. Noting now that $\nu = z/r$, quantity $S_r$ can be computed:

$$\begin{aligned} S_r &\leq rX^{-1/\varepsilon}(r/D) \int_{\partial B_r^+(\xi)} |\nabla d \cdot \nu| \mathrm{d} S_z \\ &\leq n\omega_n r^n X^{-1/\varepsilon}(r/D). \end{aligned} \quad (6.17)$$

Due to the special geometry of the ball, we may find a supplemental surface $\sigma$ in which $\nabla d \cdot \nu \equiv 0$ and such that $\partial B_r^+(\xi) \cup \sigma$ enclose a set $A \subset \overline{B}_R$. Now $T_r$ may be estimated after an integration by parts. More precisely,

$$\begin{aligned} T_r &= \int_{\partial A} |v|^p X(d/D) \nabla d \cdot \nu \mathrm{d} S_z \\ &= \int_A \mathrm{div} \left\{ X(d/D) \nabla d \right\} |v|^p \mathrm{d} z + \int_A X(d/D) \nabla d \cdot \nabla(|v|^p) \mathrm{d} z. \end{aligned}$$

A simple calculation shows that $\mathrm{div}\{X(d/D)\nabla d\} = d^{-1} X^2(d/D) - X(d/D)(-\Delta d)$, in the sense of distributions. In the second integral we compute the gradient and note that $\nabla d \cdot \nabla |v| \leq |\nabla v|$, a.e. in $A$. Thus,

$$\begin{aligned} T_r &\leq \int_A \frac{|v|^p}{d} X^2(d/D) \mathrm{d} z - \int_A |v|^p X(d/D)(-\Delta d) \mathrm{d} z + p \int_A |v|^{p-1} |\nabla v| X(d/D) \mathrm{d} z \\ &\leq \int_{B_R} \frac{|v|^p}{d} X^2(d/D) \mathrm{d} z + p \int_{B_R} |v|^{p-1} |\nabla v| X(d/D) \mathrm{d} z, \end{aligned}$$

where se have used property $(\mathfrak{C})$. We return to the original function $u$ by recalling that $v(z) = (d(z))^{1/p-1} u(z)$ in the first integral, and rearrange the integrand in the second one:

$$T_r \leq \int_{B_R} \frac{|u|^p}{d^p} X^2(d/D) \mathrm{d} z + p \int_{B_R} \left\{ d^{1/2} |v|^{p/2-1} |\nabla v| \right\} \left\{ d^{-1/2} |v|^{p/2} X(d/D) \right\} \mathrm{d} z.$$

According to Theorem 2.11-(ii), there exist constants $\theta$ and $c > 0$, depending only on $n, p$, such that for any $D \geq Re^{\theta}$, the first term on the right hand side is bounded above by $cI[u]$. For the second term we may use Cauchy-Schwarz inequality. It follows that

$$\begin{aligned} T_r &\leq cI[u] + p\left(\int_{B_R} d|v|^{p-2}|\nabla v|^2 \mathrm{d}z\right)^{1/2} \left(\int_{B_R} d^{-1}|v|^p X^2(|z|/D)\mathrm{d}z\right)^{1/2} \\ &\leq cI[u] + p(c_{n,p}^{-1}I[u])^{1/2}(cI[u])^{1/2} \\ &= C_4(n,p)I[u], \end{aligned} \quad (6.18)$$

by (2.14) and Theorem 2.11-(ii) again. Setting $\theta^+ = \max\{0, \theta\}$, by (6.17) and (6.18), estimate (6.16) implies

$$P_r \leq C_5(n,p) r^{n\varepsilon/(1+\varepsilon)} X^{-1/(1+\varepsilon)}(r/D)(I[u])^{1/(1+\varepsilon)},$$

for any $D \geq Re^{\theta^+}$. To estimate $Q_r$, we apply Hölder's inequality as follows

$$\begin{aligned} Q_r &\leq \left(\int_{B_r(\xi)} d(-\Delta d)\mathrm{d}z\right)^{\varepsilon/(1+\varepsilon)} \left(\int_{B_r(\xi)} |v|^p(-\Delta d)\mathrm{d}z\right)^{1/(1+\varepsilon)} \\ &\leq c_p\left(\mathcal{L}^n(B_r(\xi) \cap B_R) - \int_{\partial B_r(\xi)} d\nabla d \cdot \nu \mathrm{d}S_z\right)^{\varepsilon/(1+\varepsilon)} \left(\int_{B_R} |v|^p(-\Delta d)\mathrm{d}z\right)^{1/(1+\varepsilon)} \\ &\leq c_p\left(\mathcal{L}^n(B_r(\xi)) + \int_{\partial B_r(\xi)} d\mathrm{d}S_z\right)^{\varepsilon/(1+\varepsilon)} (I[u])^{1/(1+\varepsilon)} \\ &\leq c_{n,p} r^{n\varepsilon/(1+\varepsilon)}(I[u])^{1/(1+\varepsilon)} \end{aligned} \quad (6.19)$$

Inserting (6.19), (6.15) and the above inequality in (6.14), we obtain

$$\begin{aligned} \Lambda_r(x) &\leq C_6(n,p) r^{1-n/p-n\varepsilon/p}\left(r^{n\varepsilon/(1+\varepsilon)} X^{-1/(1+\varepsilon)}(r/D)(I[u])^{1/(1+\varepsilon)}\right)^{(1+\varepsilon)/p} \\ &= C_6(n,p) r^{1-n/p} X^{-1/p}(r/D)(I[u])^{1//p}, \end{aligned} \quad (6.20)$$

for any $D \geq Re^{\theta^+}$.

Applying estimates (6.20) and (6.12) to estimate (6.11), we finally conclude

$$|u(x) - u_{B_r}| \leq C_7(n,p) r^{1-n/p} X^{-1/p}(r/D)(I[u])^{1/p}, \quad (6.21)$$

for all $x \in B_r(\xi) \cap B_R$ and any $D \geq Re^{\theta^+}$. It follows from (6.21) for $x = \xi$, that

$$|u_{B_r}| \leq C_7(n,p) r^{1-n/p} X^{-1/p}(r/D)(I[u])^{1/p}.$$

Thus,

$$\begin{aligned} |u(x)| &\leq |u(x) - u_{B_r(\xi)}| + |u_{B_r(\xi)}| \\ &\leq 2C_7(n,p) r^{1-n/p} X^{-1/p}(r/D)(I[u])^{1/p}, \end{aligned}$$

for all $x \in B_r(\xi) \cap B_R$ and any $D \geq Re^{\theta^+}$. Now if $x \in B_R$, consider a ball $B_r(\xi)$ of radius $r = d(x)$, centered at $\xi = x - r\nabla d(x)$. Then the previous inequality yields

$$|u(x)| \leq C(n,p)(d(x))^{1-n/p} X^{-1/p}(d(x)/D)(I[u])^{1/p}$$

for any $D \geq Re^{\theta^+}$. Rearranging and taking the supremum over all $x \in B_R$, the result follows. ∎

**Proof of Theorem IX.** Given any two points $x, y \in B_R$ with $x \neq y$, we note that either $|x-y| \geq \frac{d(x)}{2} \wedge \frac{d(y)}{2}$, or $|x-y| \leq \frac{d(x)}{2} \vee \frac{d(y)}{2}$. In the first case we have

$$\begin{aligned} |u(x) - u(y)| &\leq |u(x)| + |u(y)| \\ &\leq C\Big((d(x))^{1-n/p} X^{-1/p}(d(x)/D) + (d(y))^{1-n/p} X^{-1/p}(d(y)/D)\Big)(I[u])^{1/p}, \end{aligned}$$

for any $D \geq e^{\Theta} R$, where we have used (6.10) twice. For any $D \geq R$, the function $f(t) := t^{1-n/p} X^{-1/p}(t/D)$ is increasing in $(0, R)$, thus $d(x) \wedge d(y) \leq 2|x-y|$ implies

$$|u(x) - u(y)| \leq 2^{2-n/p} C|x-y|^{1-n/p} X^{-1/p}(|x-y|/D)(I[u])^{1/p},$$

for any $D \geq \max\{4R, e^{\Theta} R\}$.

It remains to consider the case where $|x-y| \leq \frac{d(x)}{2} \vee \frac{d(y)}{2}$. Let for example $|x-y| \leq \frac{d(x)}{2}$. Consider the ball $B_r(x)$ with $r := 3|x-y|/2$. We have

$$\begin{aligned} |u(x) - u(y)| &\leq |u(x) - u_{B_r(x)}| + |u(y) - u_{B_r(x)}| \\ &\leq \frac{2^n}{n\omega_n}\Big\{\underbrace{\int_{B_r(x)} \frac{|\nabla u(z)|}{|x-z|^{n-1}}dz}_{=:J(x)} + \underbrace{\int_{B_r(x)} \frac{|\nabla u(z)|}{|y-z|^{n-1}}dz}_{=:J(y)}\Big\}, \end{aligned} \quad (6.22)$$

where we have used Lemma 2.19 twice. To estimate $J(x)$ and $J(y)$, we let $\eta$ to be either $x$ or $y$. We start with the substitution $u(z) = (d(z))^{1-1/p} v(z)$, to get

$$J(\eta) \leq \underbrace{\int_{B_r(x)} \frac{(d(z))^{1-1/p}|\nabla v(z)|}{|\eta-z|^{n-1}}dz}_{=:K_r(\eta)} + \frac{p-1}{p}\underbrace{\int_{B_r(x)} \frac{|v(z)|}{(d(z))^{1/p}|\eta-z|^{n-1}}dz}_{=:\Lambda_r(\eta)}. \quad (6.23)$$

We will derive suitable bounds for $K_r(\eta), \Lambda_r(\eta)$. For $K_r(\eta)$ we use Hölder's inequality

$$K_r(\eta) \leq \Big(\int_{B_r(x)} |\eta-z|^{-\frac{(n-1)p}{p-1}}dz\Big)^{1-1/p} \Big(\int_{B_r(x)} d^{p-1}|\nabla v|^p dz\Big)^{1/p}.$$

If $\eta = y$ the first integral increases if we integrate over $B_r(y)$ instead of $B_r(x)$. Thus

$$\begin{aligned} K_r(\eta) &\leq \Big(\int_{B_r(\eta)} |\eta-z|^{-\frac{(n-1)p}{p-1}}dz\Big)^{1-1/p} \Big(\int_{B_R} d^{p-1}|\nabla v|^p dz\Big)^{1/p} \\ \text{(by (2.13))} &\leq \Big(n\omega_n \frac{p-1}{p-n}\Big)^{1-1/p} r^{1-n/p} (\tfrac{1}{c_p} I[u])^{1/p} \\ &\leq (3/2)^{1-n/p} C_1(n,p)|x-y|^{1-n/p} X^{-1/p}(|x-y|/D)(I[u])^{1/p}, \quad (6.24) \end{aligned}$$

for any $D \geq 2R$, where the last inequality follows by the fact that $0 < X(t) \leq 1$ for all $t \in (0,1]$. To estimate $\Lambda_r(\eta)$ we return to the original function by $v(z) = (d(z))^{1/p-1}u(z)$, thus

$$\Lambda_r(\eta) = \int_{B_r(x)} \frac{|u(z)|}{d(z)|\eta - z|^{n-1}} dz.$$

Inserting (6.10) in $\Lambda_r(\eta)$, to obtain

$$\Lambda_r(\eta) \leq C_2(n,p)(I[u])^{1/p} \int_{B_r(x)} \frac{(d(z))^{-n/p} X^{-1/p}(d(z)/D)}{|\eta - z|^{n-1}} dz, \qquad (6.25)$$

for any $D \geq e^{\Theta} R$. To estimate the above integral we note that if $z \in B_r(x)$ and $\xi_z \in \partial B_R$ realizes the distance of $z$ to the boundary $\partial B_R$, then we obtain

$$\begin{aligned}
d(z) &= |\xi_z - z| \\
\text{(triangle inequality)} &\geq |\xi_z - x| - |z - x| \\
(d(x) := \inf_{\xi \in \partial B_R} |\xi - x|) &\geq d(x) - |z - x| \\
(z \in B_r(x)) &\geq d(x) - r \\
\text{(by hypothesis: } d(x) \geq 2|x-y| = 4r/3) &\geq r/3.
\end{aligned}$$

Thus $d(z) \geq r/3$ for all $z \in B_r(x)$, and (6.25) becomes

$$\Lambda_r(\eta) \leq C_2(n,p) r^{-n/p} X^{-1/p}(r/D)(I[u])^{1/p} \int_{B_r(x)} \frac{1}{|\eta - z|^{n-1}} dz,$$

for any $D \geq \max\{6R, e^{\Theta} R\}$. If $\eta = y$ the last integral increases if we integrate over $B_r(y)$ instead of $B_r(x)$. Thus

$$\begin{aligned}
\Lambda_r(\eta) &\leq C_2(n,p) r^{-n/p} X^{-1/p}(r/D)(I[u])^{1/p} \int_{B_r(\eta)} \frac{1}{|\eta - z|^{n-1}} dz \\
&= n\omega_n C_2(n,p) r^{1-n/p} X^{-1/p}(r/D)(I[u])^{1/p} \\
&= n\omega_n (3/2)^{1-n/p} C_2(n,p) |x-y|^{1-n/p} X^{-1/p}(|x-y|/(2D/3))(I[u])^{1/p}, \quad (6.26)
\end{aligned}$$

for any $D \geq \max\{6R, e^{\Theta} R\}$. Inserting (6.24) and (6.26) in (6.23) we get

$$J(\eta) \leq C_3(n,p)|x-y|^{1-n/p} X^{-1/p}(|x-y|/D)(I[u])^{1/p},$$

for any $D \geq \max\{6R, e^{\Theta} R\}$. Estimate (6.22) now gives

$$|u(x) - u(y)| \leq \frac{2^{n+1}}{n\omega_n} C_3(n,p)|x-y|^{1-n/p} X^{-1/p}(|x-y|/D)(I[u])^{1/p},$$

for any $D \geq \max\{6R, e^{\Theta} R\}$. ∎

# Appendix

## Some calculus lemmas

The following technical fact concerns the auxiliary function $X(t) = (1 - \log t)^{-1}$, $t \in (0, 1]$.

**Lemma A.1.** *Let $\alpha > -1$ and $\beta, R > 0$. For all $r \in (0, R]$, all $c > 1/(\alpha + 1)$ and any $D \geq e^\eta R$, where $\eta := \max\{0, \frac{(\beta-\alpha-1)c+1}{(\alpha+1)c-1}\}$, we have*

(i) $\quad \displaystyle\int_0^r t^\alpha X^{-\beta}(t/D)\mathrm{d}t \leq cr^{\alpha+1}X^{-\beta}(r/D).$

*If $\alpha$ is restricted in $(-1, 0]$ then for all $0 \leq y \leq x \leq R$ we have*

(ii) $\quad \displaystyle\int_y^x t^\alpha X^{-\beta}(t/D)\mathrm{d}t \leq c(x-y)^{\alpha+1}X^{-\beta}((x-y)/D).$

**Proof.** Let $c > 0$ and $D \geq R$. We set

$$f(r) := \int_0^r t^\alpha X^{-\beta}(t/D)\mathrm{d}t - cr^{\alpha+1}X^{-\beta}(r/D), \quad r \in (0, R].$$

To prove (i) it suffices to show that for suitable values of the parameters $c$ and $D$, we have $f(r) \leq 0$ for all $r \in (0, R)$. We have $f(0+) = 0$ and thus it is enough to choose $c$ and $D$ in such a way so that $f$ is decreasing in $(0, R)$. To this end we compute

$$f'(r) = cr^\alpha X^{-\beta}(r/D)[1/c - (\alpha + 1) + \beta X(r/D)], \quad r \in (0, R].$$

It is easy to see that for $c > 1/(\alpha + 1)$, any $D \geq e^\eta R$ is such that $f'(r) \leq 0$ for all $r \in (0, R)$. To prove (ii) we note that since $f$ in decreasing, $0 \leq y \leq x \leq R$ implies $f(y) \geq f(x)$, and so

$$\begin{aligned}
\int_y^x t^\alpha X^{-\beta}(t/D)\mathrm{d}t &\leq c[x^{\alpha+1}X^{-\beta}(x/D) - y^{\alpha+1}X^{-\beta}(y/D)] \\
&\leq c(x^{\alpha+1} - y^{\alpha+1})X^{-\beta}(x/D) \\
&\leq c(x^{\alpha+1} - y^{\alpha+1})X^{-\beta}((x-y)/D),
\end{aligned}$$

where the last two inequalities follow since $X^{-\beta}(r/D)$ is decreasing in $(0, R)$. If $\alpha \in (-1, 0]$ then $x^{\alpha+1} - y^{\alpha+1} \leq (x - y)^{\alpha+1}$, and the result follows. ∎

The next lemma was used in §2.1.2.

**Lemma A.2.** *For any $p > 1$ and all $a, b \in \mathbb{R}^n$ we have:*

(i) *if $1 < p < 2$ then*
$$|a - b|^p - |a|^p \geq \frac{3p(p-1)}{16} \frac{|b|^2}{(|a - b| + |a|)^{2-p}} - p|a|^{p-2} a \cdot b.$$

(ii) *if $p \geq 2$ then*

(a) $|a - b|^p - |a|^p \geq \dfrac{1}{2^{p-1} - 1}|b|^p - p|a|^{p-2} a \cdot b,$

(b) $|a - b|^p - |a|^p \geq \dfrac{1}{2^{p-2}(2^{p-1} - 1)}|a|^{p-2}|b|^2 - p|a|^{p-2} a \cdot b.$

**Proof.** Parts (i) and (ii)-(a) are proved in detail in the Appendix of [Lndqv]. To prove (ii)-(b), if $|b| \geq |a|/2$ then it follows from (ii)-(a) that

$$|a - b|^p - |a|^p \geq \frac{1}{2^{p-2}(2^{p-1} - 1)}|a|^{p-2}|b|^2 - p|a|^{p-2} a \cdot b.$$

On the other hand, if $|b| < |a|/2$ then $|a - \xi b| \geq |a|/2$ for all $\xi \in (0, 1)$. Hence, taking the Taylor expansion of $f(t) = |a - bt|^p$ around $t = 0$ we have (for some $\xi \in (0, 1)$)

$$\begin{aligned}
|a - b|^p &= |a|^p - p|a|^{p-2} a \cdot b + \frac{p(p-2)}{2}|a - \xi b|^{p-4}((a - \xi b) \cdot b)^2 + p|a - \xi b|^{p-2}|b|^2 \\
&\geq |a|^p - p|a|^{p-2} a \cdot b + \frac{p}{2^{p-2}}|a|^{p-2}|b|^2.
\end{aligned}$$

The constant obtained in the case $|b| < |a|/2$ is smaller and thus works simultaneously for both cases. ∎